\documentclass[11pt]{article}
 \usepackage{graphicx}		
 \usepackage{multicol}
 \usepackage{amsmath, amsthm, amssymb}	
 \usepackage
 [
  pdfpagemode={UseNone},
  pdfstartview={FitH},
  colorlinks,
  urlcolor={blue},
  pdftitle={First occurrences of square-free gaps and an algorithm for their computation},
  pdfauthor={Louis Marmet},
  pdfsubject={square-free gaps},
  pdfkeywords={number theory, square-free},
  pdfdisplaydoctitle
 ]{hyperref}

\title{First occurrences of square-free gaps and an algorithm for their computation}
\author{Louis Marmet}
\date{November 14$^{th}$, 2007}

\begin{document}
\maketitle

\section{Abstract}
  This paper reports the results of a search for first occurrences of square-free gaps using an algorithm based on the sieve of Eratosthenes.  Using $Qgap(L)$ to denote the starting number of the first gap having exactly the length~$L$, the following values were found since August 1999: $Qgap(10)=262\,315\,467$, $Qgap(12)=47\,255\,689\,915$, $Qgap(13)=82\,462\,576\,220$, $Qgap(14)=1\,043\,460\,553\,364$, $Qgap(15)=79\,180\,770\,078\,548$, $Qgap(16)=$\\
$3\,215\,226\,335\,143\,218$, $Qgap(17)=23\,742\,453\,640\,900\,972$ and $Qgap(18)=125\,781\,000\,834\,058\,568$.  No gaps longer than $18$ were found up to $N=125\,870\,000\,000\,000\,000$.

\section{Introduction}
\label{sec:introduction}

\subsection{Square-free numbers}
\label{sec:squarefree}
  A number is said to be square-free if its prime decomposition contains no repeated factors.  For example, $30$ is square-free since its prime decomposition $2\times 3\times 5$ contains no repeated factors.  However, $18$ is not square-free since the factor $3$ appears twice in its prime decomposition $2\times 3\times 3$.

  The first few square-free numbers give the sequence: $1$, $2$, $3$, $5$, $6$, $7$, $10$, $11$, $13$, $14$, $15$, $17$, $19$, $21$, $22$, $23$, $26$, $29$, $30$, $31$, $33$, $34$, $35$, $37$, $38$, $39$, $41$, $42$, $43$, $46$, $47$, $51$, $53$, $55$, $57$, $58$, $59$, $61$, $62$, $65$, $66$, $67$, $69$, $70$, $71$, $73$, $74$, $77$, $78$, $79$, $82$, $83$, $85$, $86$, $87$, $89$, $91$, $93$, $94$, $95$, $97$, etc. (Sloane's A005117~\footnote{\href{http://oeis.org/A005117} {oeis.org/A005117}}.)

  More details on square-free numbers can be found at \href{http://mathworld.wolfram.com/Squarefree.html} {mathworld.wolfram.com/Squarefree.html}.  If a square-free number is used as the argument of the M\"obius function\footnote{\href{http://mathworld.wolfram.com/MoebiusFunction.html} {mathworld.wolfram.com/MoebiusFunction.html}}, a non-zero value ($+1$ or $-1$) is obtained.

\subsection{Gaps between square-free numbers}
\label{sec:gaps}
  A square-free gap is a series of $L$ consecutive numbers missing from the sequence of square-free numbers.  The first square-free gap in the sequence of square-free numbers starts at $N=4$ and has a length of one.  The next gap starts at $N=8$ and has a length $L=2$ (since $8$ and $9$ are non-square-free).  The following table lists the first few gaps and their lengths.

\begin{table}[ht]
 \caption{First few square-free gaps and their lengths}
 \centering
 \begin{tabular}{|l|c|c|c|c|c|c|c|c|c|c|c|c|c|c|c|c|}
  \hline
  Gap starts at $N$& 4& 8& 12& 16& 18& 20& 24& 27& 32& 36& 40& 44& 48& 52& 54& 56\\
  \hline
  Length of gap $L$& 1& 2& 1& 1& 1& 1& 2& 2& 1& 1& 1& 2& 3& 1& 1& 1\\
  \hline
 \end{tabular}
 \label{table:first_gaps}
\end{table}

  For any $L$, it can be shown that there exist infinitely many gaps of length greater than $L$ in the sequence of square-free numbers.  Longer lists of square-free gaps and recent results are available in Appendices B, C and D.  These gaps are series of consecutive squareful\footnote{\href{http://mathworld.wolfram.com/Squareful.html} {mathworld.wolfram.com/Squareful.html}} numbers (Sloane's A013929~\footnote{\href{http://oeis.org/A013929} {oeis.org/A013929}}).

  Note that the term ``squarefull'' sometimes denotes a positive integer $n$ such that if $p$ is a prime dividing $n$, then $p^2$ divides $n$.\footnote{M.~Filaseta, O.~Trifonov, ``The distribution of squarefull numbers in short intervals,'' Acta Arith. \textbf{67}, pp. 323-333, 1994.}  (See the selected Preprints of Michael Filaseta\footnote{\href{http://www.math.sc.edu/~filaseta/paperindex.html} {www.math.sc.edu/$\sim$filaseta/paperindex.html}}.)

  The smallest integer of the first gap having exactly the length $L$ is denoted here as $Qgap(L)$ ($Q$ for ``Quadratfrei''\footnote{\href{http://mathworld.wolfram.com/Quadratfrei.html} {mathworld.wolfram.com/Quadratfrei.html}}, or squarefree\footnote{\href{http://mathworld.wolfram.com/Squarefree.html} {mathworld.wolfram.com/Squarefree.html}}).  Thus $Qgap(1)=4$, $Qgap(2)=8$, $Qgap(3)=48$, etc.

\subsection{Upper limits for \texorpdfstring{$Qgap(16)$}{Qgap(16)} to \texorpdfstring{$Qgap(24)$}{Qgap(24)}}
\label{sec:upperlimits}
  Erick Bryce Wong\footnote{erick at sfu.ca.} found upper limits for $Qgap(L)$ for $L=16$ to $24$.  His idea was to find them by prescribing a repeated prime factors for each term and using the Chinese Remainder Theorem\footnote{\href{http://www.cut-the-knot.org/blue/chinese.shtml} {www.cut-the-knot.org/blue/chinese.shtml} and \href{http://mathworld.wolfram.com/ChineseRemainderTheorem.html} {mathworld.wolfram.com/ChineseRemainderTheorem.html}} to obtain a number.  More precisely, he prescribed all but five of the moduli and then tested the last moduli, up to the first $1000$ primes, to check if the number is squarefree. He tried this over millions of permutations.  His impressive results are:

\begin{displaymath}
  \begin{array}{lrl}
  & \text{Upper limit}& \text{Found on}\\ 
  Qgap(16)\leq& 46\,717\,595\,829\,767\,167\,& \text{Feb. 17$^{th}$ 2000}\\
  Qgap(17)\leq& 23\,742\,453\,640\,900\,972\,& \text{Feb. 21$^{st}$ 2000}\\
  Qgap(18)\leq& 125\,781\,000\,834\,058\,568& \text{June 26$^{th}$ 2000}\\
  Qgap(19)\leq& 31\,310\,794\,237\,768\,728\,712& \text{July 18$^{th}$ 2000}\\
  Qgap(20)\leq& 148\,372\,453\,443\,663\,297\,638\,331& \text{July 10$^{th}$ 2000}\\
  Qgap(21)\leq& 321\,362\,101\,382\,225\,854\,472& \text{Feb. 17$^{th}$ 2000}\\
  Qgap(22)\leq& 213\,922\,449\,434\,979\,698\,424\,416& \text{Aug. 4$^{th}$ 2000}\\
  Qgap(23)\leq& 687\,445\,369\,966\,391\,012\,821\,156\,868& \text{July 18$^{th}$ 2000}\\
  Qgap(24)\leq& 28\,548\,715\,276\,566\,524\,078\,226\,797\,585\,011& \text{Sept. 4$^{th}$ 2000}
  \end{array}
\end{displaymath}

\subsection{First occurrence of a gap of length \texorpdfstring{$L$}{L}}
\label{sec:firstoccurence}
  $Qgap(L)$ is listed in the following table for $L$ up to $18$.  The third column gives the prime factors that are repeated for each number in the gap.  The values of $Qgap(L<10)$ and $Qgap(11)$ have been confirmed by different sources. (See for example, ``Sloane's On-Line Encyclopedia of Integer Sequences''\footnote{\href{http://oeis.org/} {oeis.org/}}, sequence A045882~\footnote{ \href{http://oeis.org/A045882} {oeis.org/A045882}.}.)  The sequence $Qgap(L)$ is listed under A051681~\footnote{\href{http://oeis.org/A051681} {oeis.org/A051681}} in the Encyclopedia of Integer sequences.

\begin{displaymath}
  \begin{array}{rrcr}
  L& Qgap(L)& \text{Repeated prime factors of each number in the gap}& \text{Gap reported by}\\			
  1& 4& 2& \text{E. Friedman}\\
  2& 8& 2,3& \text{E. Friedman}\\
  3& 48& 2,7,5& \text{E. Friedman}\\
  4& 242& 11,3,2,7& \text{E. Friedman}\\
  5& 844& 2,13,3,11,2& \text{E. Friedman}\\
  6& 22\,020& 2,19,11,3,2,5& \text{E. Friedman}\\
  7& 217\,070& 7,3,2,113,11,5,2& \text{E. Friedman}\\
  8& 1\,092\,747& 19,2,7,5,11,2,3,13& \text{E. Friedman}\\
  9& 8\,870\,024& 2,5,11,29,2,7,31,3,2& \text{P. De Geest}\\
  10& 262\,315\,467& 3,2,29,2957,79,2,7,17,5,2\times 3& \text{D. Bernier}\\
  11& 221\,167\,422& 3,31,2,5,37,13,2,7,11,3,2& \text{P. De Geest}\\
  12& 47\,255\,689\,915& 7,2,3,103,43,2,29,17,13,2,5,3& \text{L. Marmet}\\
  13& 82\,462\,576\,220& 2,3,13,23,2,5,17,41,2,19,3,7,2& \text{L. Marmet}\\
  14& 1\,043\,460\,553\,364& 2,3,7,19,2,13\times 59,67,43,2,181,3,5,2,11& \text{L. Marmet}\\
  15& 79\,180\,770\,078\,548& 2,3,5,29,2,13,17,53,2,19,3,41,2,31,67& \text{L. Marmet}\\
  16& 3\,215\,226\,335\,143\,218& 11,23,2,3,269,53,2,5,17,163,2,101,3,19,2,137& \text{Z. McGregor-Dorsey}\\
  17& 23\,742\,453\,640\,900\,972& 2,11\times 23,127,5,2,3,53,37,2,7,13,17,2,19,3,29,2& \text{E. Wong}\\
  18& 125\,781\,000\,834\,058\,568& 2,3,37,31,2,19,29,5,2,7\times 23,3,139,2,11,17,13,2,199& \text{L. Marmet}
  \end{array}
\end{displaymath}

  The first gaps reported in this work were found on the following dates.

\begin{displaymath}
  \begin{array}{lrrl}
  & \text{First gap}& \text{Found on }& \text{ Found by}\\
  Qgap(10)=& 262\,315\,467& \text{August 1999}& \text{D. Bernier},\\
  Qgap(12)=& 47\,255\,689\,915& \text{October 19$^{th}$ 1999}& \text{L. Marmet},\\
  Qgap(13)=& 82\,462\,576\,220& \text{October 20$^{th}$ 1999}& \text{L. Marmet},\\
  Qgap(14)=& 1\,043\,460\,553\,364& \text{October 25$^{th}$ 1999}& \text{L. Marmet},\\
  Qgap(15)=& 79\,180\,770\,078\,548& \text{November 29$^{th}$ 1999}& \text{L. Marmet},\\
  Qgap(16)=& 3\,215\,226\,335\,143\,218& \text{July 22$^{nd}$ 2000}& \text{Z. McGregor-Dorsey \emph{et al.}},\\
  Qgap(17)=& 23\,742\,453\,640\,900\,972& \text{July 8$^{th}$ 2001}& \text{E. Wong \emph{et al.}},\\
  Qgap(18)=& 125\,781\,000\,834\,058\,568& \text{September 9$^{th}$ 2005}& \text{L. Marmet \emph{et al.}}
  \end{array}
\end{displaymath}

\subsubsection{Basic algorithm: the sieve of Eratosthenes}
\label{sec:basicalgorithm}
  The square-free gaps can be calculated by finding consecutive numbers that are not square-free.  A simple method to show that $N$ is not square­free is to find a prime factor of $N$ whose square divides $N$.  By trying every prime up to the square root of $N$, one can establish whether $N$ is square-free or not.  However, this is a very inefficient way to test billions of numbers.

  A faster algorithm is used by ``Mathematica'' to determine if a number is square-free.  The method is quite interesting\footnote{\href{http://reference.wolfram.com/mathematica/ref/SquareFreeQ.html} {reference.wolfram.com/mathematica/ref/SquareFreeQ.html}}.

  However, to determine which of many consecutive numbers are square-free, an algorithm based on to the sieve of Eratosthenes\footnote{\href{http://mathworld.wolfram.com/SieveofEratosthenes.html} {mathworld.wolfram.com/SieveofEratosthenes.html}} is much faster.  It uses a list of numbers from which each composite number is removed.  Once the process is finished, only the prime numbers are in the list.

  To find the square-free numbers using a sieve, a similar technique is used but the algorithm eliminates numbers that are not square-free.  Starting with a list of integers, first cross out the multiples of $4$:
 
\texttt{~1 ~2 ~3 ~X ~5 ~6 ~7 ~X ~9 10 11 ~X 13 14 15 ~X 17 18 19 ~X 21 22 23 ~X 25 26 ...}

then the multiples of $9$, $25$, etc., up to the last number in the list:

\texttt{~1 ~2 ~3 ~X ~5 ~6 ~7 ~X ~X 10 11 ~X 13 14 15 ~X 17 ~X 19 ~X 21 22 23 ~X ~X 26~...}

The remaining numbers are square-free numbers; the gaps are indicated by the series of consecutive ``X''.

\subsubsection{Improvements of the algorithm}
\label{sec:improvements}
  The following improvements were implemented in a computer program and are presented in the same order they were added to the program.

\subsubsection{Improvement I}
\label{sec:improvement1}
\begin{center}
  ``Lists of squared-primes and the next non-square-free number use less memory.''
\end{center}

  To implement this algorithm on a computer, it is not necessary to keep the entire list of integers in memory.  An improvement of the algorithm uses instead two shorter arrays to calculate the next non-square-free number after $N$:
\begin{itemize}
	\item the first array, called \texttt{p2}, gives the squares of the prime numbers up to the largest number to be tested $N_{max}$,
	\item the second array, called \texttt{nsqf}, gives for each \texttt{p2[i]} the next non-square-free number, that is, the smallest number larger than $N$ that is a multiple of \texttt{p2[i]}.  This array can easily be calculated using modulo arithmetic.
\end{itemize}

  These arrays will have approximately $2\sqrt{N_{max}}/\ln{N_{max}}$ elements\footnote{\href{http://www.utm.edu/research/primes/howmany.shtml} {www.utm.edu/research/primes/howmany.shtml}}.

  To find square-free gaps, one finds sequences of non-square-free numbers.  The following example shows the arrays used to find gaps starting from $N=20$~\footnote{The notation used in the programming language C is used here, where the index of an array starts at 0.}.

\begin{displaymath}
  \begin{array}{lcrrrrrrrrrr}
  \text{Index}& \texttt{i}& 0& 1& 2& 3& 4& 5& 6& 7& 8& ...\\
  \text{Squared prime}& \texttt{p2[i]}& 4& 9& 25& 49& 121& 169& 289& 361& 529& ...\\
  \text{Next non-square-free}& \texttt{nsqf[i]}& 24& 27& 25& 49& 121& 169& 289& 361& 529& ...
  \end{array}
\end{displaymath}

  Using this table, it is easy to find the next non-square-free number: it is the smallest number in the array \texttt{nsqf}, that is, $24$.  We set $N=24$ and recalculate the array.  This is easy since the only needed operation is to add the corresponding squared-prime to the multiple: $24+4=28$.  The following table is obtained:

\begin{displaymath}
  \begin{array}{lcrrrrrrrrrr}
  \text{Index}& \texttt{i}& 0& 1& 2& 3& 4& 5& 6& 7& 8& ...\\
  \text{Squared prime}& \texttt{p2[i]}& 4& 9& 25& 49& 121& 169& 289& 361& 529& ...\\
  \text{Next non-square-free}& \texttt{nsqf[i]}& 28& 27& 25& 49& 121& 169& 289& 361& 529& ...
  \end{array}
\end{displaymath}

  Again, the next non-square-free number is the smallest \texttt{nsqf[i]}.  By repeating this procedure, $N$ takes the values of all the non-square-free numbers.  It is advantageous to sort \texttt{nsqf} in increasing order at each step.  This way, the smallest is always \texttt{nsqf[0]}.  We set $N=25$ and recalculate the array to obtain:

\begin{displaymath}
  \begin{array}{lcrrrrrrrrrr}
  \text{Index}& \texttt{i}& 0& 1& 2& 3& 4& 5& 6& 7& 8& ...\\
  \text{Squared prime}& \texttt{p2[i]}& 9& 4& 49& 25& 121& 169& 289& 361& 529& ...\\
  \text{Next non-square-free}& \texttt{nsqf[i]}& 27& 28& 49& 50& 121& 169& 289& 361& 529& ...
  \end{array}
\end{displaymath}

  The order of the array \texttt{p2} has also been changed so that each number \texttt{p2[i]} always corresponds to its multiple \texttt{nsqf[i]}.  Repeating the procedure will generate the non-square-free numbers $N=27$, $28$, $32$, $36$, $40$, etc.  Note that special care has to be taken when some numbers in \texttt{nsqf} are equal - each of these has to be increased by the value of its corresponding \texttt{p2[i]}.

  The sort is relatively efficient since after \texttt{nsqf[0]} is given its new value, \texttt{nsqf[1]}, \texttt{nsqf[2]} and the following elements are still in increasing order.  The new value is moved up the array until its proper place is found.

\subsubsection{Improvement II}
\label{sec:improvement2}
\begin{center}
  ``Many non-square-free numbers can be skipped.''
\end{center}

  If gaps of a given length $L_{min}$ or more are searched, some \texttt{nsqf[i]} can be skipped.  To show this, one finds first the minimum number of squared-primes \texttt{NP2$_{min}$}~\footnote{\href{http://oeis.org/A107079} {http://oeis.org/A107079}} required in a gap of length $L$:
 
\begin{displaymath}
  \begin{array}{lcrrrrrrrrrrrrrrrrr}
  \text{Gap length }L& 1& 2& 3& 4& 5& 6& 7& 8& 9& 10& 11& 12& 13& 14& 15& 16& 17& ...\\
  \texttt{NP2$_{min}$[$L$]}& 1& 2& 3& 4& 4& 5& 6& 7& 7& 7&  8&  9&  9& 10& 11& 12& 12& ...
  \end{array}
\end{displaymath}

  For example if we choose $L=7$, \texttt{NP2$_{min}$[$L$]$=6$} prime factors are required for the gap starting at $N=217070$ (in this case, the prime factors are $2$, $3$, $5$, $7$, $11$ and $113$).

  Continuing with the example given above, we now specifically search gaps with length $L_{min}=7$ or longer.  There is no gap of length $L_{min}=7$ in the interval starting at \texttt{nsqf[0]} and ending at \texttt{nsqf[5]}, if\\
\texttt{nsqf[5]$>$nsqf[0]$+7$}.  (In general, there is no gap of length $L_{min}$ in the interval starting at \texttt{nsqf[0]} and ending at \texttt{nsqf[NP2$_{min}$[$L_{min}-1$]]}, if \texttt{nsqf[NP2$_{min}$[$L_{min}-1$]]$>$nsqf[0]$+L_{min}$}.)

  With $N=40$, we have:
 
\begin{displaymath}
  \begin{array}{lcrrrrrrrrrr}
  \text{Index}& \texttt{i}& 0& 1& 2& 3& 4& 5& 6& 7& 8& ...\\
  \text{Squared prime}& \texttt{p2[i]}& 4& 9& 49& 25& 121& 169& 289& 361& 529& ...\\
  \text{Next non-square-free}& \texttt{nsqf[i]}& 44& 45& 49& 50& 121& 169& 289& 361& 529& ...
  \end{array}
\end{displaymath}

  Since \texttt{nsqf[5] $=169>$ nsqf[0]$+7=51$}, there is no gap of length $7$ in the interval starting at $44$ and ending at $169$.  We can therefore safely set \texttt{$N=$nsqf[5]$-L_{min}+1=163$} and recalculate the following table:

\begin{displaymath}
  \begin{array}{lcrrrrrrrrrr}
  \text{Index}& \texttt{i}& 0& 1& 2& 3& 4& 5& 6& 7& 8& ...\\
  \text{Squared prime}& \texttt{p2[i]}& 4& 169& 9& 25& 49& 121& 289& 361& 529& ...\\
  \text{Next non-square-free}& \texttt{nsqf[i]}& 164& 169& 171& 175& 196& 242& 289& 361& 529& ...
  \end{array}
\end{displaymath}

  This cuts down on the number of non-square-free that have to be tested and the speed of the calculation is increased.  A factor five in speed was obtained when this was implemented in the program which was used to find $Qgap(14)$ and $Qgap(15)$.

\subsubsection{Improvement III}
\label{sec:improvement3}
\begin{center}
  ``The smallest squared-primes are not needed to calculate the sieve.''
\end{center}

  This variation on Improvement II was suggested by Joseph Wetherell.  It turns out to be more efficient when it is combined with Improvements IV and V.  The trick is to consider the smallest squared-primes separately from the large squared-primes.  Most of the time in the algorithm is spent on the process of taking the smallest elements off of \texttt{nsqf} and sorting them back into the array.  If one can reduce the number of non-square-free numbers which are tested, the speed of the algorithm will improve.

  This is actually possible since the smallest squared-primes are not needed to calculate the sieve.  If we have a gap of, say, $L=L_{min}=7$ non-square-free numbers, then at least \texttt{NP2$_{min}$[7]$=6$} different primes are found in the gap.  This means that the smallest \texttt{$k1=$NP2$_{min}$[$L_{min}$]$-1=5$} squared-primes can be left out of the calculation.  If we search for the 6$^{th}$ squared-prime, every gap of length $L=7$ (or more) will be found.

  This method therefore separates the small squared-primes from the large ones, creating a base with $k1$ elements.  The table for $N=163$ would now look like this:

\begin{displaymath}
  \begin{array}{lcrrrrrrrrrr}
  \text{Index}& \texttt{i}	& 0& 1& 2& 3& 4& k1& 6& 7& 8& ...\\
  \text{Squared prime}& \texttt{p2[i]}& 4& 9& 25& 49& 121& 169& 289& 361& 529& ...\\
  \text{Next non-square-free}& \texttt{nsqf[i]}& --& --& --& --& --& 169& 289& 361& 529& ...
  \end{array}
\end{displaymath}

with the base shown as $--$ for the values of \texttt{nsqf[i]}.  From the table, one sees that one must test for a gap having $L_{min}=7$ around $N=169$.  If none is found, then \texttt{nsqf[$k1$]} is increased by $169$ and sorted back into the array of large squared-primes to get:
 
\begin{displaymath}
  \begin{array}{lcrrrrrrrrrr}
  \text{Index}& \texttt{i}& 0& 1& 2& 3& 4& k1& 6& 7& 8& ...\\
  \text{Squared prime}& \texttt{p2[i]}& 4& 9& 25& 49& 121& 289& 169& 361& 529& ...\\
  \text{Next non-square-free}& \texttt{nsqf[i]}& --& --& --& --& --& 289& 338& 361& 529& ...
  \end{array}
\end{displaymath}

  The sort is faster since it is only done on the large squared-primes.  The process is then continued at $N=289$.

\subsubsection{Improvement IV}
\label{sec:improvement4}
\begin{center}
  ``The values of the small modulos can be computed ahead of time.''
\end{center}

  To test if there is a gap of $L_{min}$ around $N$, one must still know about multiples of the $k1$ small squared-primes near $N$.  As suggested by Joseph Wetherell, this can be done by trial division; even if trial division is slow, it is faster than resorting the base array.  One can also optimize the trial divisions, because the trial divisions for, say, $N + 1$ and $N + 2$ are related to each other.  For each small prime $p$, compute $N\%$\texttt{p2} and store it in a list \texttt{mod} (the ``\%'' symbol is the modulo function in the C language).  Now to see if \texttt{p2} divides $N+1$, we just test if this stored value is $-1$ (mod \texttt{p2}).  To see if \texttt{p2} divides $N+2$, we test if this stored value is $-2$ (mod \texttt{p2}).  (We also precompute the value of $-1$ (mod \texttt{p2}), $-2$ (mod \texttt{p2}), etc., for the small set of values which we will possibly need to test.)  Note that it is also necessary to test $N-1$, $N-2$, etc.  For $N=289$, we have the following arrays:
 
\begin{displaymath}
  \begin{array}{lcrrrrrrrrrr}
  \text{Index}& \texttt{i}& 0& 1& 2& 3& 4& k1& 6& 7& 8& ...\\
  \text{Squared prime}& \texttt{p2[i]}& 4& 9& 25& 49& 121& 289& 169& 361& 529& ...\\
  \text{Next non-square-free}& \texttt{nsqf[i]}& --& --& --& --& --& 289& 338& 361& 529& ...\\
  N \% \texttt{p2[i]}& \texttt{mod[i]}& 1& 1& 14& 44& 47& & & & & 
  \end{array}
\end{displaymath}

  We see that \texttt{p2[0]} and \texttt{p2[1]} divide $N-1=288$, but this is the only other non-square-free number for this gap of length $2$.  We can therefore increase \texttt{nsqf[$k1$]} by $289$, sort it back into the array (reordering \texttt{p2} accordingly) and continue with $N=338$:
 
\begin{displaymath}
  \begin{array}{lcrrrrrrrrrr}
  \text{Index}& \texttt{i}& 0& 1& 2& 3& 4& k1& 6& 7& 8& ...\\
  \text{Squared prime}& \texttt{p2[i]}& 4& 9& 25& 49& 121& 169& 361& 529& 289& ...\\
  \text{Next non-square-free}& \texttt{nsqf[i]}& --& --& --& --& --& 338& 361& 529& 578& ...\\
  N \% \texttt{p2[i]}& \texttt{mod[i]}& 2& 5& 13& 44& 96& & & & & 
  \end{array}
\end{displaymath}

\subsubsection{Improvement V}
\label{sec:improvement5}
\begin{center}
  ``Two large primes-squared that are too far cannot result in a gap.''
\end{center}

  If we include \texttt{p2[4]} in the array of large squared-primes, (so there are only \texttt{$k2=$NP2$_{min}$[$L_{min}$]$-2=4$} elements in the base), then we know that a number \texttt{$N=$ nsqf[$k2$]} we are testing can be part of a gap only if the next number in \texttt{nsqf} is close to $N$, that is, if \texttt{nsqf[k2+1]} is not larger than $N+L_{min}-1$.  Based on this suggestion by Joseph Wetherell, the arrays become with $N=338$:

\begin{displaymath}
  \begin{array}{lcrrrrrrrrrr}
  \text{Index}& \texttt{i}& 0& 1& 2& 3& k2& 5& 6& 7& 8& ...\\
  \text{Squared prime}& \texttt{p2[i]}& 4& 9& 25& 49& 169& 361& 121& 529& 289& ...\\
  \text{Next non-square-free}& \texttt{nsqf[i]}& --& --& --& --& 338& 361& 363& 529& 578& ...\\
  N \% \texttt{p2[i]}& \texttt{mod[i]}& --& --& --& --& & & & & & 
  \end{array}
\end{displaymath}

  Since \texttt{nsqf[$k2+1$]$=361$} is larger than $N+L_{min}-1=344$, we can skip to $N=361$.  Since $N$ does not pass the closeness test, we do not have to do any computations with the base, saving us a lot of time.

\subsubsection{Improvement VI}
\label{sec:improvement6}
\begin{center}
  ``A chained list is faster for the sort.''
\end{center}

  A chained list can be built with a set of numbers that specify an order for the elements of an array, as suggested by Joseph Wetherell.  If we use a chained list represented by the array called \texttt{next} such that \texttt{nsqf[next[i]]$\geq$nsqf[i]}, we can sort the array \texttt{nsqf} without moving any data within the arrays \texttt{nsqf} or \texttt{p2}.  For a reason that will become obvious later, we choose to have \texttt{p2} sorted in increasing order.  With $N=361$, the arrays would be:
 
\begin{displaymath}
  \begin{array}{lcrrrrrrrrrr}
  \text{Index}& \texttt{i}& 0& 1& 2& 3& k2& 5& 6& 7& 8& ...\\
  \text{Chained list}& \texttt{next[i]}& --& --& --& k2& 5& 8& 9& \underline{\textbf{4}}& 6& ...\\
  \text{Squared prime}& \texttt{p2[i]}& 4& 9& 25& 49& 121& 169& 289& 361& 529& ...\\
  \text{Next non-square-free}& \texttt{nsqf[i]}& --& --& --& --& 363& 507& 578& 361& 529& ...\\
  N \% \texttt{p2[i]}& \texttt{mod[i]}& 1& 1& 11& 18& & & & & & 
  \end{array}
\end{displaymath}
	
  We use an additional variable, ``\texttt{head}'', which points to the smallest item in the array \texttt{nsqf}.  For the table above, we have \texttt{head $=7$} (\texttt{next[head]} is highlighted).  Since we only need to change two values in the array \texttt{next} to perform a sort, this method is faster.  Searching through the array \texttt{nsqf} now consists of going through the data in the following order:

\texttt{for (i=head; tempnsqf>=nsqf[next[i]]; i=next[i]) { ... }}

  There is another advantage to this method: since the array \texttt{p2} is sorted in increasing order, we know that if we have to sort item \texttt{nsqf[i=head]}, then \texttt{nsqf[i-1]$\leq$nsqf[i]$+$p2[i]}.  This means we can jump immediately to \texttt{i=head-1} and start searching from there.  In general, this cuts the search in half!  Searching through the list now consists of:

\texttt{for (i=head-1; tempnsqf>=nsqf[i2=next[i]]; i=i2) { ... }}

  We continue the example with the above table.  Since \texttt{nsqf[next[head]]$=363$} is not larger than $N+L_{min}-1=367$, we calculated the modulos.  Clearly, there is only a gap of length $L=2$ starting at $360$.  We therefore set $N=363$ and sort \texttt{tempnsqf$=$nsqf[head]$=$nsqf[head]$+$p2[head]$=722$}.  We start with \texttt{i=6} and the sort requires only one comparison!  The following arrays are obtained:

\begin{displaymath}
  \begin{array}{lcrrrrrrrrrr}
  \text{Index}& i& 0& 1& 2& 3& k2& 5& 6& 7& 8& ...\\
  \text{Chained list}& \texttt{next[i]}& --& --& --& k2& \underline{\textbf{5}}& 8& 7& 9& 6& ...\\
  \text{Squared prime}& \texttt{p2[i]}& 4& 9& 25& 49& 121& 169& 289& 361& 529& ...\\
  \text{Next non-square-free}& \texttt{nsqf[i]}& --& --& --& --& 363& 507& 578& 722& 529& ...\\
  N \% \texttt{p2[i]}& \texttt{mod[i]}& --& --& --& --& & & & & & 
  \end{array}
\end{displaymath}

with \texttt{head} now equal to $4$.

  A factor three in speed was obtained when this algorithm was implemented in the program.

\subsubsection{Improvement VII}
\label{sec:improvement7}
\begin{center}
  ``Look for the largest spacing between two of three large squared-primes.''
\end{center}

  Instead of having \texttt{k2} elements in the base and look for two large squared-primes that are not too far apart, we can use a base with \texttt{k3$=$NP2$_{min}$[$L_{min}$]$-3=3$} squared-primes, but look to see if \texttt{nsqf[head]} and\\
\texttt{nsqf[next[next[head]]]} are not more than $L_{min}$ apart.
 
\begin{displaymath}
  \begin{array}{lcrrrrrrrrrr}
  \text{Index}& i& 0& 1& 2& k3& 4& 5& 6& 7& 8& ...\\
  \text{Chained list}& \texttt{next[i]}& --& --& k3& 5& \underline{\textbf{3}}& 8& 7& 9& 6& ...\\
  \text{Squared prime}& \texttt{p2[i]}& 4& 9& 25& 49& 121& 169& 289& 361& 529& ...\\
  \text{Next non-square-free}& \texttt{nsqf[i]}& --& --& --& 392& 363& 507& 578& 722& 529& ...\\
  N \% \texttt{p2[i]}& \texttt{mod[i]}& --& --& --& & & & & & & 
  \end{array}
\end{displaymath}

  In this example, we look at the difference between \texttt{nsqf[head]$=363$} and \texttt{nsqf[next[next[head]]]$=507$}.  Since the difference between the two values is larger than $L_{min}-1$, there is no gap of length $L_{min}$ or longer.

  This improvement gives the program a 37\% speed increase with $L_{min}=14$.

  The program becomes slower if four or more large squared-primes are considered (this was confirmed in tests for $L>13$, $N=10^{14}$ to $10^{14}+10^9$).

\subsubsection{Evaluation of order of algorithm}
\label{sec:order}
  This evaluation applies to the first improvement of the algorithm.  It was found empirically that for a given value of $L_{min}$, the other improvements increased the speed of the calculation by a constant factor.

  To calculate if $N$ is a square-free number, the algorithm takes advantage of the known remainders for $N-1$.  Each time $N$ is tested, a new non-square-free number is calculated using \texttt{nsqf[0]=nsqf[0]+p2[0];}.  This new value has to be moved up the list to keep \texttt{nsqf} in increasing order.  It is that operation that requires most of the computation time.  To evaluate the speed of the algorithm, it is necessary to find the average number of moves $m$ that will be required to bring the new value \texttt{nsqf[0]} to its correct position in the list, above the number \texttt{nsqf[$m$]}.  This is done by evaluating, for every \texttt{i}, the probability that \texttt{nsqf[0]$>$nsqf[i]}, and then summing over \texttt{i}.

  First, consider the case when $p2[0]=4$ which occurs with a probability of $1/4$.  Since \texttt{nsqf[0]} has been increased by $4$, there is a probability of $4/9$ that it will have to be moved above \texttt{nsqf[j]} (if \texttt{p2[j]$=9$}).  There is an additional probability of $4/25$ that \texttt{nsqf[0]} will have to be moved above \texttt{nsqf[k]} (if \texttt{p2[k]$=25$}), etc.  We therefore get the average number of steps required to place the new \texttt{nsqf[0]} to its correct position:

$$S(1)=1/4\times (4/9+4/25+4/49+4/121+...)=\sum_{i=2}\frac{1}{p^2(i)}$$

where $p(i)$ is the ith prime number ($p(1)=2$, $p(2)=3$, $p(3)=5$, etc.) and $p^2(i)=p(i)\times p(i)$.

  In the case when \texttt{p2[0]$=9$} (which occurs with a probability of $1/9$), one move is always necessary to bring \texttt{nsqf[0]} above \texttt{nsqf[j]} (if \texttt{p2[j]$=4$}).  There is a probability of $9/25$ that \texttt{nsqf[0]} will have to be moved above \texttt{nsqf[k]} (if p2[k]=25), etc.  We therefore get:

$$S(2)=1/9\times (1+9/25+9/49+9/121+...)=\frac{1}{p^2(2)}+\sum_{i=3}\frac{1}{p^2(i)}$$

  In general, when \texttt{p2[0]$=p^2(m)$}, the average number of moves is:

$$S(m)=\frac{m-1}{p^2(m)} + \sum_{i=m+1}\frac{1}{p^2(i)}$$

  The sum over all the $S(m)$ gives the average number of moves required to place \texttt{nsqf[0]} to its correct position in the list:

$$\sum_{m=1} S(m)=2\times \sum_{i=2}\frac{i-1}{p^2(i)}$$

This series converges, as determined with a convergence test\footnote{``Suppose that $f(x)$ is a positive decreasing function and that $$lim_{k\rightarrow\infty}\frac{e^kf(e^k)}{f(k)}=q$$ for natural $k$.  If $q<1$, the series $\sum_{k=1}^\infty f(k)$ converges. If $q>1$, this series diverges. (Ermakov)''  Equation 0.224, ``Table of Integrals, Series, and Products,'' Gradshteyn and Ryzhik (Academic Press, Inc., p. 5)}.  It converges very slowly to approximately $1.30...$  Therefore, given the remainders for $N-1$, the number of operations required to find out if $N$ is square-free is independent of the value of $N$.

\subsection{Acknowledgements}
\label{sec:acknowledgements}
  Contributed to find $Qgap(10)$: David Bernier.

  Contributed to find $Qgap(16)$: Zach McGregor-Dorsey, Louis Marmet, Joe Wetherell, Gunnard Engebreth, D. Bernier, Erick Wong, Alan Simpson and Nicolas Marmet.

  Contributed to find $Qgap(17)$: E. Wong, Z. McGregor-Dorsey, L. Marmet, Jean-Pierre Bernier, D. Bernier, Nancy Robertson, N. Marmet, Charles Ward and G. Engebreth.

  Contributed to find $Qgap(18)$: D. Bernier, L. Marmet, E. Wong, J. Wetherell, Z. McGregor-Dorsey, G. Engebreth, A. Simpson, N. Marmet, N. Robertson, J.-P. Bernier, C.R. Ward, Bruno Le Tual and Horand Gassmann.

  This project started from an idea that was initially suggested to me by David Bernier.

\subsection{Related web pages and references}
\label{sec:links}

M.~Filaseta, O.~Trifonov, ``The distribution of squarefull numbers in short intervals,'' Acta Arith., 67 (1994), 323-333.

M.~Filaseta, ``On the distribution of gaps between squarefree numbers,'' Mathematika, 40 (1993), 88-101.

M.~Filaseta, O.~Trifonov, ``On gaps between squarefree numbers II,'' Journal of the London Math. Soc. (2), 45 (1992), 215-221.

C. Rivera, ``The Prime Puzzles and Problems Connection,'' \href{http://www.primepuzzles.net/problems/prob_028.htm} {www.primepuzzles.net/problems/prob\_028.htm} (not updated since Nov. 1999).

B.~de~Weger, C.E.~van~de~Woestijne, ``On the powerfree parts of consecutive integers,'' Acta Arithmetica 90 (1999), 387-395, \href{http://www.win.tue.nl/~bdeweger/onderzoek.html} {http://www.win.tue.nl/$\sim$bdeweger/onderzoek.html}.

A translation of this paper is available in Belorussian at\\
\href{http://www.webhostinghub.com/support/by/edu/index-marmet-be} {www.webhostinghub.com/support/by/edu/index-marmet-be}.

\vskip 0.9in

This article was first published at 
\href{http://www.marmet.org/louis/sqfgap/index.html} {www.marmet.org/louis/sqfgap/index.html}.

Copyright $\copyright $1999  \href{http://www.marmet.org/louis/index.html} {Louis Marmet}

\clearpage

\section*{Appendix A}

  The following graph shows an estimation of how large we can expect $Qgap(L)$ to be.
 
\begin{figure}[ht]
  \centerline{\includegraphics[width=5.158in]{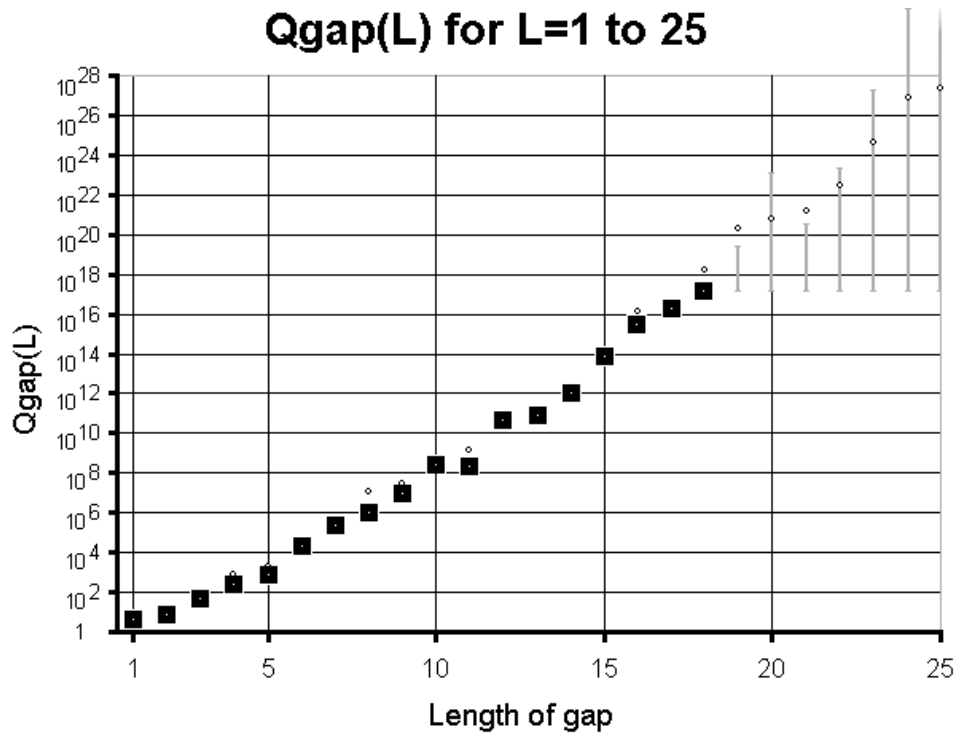}}
  \caption{Known values of $Qgap(L)$ (squares) and estimated values (empty circles).}
  \label{fig:firstsqfgaps}
\end{figure}

  The empirical estimation, based on the calculated values for $L<16$, uses an approximation of the probability of obtaining the minimum number of primes required to produce a gap with length $L$.  The upper limits for $Qgap(L>16)$ were obtained by E. Wong.  The values of $Qgap(L)$ lie within the ranges indicated by the gray lines.

\clearpage

\section*{Appendix B}

  Square-free gaps and their length $\geq 14$, up to $5\,000\,000\,000\,000\,000$

\setlength{\columnseprule}{0.1pt}
\begin{multicols}{6}
{
\tiny
\begin{flushright}

  \texttt{ 1043460553364: 14}\\
  \texttt{ 4086996180820: 14}\\
  \texttt{ 5795458534275: 14}\\
  \texttt{ 5907695879319: 14}\\
  \texttt{ 22021168374368: 14}\\
  \texttt{ 24908160544315: 14}\\
  \texttt{ 25488235201119: 14}\\
  \texttt{ 26176980660775: 14}\\
  \texttt{ 27075246766243: 14}\\
  \texttt{ 28721610171943: 14}\\
  \texttt{ 28805392031443: 14}\\
  \texttt{ 28880901889539: 14}\\
  \texttt{ 29785893889868: 14}\\
  \texttt{ 32037597678248: 14}\\
  \texttt{ 36123614624743: 14}\\
  \texttt{ 36442589727571: 14}\\
  \texttt{ 37137019337619: 14}\\
  \texttt{ 38606457834267: 14}\\
  \texttt{ 40811567562764: 14}\\
  \texttt{ 46456561416219: 14}

  \texttt{ 48975804002312: 14}\\
  \texttt{ 49245216926812: 14}\\
  \texttt{ 52333043400944: 14}\\
  \texttt{ 52363499321167: 14}\\
  \texttt{ 52765476883143: 14}\\
  \texttt{ 54814996088144: 14}\\
  \texttt{ 55548673575423: 14}\\
  \texttt{ 55865746434715: 14}\\
  \texttt{ 58852723694919: 14}\\
  \texttt{ 59068995605116: 14}\\
  \texttt{ 61747840706667: 14}\\
  \texttt{ 62149445184671: 14}\\
  \texttt{ 62618166402775: 14}\\
  \texttt{ 63330512887168: 14}\\
  \texttt{ 66917476836763: 14}\\
  \texttt{ 70458333341416: 14}\\
  \texttt{ 71559169190475: 14}\\
  \texttt{ 76319977090472: 14}\\
  \texttt{ 76381235551744: 14}\\
  \texttt{ 76622541563267: 14}

  \texttt{ 79180770078548: 15}\\
  \texttt{ 79967232391912: 14}\\
  \texttt{ 80970974553112: 14}\\
  \texttt{ 83696151329872: 14}\\
  \texttt{ 85867066913947: 14}\\
  \texttt{ 88490045527468: 14}\\
  \texttt{ 89141750469412: 14}\\
  \texttt{ 89872191351339: 14}\\
  \texttt{ 93553802716516: 14}\\
  \texttt{ 96115195097468: 14}\\
  \texttt{ 99169976343220: 14}\\
  \texttt{ 99321073884916: 14}\\
  \texttt{ 100626734394763: 14}\\
  \texttt{ 102947196430240: 14}\\
  \texttt{ 104575423448671: 14}\\
  \texttt{ 104912119695975: 14}\\
  \texttt{ 109163092681520: 14}\\
  \texttt{ 111061110643215: 14}\\
  \texttt{ 115660022682524: 14}\\
  \texttt{ 124706466248019: 14}

  \texttt{ 128908232869112: 14}\\
  \texttt{ 130275478658815: 14}\\
  \texttt{ 132787974312015: 14}\\
  \texttt{ 134933696029672: 14}\\
  \texttt{ 138101024911543: 14}\\
  \texttt{ 141834821831648: 14}\\
  \texttt{ 142482797711648: 14}\\
  \texttt{ 142540924809939: 14}\\
  \texttt{ 142650781820719: 14}\\
  \texttt{ 145065154350544: 15}\\
  \texttt{ 147261235475167: 14}\\
  \texttt{ 147913509128020: 14}\\
  \texttt{ 150621269480115: 14}\\
  \texttt{ 151135375176719: 14}\\
  \texttt{ 152547834577472: 14}\\
  \texttt{ 152579171596912: 14}\\
  \texttt{ 158762391708320: 14}\\
  \texttt{ 160033416707967: 14}\\
  \texttt{ 160507172057344: 14}\\
  \texttt{ 160797901723063: 14}

  \texttt{ 166273803771243: 14}\\
  \texttt{ 167653761845264: 14}\\
  \texttt{ 170729320777243: 14}\\
  \texttt{ 170879863247044: 14}\\
  \texttt{ 171662978357224: 14}\\
  \texttt{ 174001041813068: 14}\\
  \texttt{ 174314615292620: 14}\\
  \texttt{ 175041733150520: 14}\\
  \texttt{ 176802463160648: 14}\\
  \texttt{ 179205337927312: 14}\\
  \texttt{ 179928574705364: 14}\\
  \texttt{ 181431837731743: 14}\\
  \texttt{ 181521770360768: 14}\\
  \texttt{ 183784461015463: 14}\\
  \texttt{ 184701347515916: 14}\\
  \texttt{ 185088896366668: 14}\\
  \texttt{ 186808382492415: 14}\\
  \texttt{ 186918169826743: 14}\\
  \texttt{ 187092030144715: 14}\\
  \texttt{ 187307874141343: 14}

  \texttt{ 188944592124867: 14}\\
  \texttt{ 189448532664847: 14}\\
  \texttt{ 190570407239175: 14}\\
  \texttt{ 191220351609248: 14}\\
  \texttt{ 195542659566771: 14}\\
  \texttt{ 195939937534444: 14}\\
  \texttt{ 196004175806212: 14}\\
  \texttt{ 196499621658320: 14}\\
  \texttt{ 197871374826512: 14}\\
  \texttt{ 200556595958199: 14}\\
  \texttt{ 202739128985644: 14}\\
  \texttt{ 203063641884512: 14}\\
  \texttt{ 207190567223215: 14}\\
  \texttt{ 209923199367424: 14}\\
  \texttt{ 210142763681372: 14}\\
  \texttt{ 213260842140639: 14}\\
  \texttt{ 215009714009563: 14}\\
  \texttt{ 215848563119223: 14}\\
  \texttt{ 217133130189723: 14}\\
  \texttt{ 218265293317472: 14}

  \texttt{ 220057778398192: 14}\\
  \texttt{ 221743080587067: 14}\\
  \texttt{ 222400208995718: 15}\\
  \texttt{ 225934088816671: 14}\\
  \texttt{ 226641523017112: 14}\\
  \texttt{ 227805968776543: 14}\\
  \texttt{ 230596969918364: 14}\\
  \texttt{ 232162706681368: 14}\\
  \texttt{ 232551408971572: 14}\\
  \texttt{ 236109253309163: 14}\\
  \texttt{ 236450157663247: 14}\\
  \texttt{ 238356890069740: 14}\\
  \texttt{ 238779243514923: 14}\\
  \texttt{ 239380654978924: 14}\\
  \texttt{ 240986242178864: 14}\\
  \texttt{ 241427108122216: 14}\\
  \texttt{ 242382589845668: 14}\\
  \texttt{ 243657462703624: 14}\\
  \texttt{ 245679971474671: 14}\\
  \texttt{ 246052616141572: 14}

  \texttt{ 246557428881855: 14}\\
  \texttt{ 247170265394344: 14}\\
  \texttt{ 249047531240515: 14}\\
  \texttt{ 252568026661839: 14}\\
  \texttt{ 253128636348819: 14}\\
  \texttt{ 254382976818267: 14}\\
  \texttt{ 254918631441220: 14}\\
  \texttt{ 255612478858747: 14}\\
  \texttt{ 257496300981016: 14}\\
  \texttt{ 259036590368419: 14}\\
  \texttt{ 263154806220571: 14}\\
  \texttt{ 264727341474343: 14}\\
  \texttt{ 265524803238220: 14}\\
  \texttt{ 266614047435067: 14}\\
  \texttt{ 269658892858543: 14}\\
  \texttt{ 269733988656339: 14}\\
  \texttt{ 272497539022675: 14}\\
  \texttt{ 274782914525535: 14}\\
  \texttt{ 275458942061071: 14}\\
  \texttt{ 279987459568971: 14}

  \texttt{ 281868811510839: 14}\\
  \texttt{ 282478397137143: 14}\\
  \texttt{ 284287687133012: 14}\\
  \texttt{ 284758974645863: 14}\\
  \texttt{ 285020684617868: 14}\\
  \texttt{ 285838707650812: 14}\\
  \texttt{ 290054274502443: 14}\\
  \texttt{ 293450291054943: 14}\\
  \texttt{ 293504474627944: 14}\\
  \texttt{ 294041644697563: 14}\\
  \texttt{ 295679400256916: 14}\\
  \texttt{ 296998357254568: 14}\\
  \texttt{ 305527536805524: 14}\\
  \texttt{ 305954014084820: 14}\\
  \texttt{ 306596811037768: 14}\\
  \texttt{ 306699434957743: 14}\\
  \texttt{ 311134516185171: 14}\\
  \texttt{ 312686748592171: 14}\\
  \texttt{ 313151373856915: 14}\\
  \texttt{ 314475532576612: 15}

  \texttt{ 316675312313119: 14}\\
  \texttt{ 317407161728419: 14}\\
  \texttt{ 317601894102920: 14}\\
  \texttt{ 317636119018267: 14}\\
  \texttt{ 319079488554512: 14}\\
  \texttt{ 319233677846668: 14}\\
  \texttt{ 319258221819424: 14}\\
  \texttt{ 319700062891240: 14}\\
  \texttt{ 321777171671919: 14}\\
  \texttt{ 322050994737343: 14}\\
  \texttt{ 324149841361023: 14}\\
  \texttt{ 325408134993243: 14}\\
  \texttt{ 325594924038368: 14}\\
  \texttt{ 329052041088819: 14}\\
  \texttt{ 330161375650767: 14}\\
  \texttt{ 331707966393968: 14}\\
  \texttt{ 332494147547439: 14}\\
  \texttt{ 335574316230568: 14}\\
  \texttt{ 336091213643768: 14}\\
  \texttt{ 336496566428944: 14}

  \texttt{ 340034996103640: 14}\\
  \texttt{ 344090011660712: 14}\\
  \texttt{ 344995227576512: 14}\\
  \texttt{ 345664098408175: 14}\\
  \texttt{ 347753503338675: 14}\\
  \texttt{ 352121315027372: 14}\\
  \texttt{ 353399638051995: 14}\\
  \texttt{ 353399638051995: 14}\\
  \texttt{ 359579057071544: 14}\\
  \texttt{ 361271240652940: 14}\\
  \texttt{ 362303412820719: 14}\\
  \texttt{ 363108053058675: 14}\\
  \texttt{ 365145272865844: 14}\\
  \texttt{ 367709652111967: 14}\\
  \texttt{ 368163971088568: 14}\\
  \texttt{ 373323683338567: 14}\\
  \texttt{ 375228738577375: 14}\\
  \texttt{ 375419692667715: 14}\\
  \texttt{ 377517619259072: 14}\\
  \texttt{ 380339486735547: 14}

  \texttt{ 380405207832424: 14}\\
  \texttt{ 380922437852972: 14}\\
  \texttt{ 386244443380075: 14}\\
  \texttt{ 387895592534019: 14}\\
  \texttt{ 389120483450419: 14}\\
  \texttt{ 389202545114864: 14}\\
  \texttt{ 391605463879063: 14}\\
  \texttt{ 395508976515520: 14}\\
  \texttt{ 400234375371968: 14}\\
  \texttt{ 400671041968672: 14}\\
  \texttt{ 404918982367172: 14}\\
  \texttt{ 406076609307472: 14}\\
  \texttt{ 408595748009263: 14}\\
  \texttt{ 409425016308812: 14}\\
  \texttt{ 409976565039619: 14}\\
  \texttt{ 410911526123367: 14}\\
  \texttt{ 411103316180320: 14}\\
  \texttt{ 411213411411667: 14}\\
  \texttt{ 413206218611575: 14}\\
  \texttt{ 415118169415023: 14}

  \texttt{ 417585710635143: 14}\\
  \texttt{ 418450521382467: 14}\\
  \texttt{ 420282062364040: 14}\\
  \texttt{ 422543439918015: 14}\\
  \texttt{ 424282708059424: 14}\\
  \texttt{ 425667633054316: 14}\\
  \texttt{ 427503881624012: 14}\\
  \texttt{ 429557825820068: 14}\\
  \texttt{ 429840103258023: 14}\\
  \texttt{ 431156479024024: 14}\\
  \texttt{ 433465189033472: 14}\\
  \texttt{ 436106024204739: 14}\\
  \texttt{ 437227205770448: 14}\\
  \texttt{ 440134937104444: 14}\\
  \texttt{ 443044716678879: 14}\\
  \texttt{ 444993138872835: 14}\\
  \texttt{ 447501346479943: 14}\\
  \texttt{ 448345051339516: 14}\\
  \texttt{ 449296307366415: 14}\\
  \texttt{ 455370257941612: 14}

  \texttt{ 460816049713516: 14}\\
  \texttt{ 460944893907268: 14}\\
  \texttt{ 462223502363743: 14}\\
  \texttt{ 462623338445415: 14}\\
  \texttt{ 464337974626971: 14}\\
  \texttt{ 464730617738644: 14}\\
  \texttt{ 465319061350144: 14}\\
  \texttt{ 466156525495148: 14}\\
  \texttt{ 468422510178571: 14}\\
  \texttt{ 470249285805043: 14}\\
  \texttt{ 483264807131823: 14}\\
  \texttt{ 484967402899720: 14}\\
  \texttt{ 488337770010644: 14}\\
  \texttt{ 489313678269471: 14}\\
  \texttt{ 489328779185512: 14}\\
  \texttt{ 491989125148063: 14}\\
  \texttt{ 494110321447472: 14}\\
  \texttt{ 495744814185020: 14}\\
  \texttt{ 497497867345815: 14}\\
  \texttt{ 498144133384575: 14}

  \texttt{ 502283706795267: 14}\\
  \texttt{ 504030770086867: 14}\\
  \texttt{ 505192112013716: 14}\\
  \texttt{ 507082020633944: 14}\\
  \texttt{ 507716397212264: 14}\\
  \texttt{ 509500239568071: 14}\\
  \texttt{ 511903182768668: 14}\\
  \texttt{ 512218256743447: 14}\\
  \texttt{ 513761921520019: 14}\\
  \texttt{ 514087760078071: 14}\\
  \texttt{ 514612840136120: 14}\\
  \texttt{ 515135301354871: 14}\\
  \texttt{ 516678590189168: 14}\\
  \texttt{ 516953711998815: 14}\\
  \texttt{ 517448028800275: 14}\\
  \texttt{ 518587427394547: 14}\\
  \texttt{ 521140142483144: 14}\\
  \texttt{ 521497909467116: 14}\\
  \texttt{ 522927463507963: 15}\\
  \texttt{ 525570176949812: 14}

  \texttt{ 527169097223575: 14}\\
  \texttt{ 527253204837220: 14}\\
  \texttt{ 527471563871743: 14}\\
  \texttt{ 529895138521507: 14}\\
  \texttt{ 531074003315271: 14}\\
  \texttt{ 534741979971616: 14}\\
  \texttt{ 535649011280044: 14}\\
  \texttt{ 536083151653420: 14}\\
  \texttt{ 536398113659768: 14}\\
  \texttt{ 543959245890044: 14}\\
  \texttt{ 546752346717643: 14}\\
  \texttt{ 547052616871072: 14}\\
  \texttt{ 547628170513616: 14}\\
  \texttt{ 548725741863471: 14}\\
  \texttt{ 553399175067820: 14}\\
  \texttt{ 554621873513372: 14}\\
  \texttt{ 554802878157639: 14}\\
  \texttt{ 556922150152071: 14}\\
  \texttt{ 558443984700940: 14}\\
  \texttt{ 559066905291471: 14}

  \texttt{ 559226419173772: 14}\\
  \texttt{ 559364456899672: 14}\\
  \texttt{ 559484478417319: 14}\\
  \texttt{ 562057945117963: 14}\\
  \texttt{ 564038769575823: 14}\\
  \texttt{ 564156346792443: 14}\\
  \texttt{ 564873680347243: 14}\\
  \texttt{ 565259131036447: 14}\\
  \texttt{ 565685361038568: 14}\\
  \texttt{ 567697727628339: 14}\\
  \texttt{ 572047852440175: 14}\\
  \texttt{ 576256279996815: 14}\\
  \texttt{ 576732962577112: 14}\\
  \texttt{ 578323272578775: 15}\\
  \texttt{ 579083814768940: 14}\\
  \texttt{ 580043139705412: 14}\\
  \texttt{ 580097298445575: 14}\\
  \texttt{ 586118089764872: 14}\\
  \texttt{ 590435466998019: 14}\\
  \texttt{ 592063024855443: 14}

  \texttt{ 594204306700012: 14}\\
  \texttt{ 595198607956244: 14}\\
  \texttt{ 595584909619568: 14}\\
  \texttt{ 597567054281944: 15}\\
  \texttt{ 599951435679044: 14}\\
  \texttt{ 600172654942663: 14}\\
  \texttt{ 600648730444048: 15}\\
  \texttt{ 605198689272339: 14}\\
  \texttt{ 605404255157871: 14}\\
  \texttt{ 605549016936367: 14}\\
  \texttt{ 607515358116319: 14}\\
  \texttt{ 608733021280615: 14}\\
  \texttt{ 609167337956143: 14}\\
  \texttt{ 609251840049039: 14}\\
  \texttt{ 609544770601767: 14}\\
  \texttt{ 613760263741842: 15}\\
  \texttt{ 613902879206943: 14}\\
  \texttt{ 618265083777667: 14}\\
  \texttt{ 623681524366316: 14}\\
  \texttt{ 623688579818115: 14}

  \texttt{ 625331591972512: 14}\\
  \texttt{ 629356172734843: 14}\\
  \texttt{ 630342645922447: 14}\\
  \texttt{ 630521476160775: 14}\\
  \texttt{ 630682739088171: 14}\\
  \texttt{ 632610060871120: 14}\\
  \texttt{ 635754567283063: 14}\\
  \texttt{ 638983898877967: 14}\\
  \texttt{ 640021957862119: 14}\\
  \texttt{ 641424303929312: 14}\\
  \texttt{ 642942254430947: 14}\\
  \texttt{ 645030379406140: 14}\\
  \texttt{ 646960637293540: 14}\\
  \texttt{ 648288449798047: 14}\\
  \texttt{ 652470316244415: 14}\\
  \texttt{ 653292828286964: 15}\\
  \texttt{ 654399277520020: 14}\\
  \texttt{ 654412768469468: 14}\\
  \texttt{ 656448935132815: 14}\\
  \texttt{ 657982789717623: 14}

  \texttt{ 659286588125067: 14}\\
  \texttt{ 659950058095972: 14}\\
  \texttt{ 661539730781564: 14}\\
  \texttt{ 664951731787375: 14}\\
  \texttt{ 668884718396140: 14}\\
  \texttt{ 669202906742620: 14}\\
  \texttt{ 670091066621216: 14}\\
  \texttt{ 671422822611920: 14}\\
  \texttt{ 672051897841443: 14}\\
  \texttt{ 674093232312567: 14}\\
  \texttt{ 674566232156439: 14}\\
  \texttt{ 675225381031719: 14}\\
  \texttt{ 676589954978672: 14}\\
  \texttt{ 681024988602015: 14}\\
  \texttt{ 681063167276475: 14}\\
  \texttt{ 682380503038815: 14}\\
  \texttt{ 685855603554668: 14}\\
  \texttt{ 690055423200944: 14}\\
  \texttt{ 690626108425119: 14}\\
  \texttt{ 692276022208168: 14}

  \texttt{ 693247847733764: 14}\\
  \texttt{ 693537160604912: 14}\\
  \texttt{ 694612301366068: 14}\\
  \texttt{ 695038087115944: 14}\\
  \texttt{ 696458073649371: 14}\\
  \texttt{ 697013052150607: 14}\\
  \texttt{ 714313556294415: 14}\\
  \texttt{ 715958393580340: 14}\\
  \texttt{ 716713158909039: 14}\\
  \texttt{ 718777318719075: 14}\\
  \texttt{ 720699546486616: 14}\\
  \texttt{ 725167937040163: 14}\\
  \texttt{ 725374003550372: 14}\\
  \texttt{ 731981709091923: 14}\\
  \texttt{ 737263383596723: 14}\\
  \texttt{ 742745255280344: 14}\\
  \texttt{ 744987401798047: 14}\\
  \texttt{ 745404559407844: 14}\\
  \texttt{ 745655995850812: 14}\\
  \texttt{ 746584725308667: 14}

  \texttt{ 746867922000243: 14}\\
  \texttt{ 748746925993412: 14}\\
  \texttt{ 749541430037716: 14}\\
  \texttt{ 749983772544919: 14}\\
  \texttt{ 750046468775019: 14}\\
  \texttt{ 751148913299912: 14}\\
  \texttt{ 751407210639872: 14}\\
  \texttt{ 752174772587440: 14}\\
  \texttt{ 752981108814763: 14}\\
  \texttt{ 754770839204167: 14}\\
  \texttt{ 754824275636047: 14}\\
  \texttt{ 756143836911943: 15}\\
  \texttt{ 756251775056320: 14}\\
  \texttt{ 756344874733063: 14}\\
  \texttt{ 764344642020820: 14}\\
  \texttt{ 766225770874063: 14}\\
  \texttt{ 766893131573947: 14}\\
  \texttt{ 767966021500362: 15}\\
  \texttt{ 767968499800015: 14}\\
  \texttt{ 768674238568747: 14}

  \texttt{ 770743364924123: 14}\\
  \texttt{ 771683149300147: 14}\\
  \texttt{ 774025372373947: 14}\\
  \texttt{ 774226249534168: 14}\\
  \texttt{ 777550215542343: 14}\\
  \texttt{ 778669420964012: 14}\\
  \texttt{ 779306057871019: 14}\\
  \texttt{ 779832929844519: 14}\\
  \texttt{ 781502712052444: 14}\\
  \texttt{ 785047596738812: 14}\\
  \texttt{ 786704724452140: 14}\\
  \texttt{ 793108500249939: 14}\\
  \texttt{ 793384336642112: 14}\\
  \texttt{ 793526309144899: 14}\\
  \texttt{ 798720992419220: 14}\\
  \texttt{ 799368239412412: 15}\\
  \texttt{ 800339157545739: 14}\\
  \texttt{ 801175803467271: 14}\\
  \texttt{ 801762109679671: 14}\\
  \texttt{ 801856955957343: 14}

  \texttt{ 802913153569540: 14}\\
  \texttt{ 804589628217315: 14}\\
  \texttt{ 804971205522367: 14}\\
  \texttt{ 805150851911144: 14}\\
  \texttt{ 805863756844663: 14}\\
  \texttt{ 808004903190164: 14}\\
  \texttt{ 815534583455564: 14}\\
  \texttt{ 817210181456512: 14}\\
  \texttt{ 817717356345867: 14}\\
  \texttt{ 818380203952975: 14}\\
  \texttt{ 820493306411948: 14}\\
  \texttt{ 827322648955575: 14}\\
  \texttt{ 828523770567172: 14}\\
  \texttt{ 828960485598268: 14}\\
  \texttt{ 830446316750572: 14}\\
  \texttt{ 831857440143975: 14}\\
  \texttt{ 832362100496648: 14}\\
  \texttt{ 832835838515415: 15}\\
  \texttt{ 836989079128315: 14}\\
  \texttt{ 841345412213967: 14}

  \texttt{ 843005840098112: 14}\\
  \texttt{ 845891516833240: 15}\\
  \texttt{ 846385278325719: 14}\\
  \texttt{ 847029265989639: 14}\\
  \texttt{ 849010066174664: 14}\\
  \texttt{ 852062962760115: 14}\\
  \texttt{ 854190873768544: 14}\\
  \texttt{ 858846758497540: 14}\\
  \texttt{ 859224196190672: 14}\\
  \texttt{ 860925838900147: 14}\\
  \texttt{ 860993595744315: 14}\\
  \texttt{ 861085050909567: 14}\\
  \texttt{ 863227795245567: 14}\\
  \texttt{ 863431298181772: 14}\\
  \texttt{ 865929254115939: 14}\\
  \texttt{ 866466818659268: 14}\\
  \texttt{ 866936704513840: 14}\\
  \texttt{ 867222447771975: 14}\\
  \texttt{ 867661057013416: 14}\\
  \texttt{ 869245191151244: 14}

  \texttt{ 871324517716712: 14}\\
  \texttt{ 875444799853275: 14}\\
  \texttt{ 875839918708015: 14}\\
  \texttt{ 878265307555724: 14}\\
  \texttt{ 879322694418366: 15}\\
  \texttt{ 880525306975363: 14}\\
  \texttt{ 881640925254915: 14}\\
  \texttt{ 882268844913267: 14}\\
  \texttt{ 882682509038947: 14}\\
  \texttt{ 883206968662816: 14}\\
  \texttt{ 883542349549519: 14}\\
  \texttt{ 883556453040775: 14}\\
  \texttt{ 883881473427063: 14}\\
  \texttt{ 886198679713844: 14}\\
  \texttt{ 888942909642235: 14}\\
  \texttt{ 889338807630916: 14}\\
  \texttt{ 893148922010420: 14}\\
  \texttt{ 893453130778471: 14}\\
  \texttt{ 896074606701871: 14}\\
  \texttt{ 898960342262371: 14}

  \texttt{ 900184106292243: 14}\\
  \texttt{ 905013654483463: 14}\\
  \texttt{ 908708744098519: 14}\\
  \texttt{ 909488438749924: 14}\\
  \texttt{ 909543756082075: 14}\\
  \texttt{ 910485963438812: 14}\\
  \texttt{ 910939607548215: 14}\\
  \texttt{ 914710850261368: 14}\\
  \texttt{ 915874672535872: 14}\\
  \texttt{ 921075883910367: 14}\\
  \texttt{ 922103642978115: 14}\\
  \texttt{ 923530537312963: 14}\\
  \texttt{ 927203200039519: 14}\\
  \texttt{ 930134931993943: 14}\\
  \texttt{ 931221555860319: 14}\\
  \texttt{ 932716002461444: 14}\\
  \texttt{ 933237075697971: 14}\\
  \texttt{ 933877636816912: 14}\\
  \texttt{ 937133041050871: 14}\\
  \texttt{ 939291743594715: 14}

  \texttt{ 939785628473367: 14}\\
  \texttt{ 942328877096744: 14}\\
  \texttt{ 942574927437763: 14}\\
  \texttt{ 944975516649344: 14}\\
  \texttt{ 945088459316216: 14}\\
  \texttt{ 948680596877343: 14}\\
  \texttt{ 949309990294623: 14}\\
  \texttt{ 949746830648947: 14}\\
  \texttt{ 950061620946471: 14}\\
  \texttt{ 950601981839343: 14}\\
  \texttt{ 950793853973012: 14}\\
  \texttt{ 951517018588472: 14}\\
  \texttt{ 952385497232415: 14}\\
  \texttt{ 954139706460339: 14}\\
  \texttt{ 956902291333467: 14}\\
  \texttt{ 960447784058367: 14}\\
  \texttt{ 960665050230568: 14}\\
  \texttt{ 963244642044616: 15}\\
  \texttt{ 963573079631371: 14}\\
  \texttt{ 963907093228112: 14}

  \texttt{ 964857408750124: 14}\\
  \texttt{ 965154340800620: 14}\\
  \texttt{ 968465895652924: 14}\\
  \texttt{ 972075505941871: 14}\\
  \texttt{ 972256492049738: 15}\\
  \texttt{ 976739882486524: 14}\\
  \texttt{ 978405629137168: 15}\\
  \texttt{ 979184624144963: 14}\\
  \texttt{ 979989557540271: 14}\\
  \texttt{ 980233422668512: 14}\\
  \texttt{ 982242475921816: 14}\\
  \texttt{ 982416437892123: 14}\\
  \texttt{ 985136118263368: 14}\\
  \texttt{ 986031825937324: 14}\\
  \texttt{ 986433851901319: 14}\\
  \texttt{ 988704845283543: 14}\\
  \texttt{ 992746092578947: 14}\\
  \texttt{ 993177485791244: 14}\\
  \texttt{ 993548289819832: 14}\\
  \texttt{ 995708686411671: 14}

  \texttt{ 996117970250175: 14}\\
  \texttt{ 1001746828536667: 14}\\
  \texttt{ 1001818577007512: 14}\\
  \texttt{ 1002018841671543: 14}\\
  \texttt{ 1002034345671375: 14}\\
  \texttt{ 1002648412169043: 14}\\
  \texttt{ 1003794956189043: 14}\\
  \texttt{ 1004005499786619: 14}\\
  \texttt{ 1004209127982268: 14}\\
  \texttt{ 1008197006283124: 14}\\
  \texttt{ 1009171678516071: 14}\\
  \texttt{ 1012531696254871: 14}\\
  \texttt{ 1013653212226963: 14}\\
  \texttt{ 1015370534966667: 14}\\
  \texttt{ 1020784200859312: 14}\\
  \texttt{ 1022923129349715: 14}\\
  \texttt{ 1025990073309039: 14}\\
  \texttt{ 1026512980481275: 14}\\
  \texttt{ 1033102159338967: 14}\\
  \texttt{ 1033783215814412: 15}

  \texttt{ 1034035742415472: 14}\\
  \texttt{ 1036122258114872: 14}\\
  \texttt{ 1037509519216372: 14}\\
  \texttt{ 1038525265147767: 14}\\
  \texttt{ 1039430264702740: 14}\\
  \texttt{ 1039795988929243: 14}\\
  \texttt{ 1041478162302968: 14}\\
  \texttt{ 1041942667340175: 14}\\
  \texttt{ 1044638929095172: 14}\\
  \texttt{ 1048475532604472: 14}\\
  \texttt{ 1048883704775512: 14}\\
  \texttt{ 1051400186963619: 14}\\
  \texttt{ 1052886422194568: 14}\\
  \texttt{ 1053541499101747: 14}\\
  \texttt{ 1054747186045215: 14}\\
  \texttt{ 1054925033551467: 14}\\
  \texttt{ 1055003683796620: 14}\\
  \texttt{ 1059627027736663: 14}\\
  \texttt{ 1066286581275220: 14}\\
  \texttt{ 1069775119933575: 14}

  \texttt{ 1070357948515520: 14}\\
  \texttt{ 1072788677199772: 14}\\
  \texttt{ 1077314059441448: 14}\\
  \texttt{ 1079411027724764: 14}\\
  \texttt{ 1080634613744212: 14}\\
  \texttt{ 1081580487480871: 14}\\
  \texttt{ 1081949928009344: 14}\\
  \texttt{ 1082937408522920: 14}\\
  \texttt{ 1085507184675039: 14}\\
  \texttt{ 1087428609050372: 14}\\
  \texttt{ 1088238942140372: 14}\\
  \texttt{ 1088310312678975: 14}\\
  \texttt{ 1088555799244120: 14}\\
  \texttt{ 1095339438292840: 14}\\
  \texttt{ 1098067394516047: 14}\\
  \texttt{ 1099245724468064: 14}\\
  \texttt{ 1100484537453548: 14}\\
  \texttt{ 1103811280126971: 14}\\
  \texttt{ 1104101110479463: 14}\\
  \texttt{ 1104559024531875: 14}

  \texttt{ 1104646894965971: 14}\\
  \texttt{ 1105034346372508: 14}\\
  \texttt{ 1108864900899620: 14}\\
  \texttt{ 1109920747665247: 14}\\
  \texttt{ 1112151945632739: 14}\\
  \texttt{ 1116398771645415: 14}\\
  \texttt{ 1117897185360967: 14}\\
  \texttt{ 1117924125095716: 14}\\
  \texttt{ 1127336666115643: 14}\\
  \texttt{ 1128591791523259: 14}\\
  \texttt{ 1133510837444575: 15}\\
  \texttt{ 1133594105357175: 14}\\
  \texttt{ 1134974346752823: 14}\\
  \texttt{ 1138343206042816: 14}\\
  \texttt{ 1138725938961916: 14}\\
  \texttt{ 1140027266962816: 14}\\
  \texttt{ 1144881144995019: 14}\\
  \texttt{ 1149348516294175: 14}\\
  \texttt{ 1150052702434748: 14}\\
  \texttt{ 1154342748167367: 14}

  \texttt{ 1156228129368548: 14}\\
  \texttt{ 1161407834536904: 14}\\
  \texttt{ 1163346703776519: 14}\\
  \texttt{ 1164643648053171: 14}\\
  \texttt{ 1165975498624616: 14}\\
  \texttt{ 1167327002818520: 14}\\
  \texttt{ 1171942610165467: 14}\\
  \texttt{ 1172079922329015: 14}\\
  \texttt{ 1173834136655512: 14}\\
  \texttt{ 1175385747574816: 14}\\
  \texttt{ 1175805017397875: 14}\\
  \texttt{ 1177061121449367: 14}\\
  \texttt{ 1179809523998620: 14}\\
  \texttt{ 1179954946128248: 14}\\
  \texttt{ 1180678135294375: 14}\\
  \texttt{ 1181978268885772: 14}\\
  \texttt{ 1185337560178915: 14}\\
  \texttt{ 1186126171926711: 15}\\
  \texttt{ 1186768253123223: 14}\\
  \texttt{ 1193408620802971: 14}

  \texttt{ 1194186722888967: 14}\\
  \texttt{ 1196212699039516: 14}\\
  \texttt{ 1196600884724275: 14}\\
  \texttt{ 1202767423530144: 14}\\
  \texttt{ 1204986929396019: 15}\\
  \texttt{ 1205711529790220: 14}\\
  \texttt{ 1209777715691971: 14}\\
  \texttt{ 1212331625780667: 14}\\
  \texttt{ 1212940638467864: 14}\\
  \texttt{ 1221418539613447: 14}\\
  \texttt{ 1222744906812400: 14}\\
  \texttt{ 1224793012836868: 14}\\
  \texttt{ 1224871988215544: 14}\\
  \texttt{ 1226481168310540: 14}\\
  \texttt{ 1228197999757768: 14}\\
  \texttt{ 1228525731675044: 14}\\
  \texttt{ 1228814132845075: 14}\\
  \texttt{ 1229109121984364: 14}\\
  \texttt{ 1229788641055371: 14}\\
  \texttt{ 1232238913378063: 14}

  \texttt{ 1233348651154239: 14}\\
  \texttt{ 1234686275197964: 14}\\
  \texttt{ 1239279595984820: 14}\\
  \texttt{ 1241355851278316: 14}\\
  \texttt{ 1241535222032775: 14}\\
  \texttt{ 1241833805783144: 14}\\
  \texttt{ 1242385435271824: 14}\\
  \texttt{ 1242772714742264: 14}\\
  \texttt{ 1243088562651020: 14}\\
  \texttt{ 1244449578411475: 14}\\
  \texttt{ 1246415138138419: 14}\\
  \texttt{ 1247335327266847: 14}\\
  \texttt{ 1248232560022139: 14}\\
  \texttt{ 1248381650962268: 14}\\
  \texttt{ 1248850080568168: 14}\\
  \texttt{ 1250158695994448: 14}\\
  \texttt{ 1254926989630520: 14}\\
  \texttt{ 1259211149861775: 14}\\
  \texttt{ 1262066854796724: 14}\\
  \texttt{ 1263760954523115: 14}

  \texttt{ 1265916121650267: 14}\\
  \texttt{ 1269647734190948: 14}\\
  \texttt{ 1270753618232571: 14}\\
  \texttt{ 1273168255360443: 14}\\
  \texttt{ 1273484457601684: 14}\\
  \texttt{ 1274417009922067: 14}\\
  \texttt{ 1276501376160764: 14}\\
  \texttt{ 1277913719609872: 14}\\
  \texttt{ 1278011372665075: 14}\\
  \texttt{ 1278545995919564: 14}\\
  \texttt{ 1284135477814924: 14}\\
  \texttt{ 1284399141309475: 14}\\
  \texttt{ 1286587320067240: 14}\\
  \texttt{ 1286783060101316: 14}\\
  \texttt{ 1288385500915016: 14}\\
  \texttt{ 1288432127884840: 14}\\
  \texttt{ 1288891538318563: 14}\\
  \texttt{ 1289166067168675: 14}\\
  \texttt{ 1290226383402847: 14}\\
  \texttt{ 1293086296630672: 14}

  \texttt{ 1296090481030072: 14}\\
  \texttt{ 1299076303745419: 14}\\
  \texttt{ 1303769013983224: 14}\\
  \texttt{ 1304154331510472: 14}\\
  \texttt{ 1304995870394872: 14}\\
  \texttt{ 1305346650834644: 14}\\
  \texttt{ 1305871231904647: 14}\\
  \texttt{ 1306041418307968: 14}\\
  \texttt{ 1306915315508319: 14}\\
  \texttt{ 1308641301607615: 14}\\
  \texttt{ 1311959805470872: 14}\\
  \texttt{ 1314117174266468: 14}\\
  \texttt{ 1315236429605671: 14}\\
  \texttt{ 1319548391635868: 14}\\
  \texttt{ 1320888477379224: 14}\\
  \texttt{ 1325338093786220: 14}\\
  \texttt{ 1325835749970123: 14}\\
  \texttt{ 1327036511579948: 14}\\
  \texttt{ 1329734721503067: 14}\\
  \texttt{ 1330022636862819: 14}

  \texttt{ 1332197415874167: 14}\\
  \texttt{ 1334789276849224: 14}\\
  \texttt{ 1337690847815368: 14}\\
  \texttt{ 1346926338369267: 14}\\
  \texttt{ 1348684073081872: 14}\\
  \texttt{ 1352245036644160: 14}\\
  \texttt{ 1355700127919264: 14}\\
  \texttt{ 1357351609343775: 14}\\
  \texttt{ 1357421813474415: 14}\\
  \texttt{ 1359028536730540: 14}\\
  \texttt{ 1360309169048264: 14}\\
  \texttt{ 1360893042732975: 14}\\
  \texttt{ 1361871740828067: 14}\\
  \texttt{ 1362270466768971: 14}\\
  \texttt{ 1365366207309343: 14}\\
  \texttt{ 1366451274004323: 14}\\
  \texttt{ 1366784467548819: 14}\\
  \texttt{ 1367083790439848: 14}\\
  \texttt{ 1367383259297167: 14}\\
  \texttt{ 1368253581539444: 14}

  \texttt{ 1368336961042243: 14}\\
  \texttt{ 1369646994319539: 14}\\
  \texttt{ 1371773900857363: 14}\\
  \texttt{ 1373627370439519: 14}\\
  \texttt{ 1374401560917247: 14}\\
  \texttt{ 1377909351902671: 14}\\
  \texttt{ 1380115193843067: 14}\\
  \texttt{ 1380933079224668: 14}\\
  \texttt{ 1381119379017820: 14}\\
  \texttt{ 1382409520712043: 14}\\
  \texttt{ 1385302830539019: 14}\\
  \texttt{ 1385587774437019: 14}\\
  \texttt{ 1386960561200419: 14}\\
  \texttt{ 1387117341595923: 14}\\
  \texttt{ 1388495651104239: 14}\\
  \texttt{ 1390238414196919: 14}\\
  \texttt{ 1393698794085339: 14}\\
  \texttt{ 1393777393073216: 14}\\
  \texttt{ 1395355738900867: 14}\\
  \texttt{ 1396179883418012: 14}

  \texttt{ 1397615792814520: 14}\\
  \texttt{ 1405717696947220: 14}\\
  \texttt{ 1405786859261467: 14}\\
  \texttt{ 1407929584482375: 14}\\
  \texttt{ 1407999508293715: 15}\\
  \texttt{ 1408232649110119: 14}\\
  \texttt{ 1411298770780916: 14}\\
  \texttt{ 1413326736000616: 15}\\
  \texttt{ 1413666099541419: 14}\\
  \texttt{ 1415813318783116: 14}\\
  \texttt{ 1417874670199663: 14}\\
  \texttt{ 1418810627961943: 14}\\
  \texttt{ 1424900490712567: 14}\\
  \texttt{ 1426925298986720: 14}\\
  \texttt{ 1431261511856720: 14}\\
  \texttt{ 1433620928658940: 14}\\
  \texttt{ 1435521793104367: 14}\\
  \texttt{ 1436412472238068: 14}\\
  \texttt{ 1437381568682719: 14}\\
  \texttt{ 1437422608342743: 14}

  \texttt{ 1439503148554712: 14}\\
  \texttt{ 1442542092130268: 14}\\
  \texttt{ 1442554969087239: 14}\\
  \texttt{ 1443658505841871: 14}\\
  \texttt{ 1447396507462364: 14}\\
  \texttt{ 1447482377144672: 14}\\
  \texttt{ 1450543057901044: 14}\\
  \texttt{ 1451799493433367: 14}\\
  \texttt{ 1452057869663872: 14}\\
  \texttt{ 1452715837505835: 14}\\
  \texttt{ 1453754769560044: 14}\\
  \texttt{ 1454156022591764: 14}\\
  \texttt{ 1457807362155512: 14}\\
  \texttt{ 1463243595724267: 14}\\
  \texttt{ 1463595146758923: 14}\\
  \texttt{ 1463892212483912: 14}\\
  \texttt{ 1465336399306112: 14}\\
  \texttt{ 1465712896468112: 14}\\
  \texttt{ 1466761674009068: 14}\\
  \texttt{ 1467377292644648: 14}

  \texttt{ 1467604886854012: 14}\\
  \texttt{ 1468871459602775: 14}\\
  \texttt{ 1472107952982319: 14}\\
  \texttt{ 1473270788489344: 14}\\
  \texttt{ 1475398828787019: 14}\\
  \texttt{ 1476838803571315: 14}\\
  \texttt{ 1476966230087371: 14}\\
  \texttt{ 1477654731062072: 14}\\
  \texttt{ 1479484025622375: 14}\\
  \texttt{ 1480722764457248: 14}\\
  \texttt{ 1481047660013048: 14}\\
  \texttt{ 1481772335903756: 14}\\
  \texttt{ 1481927416642275: 14}\\
  \texttt{ 1482442771986122: 15}\\
  \texttt{ 1483245859463620: 14}\\
  \texttt{ 1484446836982148: 14}\\
  \texttt{ 1485645349660143: 14}\\
  \texttt{ 1487154736562012: 14}\\
  \texttt{ 1487173060233044: 14}\\
  \texttt{ 1487676816408175: 14}

  \texttt{ 1488963081044347: 14}\\
  \texttt{ 1489470897131072: 14}\\
  \texttt{ 1490190055193816: 14}\\
  \texttt{ 1491773597899244: 14}\\
  \texttt{ 1491797412912464: 14}\\
  \texttt{ 1494909216836871: 14}\\
  \texttt{ 1498879222815939: 14}\\
  \texttt{ 1501666773882415: 14}\\
  \texttt{ 1504379852873368: 14}\\
  \texttt{ 1504933747547523: 14}\\
  \texttt{ 1506959819306871: 14}\\
  \texttt{ 1511022491247219: 14}\\
  \texttt{ 1512135112078448: 14}\\
  \texttt{ 1516330608913016: 14}\\
  \texttt{ 1517448600782367: 14}\\
  \texttt{ 1517898009472244: 14}\\
  \texttt{ 1519515603058839: 14}\\
  \texttt{ 1521623333594347: 14}\\
  \texttt{ 1523513005371115: 14}\\
  \texttt{ 1525536143887371: 14}

  \texttt{ 1528828993593547: 14}\\
  \texttt{ 1530122386059044: 14}\\
  \texttt{ 1540010187125648: 14}\\
  \texttt{ 1545116627480920: 14}\\
  \texttt{ 1547788420772372: 14}\\
  \texttt{ 1550677127043068: 14}\\
  \texttt{ 1551482181488715: 14}\\
  \texttt{ 1553234687735416: 14}\\
  \texttt{ 1554739296399044: 14}\\
  \texttt{ 1556366410106647: 14}\\
  \texttt{ 1558028716710235: 14}\\
  \texttt{ 1558488018147027: 14}\\
  \texttt{ 1562538586065712: 14}\\
  \texttt{ 1562687462755744: 14}\\
  \texttt{ 1563210552123367: 14}\\
  \texttt{ 1563801041614120: 14}\\
  \texttt{ 1566907317991671: 14}\\
  \texttt{ 1569147026407468: 14}\\
  \texttt{ 1569161994316323: 14}\\
  \texttt{ 1570207493001319: 14}

  \texttt{ 1572913250310643: 14}\\
  \texttt{ 1573369965334420: 14}\\
  \texttt{ 1578509853004275: 14}\\
  \texttt{ 1582604057410567: 14}\\
  \texttt{ 1583548831782547: 14}\\
  \texttt{ 1589645758039575: 14}\\
  \texttt{ 1590972226727272: 14}\\
  \texttt{ 1591011230955667: 14}\\
  \texttt{ 1591269736815772: 14}\\
  \texttt{ 1592849311091115: 14}\\
  \texttt{ 1593481013140120: 14}\\
  \texttt{ 1594074945400675: 14}\\
  \texttt{ 1595628941363116: 14}\\
  \texttt{ 1597456515585544: 14}\\
  \texttt{ 1598628474432944: 14}\\
  \texttt{ 1598978039422516: 14}\\
  \texttt{ 1602472376770171: 14}\\
  \texttt{ 1604834954915120: 14}\\
  \texttt{ 1606148115010963: 14}\\
  \texttt{ 1606326677414943: 14}

  \texttt{ 1613070469720743: 14}\\
  \texttt{ 1618009556832483: 14}\\
  \texttt{ 1618107473165667: 14}\\
  \texttt{ 1619455073970723: 14}\\
  \texttt{ 1620346504333148: 14}\\
  \texttt{ 1621461996536343: 14}\\
  \texttt{ 1626554930798963: 14}\\
  \texttt{ 1627761098973472: 14}\\
  \texttt{ 1628012948512371: 14}\\
  \texttt{ 1630516325175124: 14}\\
  \texttt{ 1631413849629967: 14}\\
  \texttt{ 1633785156332619: 14}\\
  \texttt{ 1636196401149775: 14}\\
  \texttt{ 1636737342986744: 14}\\
  \texttt{ 1638598256409643: 14}\\
  \texttt{ 1639181756316039: 14}\\
  \texttt{ 1640741786210644: 15}\\
  \texttt{ 1641874595513444: 14}\\
  \texttt{ 1642243362834643: 14}\\
  \texttt{ 1643123569629175: 14}

  \texttt{ 1643864314457824: 14}\\
  \texttt{ 1643951944093120: 14}\\
  \texttt{ 1644766123620067: 14}\\
  \texttt{ 1645339558397012: 14}\\
  \texttt{ 1645799076456479: 14}\\
  \texttt{ 1647423474094275: 15}\\
  \texttt{ 1648391977873663: 14}\\
  \texttt{ 1652309411883424: 14}\\
  \texttt{ 1660800635485971: 14}\\
  \texttt{ 1660865185948148: 14}\\
  \texttt{ 1661728614821043: 14}\\
  \texttt{ 1666823000814115: 14}\\
  \texttt{ 1666895115919911: 15}\\
  \texttt{ 1667064090629348: 14}\\
  \texttt{ 1670379746433271: 14}\\
  \texttt{ 1670675620743591: 14}\\
  \texttt{ 1675833317995947: 14}\\
  \texttt{ 1676071949560072: 14}\\
  \texttt{ 1676244928307019: 14}\\
  \texttt{ 1679577622169947: 14}

  \texttt{ 1684825012476548: 14}\\
  \texttt{ 1686477047556615: 14}\\
  \texttt{ 1687213902823575: 14}\\
  \texttt{ 1690496783493243: 14}\\
  \texttt{ 1692972969824716: 14}\\
  \texttt{ 1693784196921472: 14}\\
  \texttt{ 1694569833029715: 14}\\
  \texttt{ 1698689603893364: 14}\\
  \texttt{ 1699499966949015: 14}\\
  \texttt{ 1702457473333067: 14}\\
  \texttt{ 1702962430850564: 14}\\
  \texttt{ 1703798983604968: 14}\\
  \texttt{ 1707674023177368: 14}\\
  \texttt{ 1709026507289348: 14}\\
  \texttt{ 1711127138739271: 14}\\
  \texttt{ 1712958638298775: 14}\\
  \texttt{ 1714042192544468: 14}\\
  \texttt{ 1714904007298172: 14}\\
  \texttt{ 1718696770586467: 14}\\
  \texttt{ 1719525425635816: 14}

  \texttt{ 1719536777409868: 14}\\
  \texttt{ 1720931061947775: 14}\\
  \texttt{ 1721331360054267: 14}\\
  \texttt{ 1721645029117168: 14}\\
  \texttt{ 1722044570499944: 14}\\
  \texttt{ 1722488125302763: 14}\\
  \texttt{ 1724152140929824: 14}\\
  \texttt{ 1724933500400068: 14}\\
  \texttt{ 1727557469766843: 14}\\
  \texttt{ 1728562371923214: 15}\\
  \texttt{ 1730031047756619: 14}\\
  \texttt{ 1731484894050844: 14}\\
  \texttt{ 1732747657503724: 14}\\
  \texttt{ 1736854893619216: 14}\\
  \texttt{ 1737788100001216: 14}\\
  \texttt{ 1738761417495772: 14}\\
  \texttt{ 1739512267563471: 14}\\
  \texttt{ 1740418546481143: 14}\\
  \texttt{ 1741323697328715: 14}\\
  \texttt{ 1741499196595323: 14}

  \texttt{ 1744207771889420: 14}\\
  \texttt{ 1747065344142716: 14}\\
  \texttt{ 1753727494321671: 14}\\
  \texttt{ 1753756749977912: 14}\\
  \texttt{ 1755893286256075: 14}\\
  \texttt{ 1758403224580219: 14}\\
  \texttt{ 1759766179304572: 14}\\
  \texttt{ 1761631458029215: 14}\\
  \texttt{ 1762777026874844: 14}\\
  \texttt{ 1763896767879067: 14}\\
  \texttt{ 1764028782444039: 14}\\
  \texttt{ 1769412815868571: 14}\\
  \texttt{ 1770286724382939: 14}\\
  \texttt{ 1770319447463664: 14}\\
  \texttt{ 1771160152656471: 14}\\
  \texttt{ 1771531434341871: 14}\\
  \texttt{ 1773666944874367: 14}\\
  \texttt{ 1774745404186363: 14}\\
  \texttt{ 1776552305734168: 14}\\
  \texttt{ 1776995556595324: 14}

  \texttt{ 1777384746760544: 14}\\
  \texttt{ 1777726915270964: 14}\\
  \texttt{ 1778062194470672: 14}\\
  \texttt{ 1778606864930643: 14}\\
  \texttt{ 1778730977125411: 15}\\
  \texttt{ 1778874147975968: 14}\\
  \texttt{ 1778898738878468: 14}\\
  \texttt{ 1780258245716167: 14}\\
  \texttt{ 1780352824906240: 14}\\
  \texttt{ 1783487063119147: 14}\\
  \texttt{ 1785705624752947: 14}\\
  \texttt{ 1785754281958172: 14}\\
  \texttt{ 1790762526073843: 14}\\
  \texttt{ 1790782158828175: 14}\\
  \texttt{ 1791967641810820: 14}\\
  \texttt{ 1792108460296772: 14}\\
  \texttt{ 1792985029120744: 14}\\
  \texttt{ 1793403628397612: 14}\\
  \texttt{ 1793708000423216: 14}\\
  \texttt{ 1796513848368019: 14}

  \texttt{ 1796537601066220: 14}\\
  \texttt{ 1796616152498420: 14}\\
  \texttt{ 1798432051651575: 14}\\
  \texttt{ 1798570241557612: 14}\\
  \texttt{ 1800543530727764: 14}\\
  \texttt{ 1801132646196247: 14}\\
  \texttt{ 1801737938341672: 14}\\
  \texttt{ 1803300561512475: 14}\\
  \texttt{ 1806229327602544: 14}\\
  \texttt{ 1806420007387664: 14}\\
  \texttt{ 1807179445143915: 14}\\
  \texttt{ 1807359570051568: 14}\\
  \texttt{ 1807978378932943: 14}\\
  \texttt{ 1808697633654843: 14}\\
  \texttt{ 1810475285639024: 14}\\
  \texttt{ 1811392774777419: 14}\\
  \texttt{ 1817141712456464: 14}\\
  \texttt{ 1819460814609471: 14}\\
  \texttt{ 1821683915682315: 14}\\
  \texttt{ 1823992306125423: 14}

  \texttt{ 1824033853889863: 14}\\
  \texttt{ 1828702594086343: 14}\\
  \texttt{ 1829422165038172: 14}\\
  \texttt{ 1831366501093168: 14}\\
  \texttt{ 1831510735969220: 14}\\
  \texttt{ 1836164798052819: 14}\\
  \texttt{ 1837512773264067: 14}\\
  \texttt{ 1841843231840643: 14}\\
  \texttt{ 1843633623956920: 14}\\
  \texttt{ 1846448459520819: 14}\\
  \texttt{ 1847034554672671: 14}\\
  \texttt{ 1847639382148744: 14}\\
  \texttt{ 1849164474145875: 14}\\
  \texttt{ 1851544395975319: 14}\\
  \texttt{ 1851651341446516: 14}\\
  \texttt{ 1853911768061815: 14}\\
  \texttt{ 1854101482412768: 14}\\
  \texttt{ 1856947003600148: 14}\\
  \texttt{ 1861957899199844: 14}\\
  \texttt{ 1864308829729372: 14}

  \texttt{ 1864372578749019: 14}\\
  \texttt{ 1864702872045844: 14}\\
  \texttt{ 1866902589587815: 14}\\
  \texttt{ 1868217355698271: 14}\\
  \texttt{ 1872615107767148: 14}\\
  \texttt{ 1873279953375919: 14}\\
  \texttt{ 1876341001165323: 14}\\
  \texttt{ 1877745174871335: 14}\\
  \texttt{ 1880150683734943: 14}\\
  \texttt{ 1880612263510244: 14}\\
  \texttt{ 1881148336659339: 14}\\
  \texttt{ 1881913094352044: 14}\\
  \texttt{ 1883372612274916: 14}\\
  \texttt{ 1884208591461212: 14}\\
  \texttt{ 1884806123254216: 14}\\
  \texttt{ 1884842437923520: 14}\\
  \texttt{ 1887179922247112: 15}\\
  \texttt{ 1888059957636615: 14}\\
  \texttt{ 1889426811739516: 14}\\
  \texttt{ 1889545086674440: 14}

  \texttt{ 1889852340352515: 14}\\
  \texttt{ 1891143421462772: 15}\\
  \texttt{ 1891633516219419: 14}\\
  \texttt{ 1892915579135312: 14}\\
  \texttt{ 1894280094994768: 14}\\
  \texttt{ 1895290201726568: 14}\\
  \texttt{ 1895702499666115: 14}\\
  \texttt{ 1898873472817275: 14}\\
  \texttt{ 1899099325766572: 14}\\
  \texttt{ 1902328905593824: 14}\\
  \texttt{ 1903285470442364: 14}\\
  \texttt{ 1904508224703764: 14}\\
  \texttt{ 1905086633797843: 14}\\
  \texttt{ 1905211817914112: 14}\\
  \texttt{ 1906866902453948: 14}\\
  \texttt{ 1908828085396472: 14}\\
  \texttt{ 1910781516402620: 14}\\
  \texttt{ 1912507713241068: 14}\\
  \texttt{ 1912554427994943: 14}\\
  \texttt{ 1916240172431524: 14}

  \texttt{ 1918613887955467: 14}\\
  \texttt{ 1919655455212924: 14}\\
  \texttt{ 1923041569781572: 14}\\
  \texttt{ 1928597441899112: 14}\\
  \texttt{ 1928966323590043: 14}\\
  \texttt{ 1930461443301175: 14}\\
  \texttt{ 1932916763495672: 14}\\
  \texttt{ 1936521004957972: 14}\\
  \texttt{ 1937398272671643: 14}\\
  \texttt{ 1940303735659419: 14}\\
  \texttt{ 1942713964651444: 14}\\
  \texttt{ 1948612594475212: 14}\\
  \texttt{ 1955830409450720: 14}\\
  \texttt{ 1958956603200915: 14}\\
  \texttt{ 1960564061350167: 14}\\
  \texttt{ 1960904663068839: 14}\\
  \texttt{ 1964240880507548: 14}\\
  \texttt{ 1969362031335667: 14}\\
  \texttt{ 1970321128365723: 14}\\
  \texttt{ 1970897048504943: 14}

  \texttt{ 1970944607508243: 14}\\
  \texttt{ 1971718226448640: 14}\\
  \texttt{ 1976165854011303: 14}\\
  \texttt{ 1976269537077423: 14}\\
  \texttt{ 1977063276197420: 14}\\
  \texttt{ 1983741593190243: 14}\\
  \texttt{ 1984238009147140: 14}\\
  \texttt{ 1991552006087072: 14}\\
  \texttt{ 1995415129656548: 14}\\
  \texttt{ 1998379505077972: 14}\\
  \texttt{ 2001256697287419: 14}\\
  \texttt{ 2002637251645719: 14}\\
  \texttt{ 2002761684491444: 14}\\
  \texttt{ 2002875318625875: 14}\\
  \texttt{ 2003220183417164: 14}\\
  \texttt{ 2006150963699068: 14}\\
  \texttt{ 2006731022695324: 14}\\
  \texttt{ 2007387550069675: 14}\\
  \texttt{ 2007526171435071: 14}\\
  \texttt{ 2008431846799720: 14}

  \texttt{ 2009463893507312: 14}\\
  \texttt{ 2009611919686768: 14}\\
  \texttt{ 2012399290118320: 14}\\
  \texttt{ 2012401120535663: 14}\\
  \texttt{ 2012919696621296: 14}\\
  \texttt{ 2013809290522468: 14}\\
  \texttt{ 2014590132458647: 14}\\
  \texttt{ 2017870640442619: 14}\\
  \texttt{ 2023463688325575: 14}\\
  \texttt{ 2023646588688219: 14}\\
  \texttt{ 2024193859666144: 14}\\
  \texttt{ 2028313746275444: 14}\\
  \texttt{ 2032854454292343: 14}\\
  \texttt{ 2033310459121171: 14}\\
  \texttt{ 2033665494948368: 14}\\
  \texttt{ 2034215622153548: 14}\\
  \texttt{ 2035605289694047: 14}\\
  \texttt{ 2037431740482268: 14}\\
  \texttt{ 2038198243357744: 14}\\
  \texttt{ 2039873999933644: 14}

  \texttt{ 2041863790689340: 14}\\
  \texttt{ 2042550262864419: 14}\\
  \texttt{ 2045391215653663: 14}\\
  \texttt{ 2045972062427371: 15}\\
  \texttt{ 2047496287236019: 14}\\
  \texttt{ 2048447151132775: 14}\\
  \texttt{ 2049467252699312: 14}\\
  \texttt{ 2051558678845515: 14}\\
  \texttt{ 2052363765864724: 14}\\
  \texttt{ 2052655258180363: 14}\\
  \texttt{ 2052938491183023: 14}\\
  \texttt{ 2053353206472315: 14}\\
  \texttt{ 2055636425706271: 14}\\
  \texttt{ 2058627635625724: 14}\\
  \texttt{ 2063336136540424: 14}\\
  \texttt{ 2063694280707338: 15}\\
  \texttt{ 2064190884746875: 14}\\
  \texttt{ 2066630553518668: 14}\\
  \texttt{ 2069644073574068: 14}\\
  \texttt{ 2070893574483415: 14}

  \texttt{ 2072522431657023: 14}\\
  \texttt{ 2072695070413448: 14}\\
  \texttt{ 2074659093057067: 14}\\
  \texttt{ 2075992189719163: 14}\\
  \texttt{ 2077016306119063: 14}\\
  \texttt{ 2080143676535812: 14}\\
  \texttt{ 2082990511757067: 14}\\
  \texttt{ 2084601015148767: 14}\\
  \texttt{ 2086439630020964: 14}\\
  \texttt{ 2086860307930216: 14}\\
  \texttt{ 2088397237257771: 14}\\
  \texttt{ 2088403713602468: 14}\\
  \texttt{ 2088784054551463: 14}\\
  \texttt{ 2090125006694620: 14}\\
  \texttt{ 2090363746601175: 14}\\
  \texttt{ 2091096982924143: 14}\\
  \texttt{ 2094572183372067: 14}\\
  \texttt{ 2099070564541675: 14}\\
  \texttt{ 2100613839529216: 15}\\
  \texttt{ 2101499885298248: 14}

  \texttt{ 2103568614391616: 14}\\
  \texttt{ 2103605342143316: 14}\\
  \texttt{ 2104519001695747: 14}\\
  \texttt{ 2106313655656219: 14}\\
  \texttt{ 2107465770178844: 14}\\
  \texttt{ 2108632768935675: 14}\\
  \texttt{ 2109734217703168: 14}\\
  \texttt{ 2114427737578544: 14}\\
  \texttt{ 2115726952146424: 14}\\
  \texttt{ 2117625747868772: 14}\\
  \texttt{ 2118743875265263: 14}\\
  \texttt{ 2121976644943516: 14}\\
  \texttt{ 2125663855802215: 14}\\
  \texttt{ 2129116015084220: 14}\\
  \texttt{ 2129724903314068: 14}\\
  \texttt{ 2135910972244120: 14}\\
  \texttt{ 2136243257632663: 14}\\
  \texttt{ 2137584978892447: 14}\\
  \texttt{ 2139649840224015: 14}\\
  \texttt{ 2142912970530512: 14}

  \texttt{ 2146197499364672: 14}\\
  \texttt{ 2146574044249144: 14}\\
  \texttt{ 2146958283981247: 14}\\
  \texttt{ 2149597852664372: 14}\\
  \texttt{ 2152018359078416: 14}\\
  \texttt{ 2154871323574999: 14}\\
  \texttt{ 2156289455973868: 14}\\
  \texttt{ 2164704515128167: 14}\\
  \texttt{ 2164875105133771: 14}\\
  \texttt{ 2167905752495619: 14}\\
  \texttt{ 2169293537330620: 14}\\
  \texttt{ 2169499273571931: 14}\\
  \texttt{ 2173002772608340: 14}\\
  \texttt{ 2174287510729168: 14}\\
  \texttt{ 2177989619519744: 14}\\
  \texttt{ 2178455860754116: 14}\\
  \texttt{ 2179148650658044: 14}\\
  \texttt{ 2180340668491972: 14}\\
  \texttt{ 2185992033453124: 14}\\
  \texttt{ 2187143654743168: 14}

  \texttt{ 2190716067196023: 14}\\
  \texttt{ 2191481925331219: 14}\\
  \texttt{ 2191493703845716: 14}\\
  \texttt{ 2191511513192119: 14}\\
  \texttt{ 2191609585248775: 14}\\
  \texttt{ 2193642202163919: 14}\\
  \texttt{ 2195817775767619: 14}\\
  \texttt{ 2196913882832575: 14}\\
  \texttt{ 2197967379457364: 14}\\
  \texttt{ 2201721698188743: 14}\\
  \texttt{ 2203405270880768: 14}\\
  \texttt{ 2209961679974647: 14}\\
  \texttt{ 2210476235514775: 14}\\
  \texttt{ 2211455460099975: 14}\\
  \texttt{ 2213636019248167: 14}\\
  \texttt{ 2216030977599272: 14}\\
  \texttt{ 2221533412609372: 14}\\
  \texttt{ 2224252385965519: 14}\\
  \texttt{ 2225144438878516: 14}\\
  \texttt{ 2227322596148143: 14}

  \texttt{ 2227348282719068: 15}\\
  \texttt{ 2230417682683323: 14}\\
  \texttt{ 2235439759782248: 14}\\
  \texttt{ 2238745595588684: 14}\\
  \texttt{ 2239527109197075: 14}\\
  \texttt{ 2240409885171848: 14}\\
  \texttt{ 2241130372596367: 14}\\
  \texttt{ 2242961918497312: 14}\\
  \texttt{ 2244255484668271: 14}\\
  \texttt{ 2244469444899639: 14}\\
  \texttt{ 2244895371983715: 14}\\
  \texttt{ 2245323849310239: 14}\\
  \texttt{ 2249080116729919: 14}\\
  \texttt{ 2249915107141472: 14}\\
  \texttt{ 2250095101089343: 14}\\
  \texttt{ 2251112579233567: 14}\\
  \texttt{ 2252050735467471: 14}\\
  \texttt{ 2253095349006616: 14}\\
  \texttt{ 2254513010999115: 14}\\
  \texttt{ 2255046185932722: 15}

  \texttt{ 2257472370333564: 14}\\
  \texttt{ 2257760158947043: 14}\\
  \texttt{ 2266274792223015: 14}\\
  \texttt{ 2267486845055612: 14}\\
  \texttt{ 2268274802004267: 14}\\
  \texttt{ 2268975651753644: 14}\\
  \texttt{ 2276765335853943: 14}\\
  \texttt{ 2277738474733267: 14}\\
  \texttt{ 2278772355388120: 14}\\
  \texttt{ 2282944765484268: 14}\\
  \texttt{ 2284996884920072: 14}\\
  \texttt{ 2288289389582968: 14}\\
  \texttt{ 2291538049537220: 14}\\
  \texttt{ 2300702726707064: 14}\\
  \texttt{ 2301158097889144: 14}\\
  \texttt{ 2301941812156719: 14}\\
  \texttt{ 2304987456364567: 14}\\
  \texttt{ 2306100268278843: 14}\\
  \texttt{ 2306736740294872: 14}\\
  \texttt{ 2308900763361868: 14}

  \texttt{ 2309324366886920: 14}\\
  \texttt{ 2311226308247120: 14}\\
  \texttt{ 2312623181427819: 14}\\
  \texttt{ 2313259643695119: 14}\\
  \texttt{ 2315287496729175: 14}\\
  \texttt{ 2317004710517716: 14}\\
  \texttt{ 2317487745520448: 14}\\
  \texttt{ 2318011375658372: 14}\\
  \texttt{ 2320780993797068: 14}\\
  \texttt{ 2321080369739215: 14}\\
  \texttt{ 2322401826878071: 14}\\
  \texttt{ 2323456017056108: 14}\\
  \texttt{ 2328161669757015: 14}\\
  \texttt{ 2329676050279120: 14}\\
  \texttt{ 2331543076509248: 14}\\
  \texttt{ 2331724090310439: 14}\\
  \texttt{ 2335223716675172: 14}\\
  \texttt{ 2336453182902039: 14}\\
  \texttt{ 2337016369924240: 14}\\
  \texttt{ 2337476146294120: 14}

  \texttt{ 2337487375443163: 14}\\
  \texttt{ 2337887973246163: 14}\\
  \texttt{ 2338961877517148: 14}\\
  \texttt{ 2341515488015671: 14}\\
  \texttt{ 2342663210984740: 14}\\
  \texttt{ 2347972209573171: 15}\\
  \texttt{ 2352533147125064: 14}\\
  \texttt{ 2352672850946668: 14}\\
  \texttt{ 2354728960111148: 14}\\
  \texttt{ 2355690790768363: 14}\\
  \texttt{ 2358879855857072: 14}\\
  \texttt{ 2361875380381419: 14}\\
  \texttt{ 2363900802934615: 14}\\
  \texttt{ 2364970277291919: 14}\\
  \texttt{ 2369800096881171: 14}\\
  \texttt{ 2369907350225071: 15}\\
  \texttt{ 2370689627114175: 14}\\
  \texttt{ 2370773842633268: 14}\\
  \texttt{ 2372286923874124: 14}\\
  \texttt{ 2372640704766267: 14}

  \texttt{ 2376291067416916: 14}\\
  \texttt{ 2380318374521212: 14}\\
  \texttt{ 2381241548297563: 14}\\
  \texttt{ 2381972780187368: 14}\\
  \texttt{ 2382891774798112: 14}\\
  \texttt{ 2384248736310424: 14}\\
  \texttt{ 2384922510668067: 14}\\
  \texttt{ 2387745162693548: 14}\\
  \texttt{ 2387842010066744: 14}\\
  \texttt{ 2390772039893368: 14}\\
  \texttt{ 2391606735378920: 14}\\
  \texttt{ 2392299717072675: 14}\\
  \texttt{ 2392436730981112: 14}\\
  \texttt{ 2395777956673624: 14}\\
  \texttt{ 2395940439423544: 14}\\
  \texttt{ 2397118860727512: 14}\\
  \texttt{ 2398559938906075: 14}\\
  \texttt{ 2398811085184767: 14}\\
  \texttt{ 2408421042184839: 14}\\
  \texttt{ 2409602164565263: 14}

  \texttt{ 2409651282786775: 14}\\
  \texttt{ 2409675033576967: 14}\\
  \texttt{ 2412247712018120: 14}\\
  \texttt{ 2413722785012620: 14}\\
  \texttt{ 2414406811493716: 14}\\
  \texttt{ 2414611458380139: 14}\\
  \texttt{ 2414634420258044: 14}\\
  \texttt{ 2418970465788668: 14}\\
  \texttt{ 2419001572650171: 14}\\
  \texttt{ 2419908720643420: 14}\\
  \texttt{ 2425851496456963: 14}\\
  \texttt{ 2427606621970575: 14}\\
  \texttt{ 2428107846935019: 14}\\
  \texttt{ 2429110412249012: 14}\\
  \texttt{ 2433195473807571: 14}\\
  \texttt{ 2433824578131424: 14}\\
  \texttt{ 2435838837492967: 14}\\
  \texttt{ 2438162453360139: 14}\\
  \texttt{ 2438638279167615: 14}\\
  \texttt{ 2439985084861867: 14}

  \texttt{ 2440160469030968: 14}\\
  \texttt{ 2440829532744268: 14}\\
  \texttt{ 2441371979063175: 14}\\
  \texttt{ 2443689908334763: 14}\\
  \texttt{ 2445911852586712: 14}\\
  \texttt{ 2446173351251260: 14}\\
  \texttt{ 2448677625365140: 14}\\
  \texttt{ 2449979385759339: 14}\\
  \texttt{ 2453261241599872: 14}\\
  \texttt{ 2455366520563712: 14}\\
  \texttt{ 2456439669227571: 14}\\
  \texttt{ 2456803397790475: 14}\\
  \texttt{ 2459635732247947: 14}\\
  \texttt{ 2461157052484867: 14}\\
  \texttt{ 2463681969241819: 14}\\
  \texttt{ 2464998510231015: 14}\\
  \texttt{ 2473257193967563: 14}\\
  \texttt{ 2474311002101467: 14}\\
  \texttt{ 2475619667420768: 14}\\
  \texttt{ 2476247498049247: 14}

  \texttt{ 2477535992303739: 14}\\
  \texttt{ 2477829728922416: 14}\\
  \texttt{ 2478726465266744: 14}\\
  \texttt{ 2479726898087012: 14}\\
  \texttt{ 2481163780412672: 14}\\
  \texttt{ 2483907327598372: 14}\\
  \texttt{ 2484780146330443: 14}\\
  \texttt{ 2488458809683239: 14}\\
  \texttt{ 2492317253656916: 15}\\
  \texttt{ 2492695466950623: 14}\\
  \texttt{ 2494791005973171: 14}\\
  \texttt{ 2494973083162219: 14}\\
  \texttt{ 2499587769228712: 14}\\
  \texttt{ 2500154009923912: 14}\\
  \texttt{ 2500666190283112: 14}\\
  \texttt{ 2508082816027119: 14}\\
  \texttt{ 2511183514498323: 14}\\
  \texttt{ 2511312576215524: 14}\\
  \texttt{ 2511380682614564: 14}\\
  \texttt{ 2511730307589343: 14}

  \texttt{ 2511836591550416: 14}\\
  \texttt{ 2512437147944871: 14}\\
  \texttt{ 2515246419125312: 14}\\
  \texttt{ 2516348031055820: 14}\\
  \texttt{ 2517528451630120: 14}\\
  \texttt{ 2519750836305620: 14}\\
  \texttt{ 2524025974773248: 14}\\
  \texttt{ 2525752423108748: 14}\\
  \texttt{ 2527858408363243: 14}\\
  \texttt{ 2528737725298075: 14}\\
  \texttt{ 2529225603791719: 14}\\
  \texttt{ 2529240557714475: 14}\\
  \texttt{ 2532302183353912: 14}\\
  \texttt{ 2536827809549444: 14}\\
  \texttt{ 2537767309477516: 14}\\
  \texttt{ 2538799162230220: 14}\\
  \texttt{ 2540089624072119: 14}\\
  \texttt{ 2541568109762716: 14}\\
  \texttt{ 2542769879498440: 14}\\
  \texttt{ 2553572674105963: 14}

  \texttt{ 2555656997828048: 14}\\
  \texttt{ 2556359042331175: 14}\\
  \texttt{ 2557778592591772: 14}\\
  \texttt{ 2561593663269763: 14}\\
  \texttt{ 2567738389705144: 14}\\
  \texttt{ 2570296944367275: 14}\\
  \texttt{ 2575565057365215: 14}\\
  \texttt{ 2579202227326624: 14}\\
  \texttt{ 2580637083699724: 14}\\
  \texttt{ 2585820029736968: 14}\\
  \texttt{ 2586708125800315: 14}\\
  \texttt{ 2587119177249916: 14}\\
  \texttt{ 2590615352654367: 15}\\
  \texttt{ 2593553478124768: 14}\\
  \texttt{ 2594468609342719: 14}\\
  \texttt{ 2594745531114243: 14}\\
  \texttt{ 2595870557939320: 14}\\
  \texttt{ 2597051740620464: 14}\\
  \texttt{ 2599286446757419: 14}\\
  \texttt{ 2603041108983843: 14}

  \texttt{ 2609579851101368: 14}\\
  \texttt{ 2610793953128419: 14}\\
  \texttt{ 2611259946563667: 14}\\
  \texttt{ 2612663241349367: 14}\\
  \texttt{ 2613121379382544: 14}\\
  \texttt{ 2617633960864216: 14}\\
  \texttt{ 2618825667215920: 14}\\
  \texttt{ 2621397959691171: 14}\\
  \texttt{ 2621558324019319: 14}\\
  \texttt{ 2625690923674071: 14}\\
  \texttt{ 2631470008034612: 14}\\
  \texttt{ 2634158551803639: 14}\\
  \texttt{ 2636011455733664: 14}\\
  \texttt{ 2638621748249971: 14}\\
  \texttt{ 2641989273970264: 14}\\
  \texttt{ 2644058830499043: 14}\\
  \texttt{ 2644342472151980: 14}\\
  \texttt{ 2644991500984502: 15}\\
  \texttt{ 2645677237038164: 14}\\
  \texttt{ 2647925387919044: 14}

  \texttt{ 2651908745869244: 14}\\
  \texttt{ 2652716809728044: 14}\\
  \texttt{ 2658491914380172: 14}\\
  \texttt{ 2658570449708716: 14}\\
  \texttt{ 2660794196408368: 14}\\
  \texttt{ 2662648951696719: 14}\\
  \texttt{ 2663297203185016: 14}\\
  \texttt{ 2664422539791812: 14}\\
  \texttt{ 2666542664654720: 14}\\
  \texttt{ 2666814782261372: 14}\\
  \texttt{ 2667207894468020: 14}\\
  \texttt{ 2668788975000968: 14}\\
  \texttt{ 2670476002153612: 14}\\
  \texttt{ 2673275790743359: 14}\\
  \texttt{ 2676253947366615: 14}\\
  \texttt{ 2679147152593243: 14}\\
  \texttt{ 2682589763221964: 14}\\
  \texttt{ 2683789624068712: 14}\\
  \texttt{ 2689538000774691: 14}\\
  \texttt{ 2689966510578944: 14}

  \texttt{ 2691945282243644: 14}\\
  \texttt{ 2692776998404743: 15}\\
  \texttt{ 2693615320572668: 14}\\
  \texttt{ 2695568213779124: 14}\\
  \texttt{ 2697839600486019: 14}\\
  \texttt{ 2699644120756720: 14}\\
  \texttt{ 2699858726328471: 14}\\
  \texttt{ 2702229285212824: 15}\\
  \texttt{ 2702801090229640: 14}\\
  \texttt{ 2703272430824375: 14}\\
  \texttt{ 2703374509369515: 14}\\
  \texttt{ 2705370353669543: 14}\\
  \texttt{ 2710260123179012: 14}\\
  \texttt{ 2710659025317724: 14}\\
  \texttt{ 2713153465020968: 14}\\
  \texttt{ 2714915600362815: 14}\\
  \texttt{ 2718390927626619: 14}\\
  \texttt{ 2722591945865715: 14}\\
  \texttt{ 2722864846007212: 14}\\
  \texttt{ 2724462667891672: 14}

  \texttt{ 2724990457526440: 14}\\
  \texttt{ 2728599569145916: 14}\\
  \texttt{ 2730369162211964: 14}\\
  \texttt{ 2731162364046764: 14}\\
  \texttt{ 2733410445969067: 14}\\
  \texttt{ 2734122664202524: 14}\\
  \texttt{ 2737178363940775: 14}\\
  \texttt{ 2737504083846463: 14}\\
  \texttt{ 2740587694413667: 14}\\
  \texttt{ 2741389026472915: 14}\\
  \texttt{ 2744854868036823: 14}\\
  \texttt{ 2747728290835275: 14}\\
  \texttt{ 2748005088413416: 14}\\
  \texttt{ 2748636635984467: 14}\\
  \texttt{ 2750863844195863: 14}\\
  \texttt{ 2751852964985924: 14}\\
  \texttt{ 2752628301384344: 14}\\
  \texttt{ 2752760081188923: 14}\\
  \texttt{ 2755703685381424: 14}\\
  \texttt{ 2755777432133315: 14}

  \texttt{ 2757924453127875: 14}\\
  \texttt{ 2758717760062975: 14}\\
  \texttt{ 2760906765445244: 14}\\
  \texttt{ 2764063795991524: 14}\\
  \texttt{ 2764165497851563: 14}\\
  \texttt{ 2766367169136620: 14}\\
  \texttt{ 2769589031340544: 14}\\
  \texttt{ 2770665092226115: 14}\\
  \texttt{ 2771342322140347: 14}\\
  \texttt{ 2771747111215815: 14}\\
  \texttt{ 2774059093456820: 14}\\
  \texttt{ 2775812183243971: 14}\\
  \texttt{ 2777372662869115: 14}\\
  \texttt{ 2784965737303275: 14}\\
  \texttt{ 2786307189831763: 14}\\
  \texttt{ 2786399556697772: 14}\\
  \texttt{ 2786729712749863: 14}\\
  \texttt{ 2788685982968012: 14}\\
  \texttt{ 2789307338908875: 14}\\
  \texttt{ 2789876651186175: 14}

  \texttt{ 2790407247168975: 14}\\
  \texttt{ 2790510426428012: 14}\\
  \texttt{ 2791756575403071: 14}\\
  \texttt{ 2798386571336415: 14}\\
  \texttt{ 2799118579629463: 14}\\
  \texttt{ 2799954790395423: 14}\\
  \texttt{ 2800921479542043: 14}\\
  \texttt{ 2802023474882515: 14}\\
  \texttt{ 2802643150048072: 14}\\
  \texttt{ 2802763793269412: 14}\\
  \texttt{ 2804868315347740: 14}\\
  \texttt{ 2805157935791440: 14}\\
  \texttt{ 2806242713969775: 14}\\
  \texttt{ 2807539417835468: 14}\\
  \texttt{ 2809430736346815: 14}\\
  \texttt{ 2810527442070640: 14}\\
  \texttt{ 2811892375294171: 14}\\
  \texttt{ 2813269473765043: 14}\\
  \texttt{ 2815017654663944: 14}\\
  \texttt{ 2816648298458943: 14}

  \texttt{ 2820926323469671: 14}\\
  \texttt{ 2823012620135108: 14}\\
  \texttt{ 2826123914789871: 14}\\
  \texttt{ 2829883248928843: 14}\\
  \texttt{ 2829939122180748: 14}\\
  \texttt{ 2832442820527864: 15}\\
  \texttt{ 2834319755090672: 14}\\
  \texttt{ 2836136873639744: 14}\\
  \texttt{ 2837473623058075: 14}\\
  \texttt{ 2841689545574444: 14}\\
  \texttt{ 2843352365356943: 14}\\
  \texttt{ 2844943274028974: 15}\\
  \texttt{ 2845290804581175: 14}\\
  \texttt{ 2847535506542419: 14}\\
  \texttt{ 2850043151150447: 14}\\
  \texttt{ 2850916035708619: 14}\\
  \texttt{ 2851092223848548: 14}\\
  \texttt{ 2852206214101244: 14}\\
  \texttt{ 2852697670270315: 14}\\
  \texttt{ 2857206549797647: 14}

  \texttt{ 2857976980783720: 14}\\
  \texttt{ 2860464035277368: 14}\\
  \texttt{ 2863980886890175: 14}\\
  \texttt{ 2864243655041343: 14}\\
  \texttt{ 2864451352768324: 14}\\
  \texttt{ 2864921232656347: 14}\\
  \texttt{ 2867203377560348: 14}\\
  \texttt{ 2871889476266672: 14}\\
  \texttt{ 2874661037073044: 14}\\
  \texttt{ 2876506355409163: 14}\\
  \texttt{ 2878154923094368: 14}\\
  \texttt{ 2878519156604943: 14}\\
  \texttt{ 2887235379180268: 14}\\
  \texttt{ 2887442431737568: 14}\\
  \texttt{ 2891469277369275: 14}\\
  \texttt{ 2892343649797420: 14}\\
  \texttt{ 2893454393402275: 14}\\
  \texttt{ 2893768196727572: 14}\\
  \texttt{ 2895812426398867: 14}\\
  \texttt{ 2895914657869747: 14}

  \texttt{ 2897411769342871: 14}\\
  \texttt{ 2900566363967439: 14}\\
  \texttt{ 2902471702741663: 14}\\
  \texttt{ 2900566363967439: 14}\\
  \texttt{ 2902471702741663: 14}\\
  \texttt{ 2904130500284139: 14}\\
  \texttt{ 2904788136137912: 14}\\
  \texttt{ 2907405665395963: 14}\\
  \texttt{ 2908951737258572: 14}\\
  \texttt{ 2909170002517244: 14}\\
  \texttt{ 2909935116065140: 14}\\
  \texttt{ 2909961461164446: 15}\\
  \texttt{ 2910198498999415: 14}\\
  \texttt{ 2911047315989012: 14}\\
  \texttt{ 2911089562877168: 14}\\
  \texttt{ 2911237702007947: 14}\\
  \texttt{ 2911390529486523: 14}\\
  \texttt{ 2912145822645020: 14}\\
  \texttt{ 2913111930559912: 15}\\
  \texttt{ 2914687739614112: 14}

  \texttt{ 2914894892202667: 14}\\
  \texttt{ 2915331463894372: 14}\\
  \texttt{ 2915462708295243: 14}\\
  \texttt{ 2917807137051423: 14}\\
  \texttt{ 2918164202548467: 14}\\
  \texttt{ 2923636010373016: 14}\\
  \texttt{ 2925339649292068: 14}\\
  \texttt{ 2929717586624372: 14}\\
  \texttt{ 2929889523299343: 14}\\
  \texttt{ 2930334878629843: 14}\\
  \texttt{ 2931586825376643: 14}\\
  \texttt{ 2931869512713463: 14}\\
  \texttt{ 2936807265402368: 15}\\
  \texttt{ 2937071144945672: 14}\\
  \texttt{ 2938419444863320: 14}\\
  \texttt{ 2938998334596639: 14}\\
  \texttt{ 2939235802880919: 14}\\
  \texttt{ 2939323755810075: 14}\\
  \texttt{ 2939458907662612: 14}\\
  \texttt{ 2944626903462867: 14}

  \texttt{ 2945074150470319: 14}\\
  \texttt{ 2947566103352647: 14}\\
  \texttt{ 2958383424295447: 14}\\
  \texttt{ 2959655754403664: 14}\\
  \texttt{ 2959892893442347: 14}\\
  \texttt{ 2959942263584047: 14}\\
  \texttt{ 2959944166955672: 14}\\
  \texttt{ 2960866051123122: 15}\\
  \texttt{ 2963362128160372: 14}\\
  \texttt{ 2963532135651343: 14}\\
  \texttt{ 2964703068067767: 14}\\
  \texttt{ 2966275358199844: 14}\\
  \texttt{ 2967797082281048: 14}\\
  \texttt{ 2969263560447115: 14}\\
  \texttt{ 2971036963029724: 14}\\
  \texttt{ 2971518185490819: 14}\\
  \texttt{ 2972532423738571: 14}\\
  \texttt{ 2976990118229715: 14}\\
  \texttt{ 2977640320461316: 14}\\
  \texttt{ 2979488195414224: 14}

  \texttt{ 2981870327479075: 14}\\
  \texttt{ 2983724975411144: 14}\\
  \texttt{ 2986513578401812: 14}\\
  \texttt{ 2988047326991344: 14}\\
  \texttt{ 2988814297511744: 14}\\
  \texttt{ 2990112016603543: 14}\\
  \texttt{ 2991920378291864: 14}\\
  \texttt{ 2996746935936639: 14}\\
  \texttt{ 2990112016603543: 14}\\
  \texttt{ 2991920378291864: 14}\\
  \texttt{ 2996746935936639: 14}\\
  \texttt{ 2999629540747240: 14}\\
  \texttt{ 3001121100478519: 14}\\
  \texttt{ 3004439580454624: 14}\\
  \texttt{ 3006952767069128: 14}\\
  \texttt{ 3008030835770863: 14}\\
  \texttt{ 3008880384513172: 14}\\
  \texttt{ 3009379507927468: 14}\\
  \texttt{ 3009534553016115: 14}\\
  \texttt{ 3009594651846668: 14}

  \texttt{ 3011084396685520: 14}\\
  \texttt{ 3011877239085043: 14}\\
  \texttt{ 3014753614204975: 14}\\
  \texttt{ 3016785100847020: 14}\\
  \texttt{ 3019513500621368: 14}\\
  \texttt{ 3020554352780072: 14}\\
  \texttt{ 3020945589562839: 14}\\
  \texttt{ 3021380027447020: 14}\\
  \texttt{ 3028116336015939: 14}\\
  \texttt{ 3031000180856944: 14}\\
  \texttt{ 3033355805679472: 14}\\
  \texttt{ 3034091779290944: 14}\\
  \texttt{ 3034332515503419: 14}\\
  \texttt{ 3035065451166546: 15}\\
  \texttt{ 3035737754832220: 14}\\
  \texttt{ 3037727891767472: 14}\\
  \texttt{ 3038449726298068: 14}\\
  \texttt{ 3038618936751567: 14}\\
  \texttt{ 3040563252275224: 14}\\
  \texttt{ 3040860176563616: 14}

  \texttt{ 3041786812959512: 14}\\
  \texttt{ 3049011602916019: 14}\\
  \texttt{ 3049477928756143: 14}\\
  \texttt{ 3053031831732967: 14}\\
  \texttt{ 3053740655136219: 14}\\
  \texttt{ 3054685594568416: 14}\\
  \texttt{ 3058317133643167: 14}\\
  \texttt{ 3059969146941275: 14}\\
  \texttt{ 3060784644109172: 14}\\
  \texttt{ 3064867491760815: 14}\\
  \texttt{ 3065442224714944: 14}\\
  \texttt{ 3067584676291923: 14}\\
  \texttt{ 3074034112936612: 14}\\
  \texttt{ 3079218306084848: 14}\\
  \texttt{ 3080939921736343: 14}\\
  \texttt{ 3081138953412219: 14}\\
  \texttt{ 3082840973371599: 14}\\
  \texttt{ 3088816098281071: 14}\\
  \texttt{ 3092151857330419: 14}\\
  \texttt{ 3092459635183444: 14}

  \texttt{ 3093232947312068: 14}\\
  \texttt{ 3093649058821840: 14}\\
  \texttt{ 3098534992603719: 14}\\
  \texttt{ 3099747468395312: 14}\\
  \texttt{ 3100696480441924: 14}\\
  \texttt{ 3104374913717768: 14}\\
  \texttt{ 3108498999245647: 14}\\
  \texttt{ 3109714735624540: 14}\\
  \texttt{ 3110593252313275: 14}\\
  \texttt{ 3115938999733540: 14}\\
  \texttt{ 3116767542071443: 14}\\
  \texttt{ 3121259445252547: 14}\\
  \texttt{ 3121504722126267: 14}\\
  \texttt{ 3128256852700416: 14}\\
  \texttt{ 3130381562538220: 14}\\
  \texttt{ 3131870650436716: 14}\\
  \texttt{ 3131931238587339: 14}\\
  \texttt{ 3133439263161124: 14}\\
  \texttt{ 3135873466148740: 14}\\
  \texttt{ 3138068460856472: 14}

  \texttt{ 3139076462751044: 14}\\
  \texttt{ 3139076462751044: 14}\\
  \texttt{ 3140577772600867: 14}\\
  \texttt{ 3140878571432943: 14}\\
  \texttt{ 3143414172046444: 14}\\
  \texttt{ 3143435170285743: 14}\\
  \texttt{ 3147538341610744: 14}\\
  \texttt{ 3149774468031067: 14}\\
  \texttt{ 3149938069890019: 14}\\
  \texttt{ 3150210999594940: 14}\\
  \texttt{ 3156119612382215: 14}\\
  \texttt{ 3159073337100668: 14}\\
  \texttt{ 3159525430406419: 14}\\
  \texttt{ 3161260768558795: 14}\\
  \texttt{ 3162558900413312: 14}\\
  \texttt{ 3163072686881920: 14}\\
  \texttt{ 3167242115845119: 14}\\
  \texttt{ 3168359769426415: 14}\\
  \texttt{ 3168499154585919: 14}\\
  \texttt{ 3168841929541316: 14}

  \texttt{ 3169819866038812: 14}\\
  \texttt{ 3171195545056371: 14}\\
  \texttt{ 3171516208438120: 14}\\
  \texttt{ 3171583108535264: 14}\\
  \texttt{ 3171583108535264: 14}\\
  \texttt{ 3173300123265871: 14}\\
  \texttt{ 3173305330211812: 14}\\
  \texttt{ 3175823171725720: 14}\\
  \texttt{ 3178593518861715: 14}\\
  \texttt{ 3179391992963144: 14}\\
  \texttt{ 3182265925005571: 14}\\
  \texttt{ 3184447249706668: 14}\\
  \texttt{ 3187978441849312: 14}\\
  \texttt{ 3188451737037039: 14}\\
  \texttt{ 3188645572390747: 14}\\
  \texttt{ 3189356303747775: 14}\\
  \texttt{ 3192696854754775: 14}\\
  \texttt{ 3192848572997312: 14}\\
  \texttt{ 3194958259787120: 14}\\
  \texttt{ 3195372586921963: 14}

  \texttt{ 3197646964352572: 14}\\
  \texttt{ 3198267768119740: 14}\\
  \texttt{ 3199119738608563: 14}\\
  \texttt{ 3199153368836572: 14}\\
  \texttt{ 3199433724539368: 14}\\
  \texttt{ 3199658816229018: 15}\\
  \texttt{ 3204152298207963: 14}\\
  \texttt{ 3204743122072420: 14}\\
  \texttt{ 3205782956790619: 14}\\
  \texttt{ 3206957946195175: 14}\\
  \texttt{ 3207051202230043: 14}\\
  \texttt{ 3207304336013944: 14}\\
  \texttt{ 3212805940336311: 15}\\
  \texttt{ 3213097301495211: 15}\\
  \texttt{ 3215226335143218: 16}\\
  \texttt{ 3219554707141219: 14}\\
  \texttt{ 3221091577788243: 14}\\
  \texttt{ 3223260596917520: 14}\\
  \texttt{ 3225661593540343: 14}\\
  \texttt{ 3226237856529919: 14}

  \texttt{ 3227766236184643: 14}\\
  \texttt{ 3227919643747267: 14}\\
  \texttt{ 3237836828174143: 14}\\
  \texttt{ 3238484447117564: 14}\\
  \texttt{ 3241353846079420: 14}\\
  \texttt{ 3242096083478216: 14}\\
  \texttt{ 3243163631424820: 14}\\
  \texttt{ 3244626803678875: 14}\\
  \texttt{ 3245802312052268: 14}\\
  \texttt{ 3245914385004075: 14}\\
  \texttt{ 3248031537275120: 14}\\
  \texttt{ 3249448167443019: 14}\\
  \texttt{ 3251352294683223: 14}\\
  \texttt{ 3252912278861320: 14}\\
  \texttt{ 3261399485476112: 14}\\
  \texttt{ 3261786051579020: 14}\\
  \texttt{ 3270819242366620: 14}\\
  \texttt{ 3274061627441644: 14}\\
  \texttt{ 3277263843303664: 15}\\
  \texttt{ 3279170169743919: 14}

  \texttt{ 3281413877959012: 14}\\
  \texttt{ 3282193324090015: 14}\\
  \texttt{ 3284326533880843: 14}\\
  \texttt{ 3286881346289775: 14}\\
  \texttt{ 3288545722763115: 14}\\
  \texttt{ 3291933254260112: 14}\\
  \texttt{ 3295151573826039: 14}\\
  \texttt{ 3295905730121715: 14}\\
  \texttt{ 3298191933079743: 14}\\
  \texttt{ 3299255557509319: 14}\\
  \texttt{ 3299902103580244: 14}\\
  \texttt{ 3305879581285520: 14}\\
  \texttt{ 3307095108021643: 14}\\
  \texttt{ 3307829195721715: 14}\\
  \texttt{ 3312504617745464: 14}\\
  \texttt{ 3315371283928816: 14}\\
  \texttt{ 3318367696533375: 14}\\
  \texttt{ 3320412587654120: 14}\\
  \texttt{ 3320563792818544: 15}\\
  \texttt{ 3321109866856472: 14}

  \texttt{ 3323301229656519: 14}\\
  \texttt{ 3324460376737743: 14}\\
  \texttt{ 3325183465082343: 14}\\
  \texttt{ 3325234444088667: 14}\\
  \texttt{ 3325962666679971: 14}\\
  \texttt{ 3326757143066619: 14}\\
  \texttt{ 3330073683916568: 14}\\
  \texttt{ 3331200704003440: 14}\\
  \texttt{ 3340769917339015: 14}\\
  \texttt{ 3343292333229248: 14}\\
  \texttt{ 3343295076307963: 14}\\
  \texttt{ 3345990317771920: 14}\\
  \texttt{ 3346157812769863: 14}\\
  \texttt{ 3348968726551471: 14}\\
  \texttt{ 3351236732834515: 14}\\
  \texttt{ 3353964267359815: 14}\\
  \texttt{ 3358106557882675: 14}\\
  \texttt{ 3358113664083568: 14}\\
  \texttt{ 3358169243372815: 14}\\
  \texttt{ 3358657022096816: 14}

  \texttt{ 3360151989824371: 14}\\
  \texttt{ 3361365487657240: 14}\\
  \texttt{ 3364641260010464: 14}\\
  \texttt{ 3365957619792123: 14}\\
  \texttt{ 3366879152279371: 14}\\
  \texttt{ 3370057043530868: 14}\\
  \texttt{ 3370165953189819: 14}\\
  \texttt{ 3373949254915671: 14}\\
  \texttt{ 3374921357265796: 14}\\
  \texttt{ 3378742930801843: 14}\\
  \texttt{ 3379740716153968: 14}\\
  \texttt{ 3379938295779116: 14}\\
  \texttt{ 3381275073267368: 14}\\
  \texttt{ 3383078588047671: 14}\\
  \texttt{ 3384893264632216: 14}\\
  \texttt{ 3385099150177767: 14}\\
  \texttt{ 3386257822986172: 14}\\
  \texttt{ 3386624332631343: 14}\\
  \texttt{ 3392365126096239: 14}\\
  \texttt{ 3393335565929824: 14}

  \texttt{ 3394204602769275: 14}\\
  \texttt{ 3398555457858464: 14}\\
  \texttt{ 3401109619293940: 14}\\
  \texttt{ 3404888784590019: 14}\\
  \texttt{ 3405754212106912: 14}\\
  \texttt{ 3406429676242744: 14}\\
  \texttt{ 3406651841709044: 14}\\
  \texttt{ 3407196263917912: 14}\\
  \texttt{ 3410026966313223: 14}\\
  \texttt{ 3410881694413215: 14}\\
  \texttt{ 3412532160815416: 14}\\
  \texttt{ 3413950906373439: 14}\\
  \texttt{ 3414906230439572: 14}\\
  \texttt{ 3416946798097015: 14}\\
  \texttt{ 3417721620497363: 14}\\
  \texttt{ 3417926083198675: 14}\\
  \texttt{ 3418009019023518: 15}\\
  \texttt{ 3422591417989671: 14}\\
  \texttt{ 3424567424035868: 14}\\
  \texttt{ 3424768586876968: 14}

  \texttt{ 3424818484382524: 14}\\
  \texttt{ 3428958266104119: 14}\\
  \texttt{ 3433319730966723: 14}\\
  \texttt{ 3434237672464375: 14}\\
  \texttt{ 3435858598541012: 14}\\
  \texttt{ 3442735798785512: 14}\\
  \texttt{ 3445355224583863: 14}\\
  \texttt{ 3445460955882944: 14}\\
  \texttt{ 3446292250040120: 14}\\
  \texttt{ 3447280473741619: 14}\\
  \texttt{ 3452503477650019: 14}\\
  \texttt{ 3452576341626520: 14}\\
  \texttt{ 3459383742451840: 14}\\
  \texttt{ 3460363358703640: 15}\\
  \texttt{ 3461406657698419: 14}\\
  \texttt{ 3461931554112115: 14}\\
  \texttt{ 3464331466867964: 14}\\
  \texttt{ 3464841053364639: 14}\\
  \texttt{ 3468710031782524: 14}\\
  \texttt{ 3468985246308543: 14}

  \texttt{ 3469669739850967: 14}\\
  \texttt{ 3470411951697040: 14}\\
  \texttt{ 3470624611245944: 14}\\
  \texttt{ 3470880965218324: 14}\\
  \texttt{ 3471966699181739: 14}\\
  \texttt{ 3473348746176771: 14}\\
  \texttt{ 3475063189294712: 14}\\
  \texttt{ 3475450093120468: 14}\\
  \texttt{ 3478178682085975: 14}\\
  \texttt{ 3478643961438724: 14}\\
  \texttt{ 3479343466096275: 14}\\
  \texttt{ 3480544396159719: 14}\\
  \texttt{ 3480702313680522: 15}\\
  \texttt{ 3481649939500720: 14}\\
  \texttt{ 3482023737272372: 14}\\
  \texttt{ 3483145838636475: 14}\\
  \texttt{ 3486548325744367: 14}\\
  \texttt{ 3487822146384919: 14}\\
  \texttt{ 3491912987666719: 14}\\
  \texttt{ 3493756169857420: 14}

  \texttt{ 3495721639589043: 14}\\
  \texttt{ 3496129564588464: 14}\\
  \texttt{ 3496588122484743: 14}\\
  \texttt{ 3500474050440340: 14}\\
  \texttt{ 3502595543448320: 14}\\
  \texttt{ 3509290176124219: 14}\\
  \texttt{ 3511503508311772: 14}\\
  \texttt{ 3512025335803912: 14}\\
  \texttt{ 3514022591042715: 14}\\
  \texttt{ 3526673716637143: 14}\\
  \texttt{ 3527448427439516: 14}\\
  \texttt{ 3528916269206920: 14}\\
  \texttt{ 3530531458245464: 14}\\
  \texttt{ 3531176253491272: 14}\\
  \texttt{ 3531223260457964: 14}\\
  \texttt{ 3533088420781575: 14}\\
  \texttt{ 3533457580022312: 14}\\
  \texttt{ 3534050567526416: 14}\\
  \texttt{ 3535041324799923: 14}\\
  \texttt{ 3535424978917840: 14}

  \texttt{ 3535540527585868: 14}\\
  \texttt{ 3537544826296743: 14}\\
  \texttt{ 3538503591424035: 14}\\
  \texttt{ 3542491326156712: 14}\\
  \texttt{ 3549042646712115: 14}\\
  \texttt{ 3551659762203771: 14}\\
  \texttt{ 3556164324013316: 14}\\
  \texttt{ 3556517590866607: 14}\\
  \texttt{ 3558748762728724: 15}\\
  \texttt{ 3559102084995171: 14}\\
  \texttt{ 3559405368481467: 14}\\
  \texttt{ 3560222619428443: 14}\\
  \texttt{ 3560485476483519: 14}\\
  \texttt{ 3562447292251820: 14}\\
  \texttt{ 3563885364006068: 14}\\
  \texttt{ 3564078833275244: 14}\\
  \texttt{ 3564192185623312: 14}\\
  \texttt{ 3566350769847640: 14}\\
  \texttt{ 3566901852552039: 14}\\
  \texttt{ 3568584694008543: 14}

  \texttt{ 3569915213064519: 14}\\
  \texttt{ 3570935569653319: 14}\\
  \texttt{ 3571412319184647: 14}\\
  \texttt{ 3572578027922444: 14}\\
  \texttt{ 3573503771782867: 14}\\
  \texttt{ 3576074237250820: 14}\\
  \texttt{ 3577008167898640: 14}\\
  \texttt{ 3579815947845471: 14}\\
  \texttt{ 3582981863776707: 14}\\
  \texttt{ 3583268907965972: 14}\\
  \texttt{ 3584077896659443: 14}\\
  \texttt{ 3584157346962774: 15}\\
  \texttt{ 3584885346039472: 14}\\
  \texttt{ 3586994439182115: 14}\\
  \texttt{ 3589259182889048: 14}\\
  \texttt{ 3590863980692019: 14}\\
  \texttt{ 3591248687841067: 14}\\
  \texttt{ 3594180583011544: 14}\\
  \texttt{ 3596921226190312: 14}\\
  \texttt{ 3596988422420072: 14}

  \texttt{ 3598491337267012: 14}\\
  \texttt{ 3599851014538172: 14}\\
  \texttt{ 3599896191802388: 14}\\
  \texttt{ 3600170791252971: 14}\\
  \texttt{ 3600634511741667: 14}\\
  \texttt{ 3602337490895515: 14}\\
  \texttt{ 3606030837007171: 14}\\
  \texttt{ 3606241828891112: 14}\\
  \texttt{ 3611414211587620: 14}\\
  \texttt{ 3612059998718712: 14}\\
  \texttt{ 3612079536916867: 14}\\
  \texttt{ 3612783373061912: 14}\\
  \texttt{ 3614042721967243: 14}\\
  \texttt{ 3614610906378764: 14}\\
  \texttt{ 3614808857400340: 14}\\
  \texttt{ 3615675479168415: 14}\\
  \texttt{ 3615686396898640: 14}\\
  \texttt{ 3617416318242220: 14}\\
  \texttt{ 3617571421600712: 14}\\
  \texttt{ 3618254640643243: 14}

  \texttt{ 3619750044445063: 14}\\
  \texttt{ 3620669860760271: 14}\\
  \texttt{ 3625110232057140: 14}\\
  \texttt{ 3626289253645712: 14}\\
  \texttt{ 3632356483819316: 14}\\
  \texttt{ 3635517178737843: 14}\\
  \texttt{ 3636245791632616: 14}\\
  \texttt{ 3636250729963467: 14}\\
  \texttt{ 3639459188438619: 14}\\
  \texttt{ 3641070195132244: 14}\\
  \texttt{ 3642744468795667: 14}\\
  \texttt{ 3649298399025712: 14}\\
  \texttt{ 3649831428505464: 14}\\
  \texttt{ 3651749924229944: 14}\\
  \texttt{ 3658766485879240: 14}\\
  \texttt{ 3661164574239464: 14}\\
  \texttt{ 3664083965478668: 14}\\
  \texttt{ 3664136530483072: 14}\\
  \texttt{ 3664324250926568: 14}\\
  \texttt{ 3665823062184116: 14}

  \texttt{ 3666992366516444: 14}\\
  \texttt{ 3667996413943819: 14}\\
  \texttt{ 3668897808223568: 14}\\
  \texttt{ 3668975963663416: 14}\\
  \texttt{ 3671197745423739: 14}\\
  \texttt{ 3673984494594212: 14}\\
  \texttt{ 3675427705073275: 14}\\
  \texttt{ 3682853397296043: 14}\\
  \texttt{ 3683096408939047: 14}\\
  \texttt{ 3683753258305543: 14}\\
  \texttt{ 3684072223968712: 14}\\
  \texttt{ 3684960340567071: 14}\\
  \texttt{ 3684976294651612: 14}\\
  \texttt{ 3685012627619716: 14}\\
  \texttt{ 3686499603597543: 14}\\
  \texttt{ 3686592535822131: 14}\\
  \texttt{ 3689017823664916: 15}\\
  \texttt{ 3692009382385867: 14}\\
  \texttt{ 3692701222549171: 14}\\
  \texttt{ 3694144197281740: 14}

  \texttt{ 3694836267399975: 14}\\
  \texttt{ 3697588409205572: 14}\\
  \texttt{ 3700715885989064: 14}\\
  \texttt{ 3701101642181612: 14}\\
  \texttt{ 3704032813186924: 14}\\
  \texttt{ 3705765261011223: 14}\\
  \texttt{ 3705918417055375: 14}\\
  \texttt{ 3713708138505820: 14}\\
  \texttt{ 3715050705180675: 14}\\
  \texttt{ 3715194872515372: 14}\\
  \texttt{ 3715818780200919: 14}\\
  \texttt{ 3718904052512168: 14}\\
  \texttt{ 3719898725024140: 14}\\
  \texttt{ 3720371874324568: 14}\\
  \texttt{ 3723020688774716: 14}\\
  \texttt{ 3723378465971175: 14}\\
  \texttt{ 3723387862706467: 14}\\
  \texttt{ 3726309352020172: 14}\\
  \texttt{ 3726319657132664: 14}\\
  \texttt{ 3726597368010472: 14}

  \texttt{ 3728829523405543: 14}\\
  \texttt{ 3728991750334244: 14}\\
  \texttt{ 3729529890895616: 14}\\
  \texttt{ 3736091074982739: 14}\\
  \texttt{ 3737976584203771: 14}\\
  \texttt{ 3738202791259443: 14}\\
  \texttt{ 3739202605376116: 14}\\
  \texttt{ 3739309665446716: 14}\\
  \texttt{ 3741087049066420: 14}\\
  \texttt{ 3741310343778919: 14}\\
  \texttt{ 3742600406517916: 14}\\
  \texttt{ 3742923073585623: 14}\\
  \texttt{ 3743091275314875: 14}\\
  \texttt{ 3743644706657871: 14}\\
  \texttt{ 3746333035131720: 14}\\
  \texttt{ 3746747699214171: 14}\\
  \texttt{ 3746891575344463: 14}\\
  \texttt{ 3751014691621671: 15}\\
  \texttt{ 3751251502649443: 14}\\
  \texttt{ 3751360495919152: 14}

  \texttt{ 3751444078114616: 14}\\
  \texttt{ 3751858856541464: 14}\\
  \texttt{ 3753299515818039: 14}\\
  \texttt{ 3754754965954672: 14}\\
  \texttt{ 3755498889265539: 14}\\
  \texttt{ 3755504514040316: 15}\\
  \texttt{ 3756405470452216: 14}\\
  \texttt{ 3759083703916216: 14}\\
  \texttt{ 3759368653846916: 14}\\
  \texttt{ 3760191221335747: 14}\\
  \texttt{ 3760436862221263: 14}\\
  \texttt{ 3762583070205615: 14}\\
  \texttt{ 3766172893141443: 14}\\
  \texttt{ 3767507726267967: 14}\\
  \texttt{ 3768991957201444: 14}\\
  \texttt{ 3778108132308020: 14}\\
  \texttt{ 3778186756553872: 14}\\
  \texttt{ 3778980737886123: 14}\\
  \texttt{ 3779022638195071: 14}\\
  \texttt{ 3781986828343964: 14}

  \texttt{ 3784900143362912: 14}\\
  \texttt{ 3785823293774344: 14}\\
  \texttt{ 3787417996815571: 14}\\
  \texttt{ 3789653814633615: 14}\\
  \texttt{ 3792467066875747: 14}\\
  \texttt{ 3792811044784616: 14}\\
  \texttt{ 3795604864187371: 14}\\
  \texttt{ 3796267726981840: 15}\\
  \texttt{ 3796550724177115: 14}\\
  \texttt{ 3797863276514668: 14}\\
  \texttt{ 3798590647256223: 14}\\
  \texttt{ 3798692256861016: 14}\\
  \texttt{ 3802006003073443: 14}\\
  \texttt{ 3803711037852212: 14}\\
  \texttt{ 3804326459275867: 14}\\
  \texttt{ 3804896807088244: 14}\\
  \texttt{ 3807132040914775: 14}\\
  \texttt{ 3808442809549971: 14}\\
  \texttt{ 3809320130226075: 14}\\
  \texttt{ 3810443319111715: 14}

  \texttt{ 3813840502177075: 14}\\
  \texttt{ 3814581847139944: 14}\\
  \texttt{ 3815016929827071: 14}\\
  \texttt{ 3815987415083467: 14}\\
  \texttt{ 3816248065653616: 14}\\
  \texttt{ 3816265002229268: 14}\\
  \texttt{ 3820308952188716: 14}\\
  \texttt{ 3824536557687424: 14}\\
  \texttt{ 3825843382173212: 14}\\
  \texttt{ 3831288814546144: 14}\\
  \texttt{ 3835081347460243: 14}\\
  \texttt{ 3836829700996143: 14}\\
  \texttt{ 3837040152131319: 14}\\
  \texttt{ 3842756788859043: 14}\\
  \texttt{ 3842762640842320: 14}\\
  \texttt{ 3843357885642248: 14}\\
  \texttt{ 3844235251192652: 14}\\
  \texttt{ 3848836239718064: 14}\\
  \texttt{ 3856238371478571: 14}\\
  \texttt{ 3856714184052772: 14}

  \texttt{ 3859968783398516: 14}\\
  \texttt{ 3863713004130067: 14}\\
  \texttt{ 3866665652310248: 14}\\
  \texttt{ 3867721616102667: 14}\\
  \texttt{ 3868686520979119: 14}\\
  \texttt{ 3870556876108662: 15}\\
  \texttt{ 3875487352895823: 14}\\
  \texttt{ 3875891236108820: 14}\\
  \texttt{ 3878364702828571: 14}\\
  \texttt{ 3880694131441172: 14}\\
  \texttt{ 3882342154079068: 14}\\
  \texttt{ 3883370451529219: 14}\\
  \texttt{ 3883764066178543: 14}\\
  \texttt{ 3885604725817024: 14}\\
  \texttt{ 3888952292283639: 14}\\
  \texttt{ 3889861917379168: 15}\\
  \texttt{ 3891021996852871: 14}\\
  \texttt{ 3891546164876812: 14}\\
  \texttt{ 3893357001288716: 14}\\
  \texttt{ 3894277806887344: 14}

  \texttt{ 3895916598344715: 14}\\
  \texttt{ 3898381537956867: 14}\\
  \texttt{ 3898915622028339: 14}\\
  \texttt{ 3898971358941020: 15}\\
  \texttt{ 3899680111340620: 14}\\
  \texttt{ 3903708015649663: 14}\\
  \texttt{ 3904881140065371: 14}\\
  \texttt{ 3905010893844124: 14}\\
  \texttt{ 3905933165587412: 14}\\
  \texttt{ 3906589422273271: 14}\\
  \texttt{ 3908065156181144: 14}\\
  \texttt{ 3910931002856212: 14}\\
  \texttt{ 3911761554147547: 14}\\
  \texttt{ 3912039037425616: 14}\\
  \texttt{ 3913792712970820: 14}\\
  \texttt{ 3914729481319612: 14}\\
  \texttt{ 3917070228450168: 14}\\
  \texttt{ 3917182870313224: 14}\\
  \texttt{ 3919127314455763: 14}\\
  \texttt{ 3922304919994924: 14}

  \texttt{ 3922433736496675: 14}\\
  \texttt{ 3923773346700223: 14}\\
  \texttt{ 3924747478581163: 14}\\
  \texttt{ 3924988332075915: 14}\\
  \texttt{ 3925162024857764: 14}\\
  \texttt{ 3926267751909867: 14}\\
  \texttt{ 3926809139581024: 14}\\
  \texttt{ 3926926482008264: 14}\\
  \texttt{ 3933347303637512: 14}\\
  \texttt{ 3934806017553015: 14}\\
  \texttt{ 3934846874334247: 14}\\
  \texttt{ 3935778155217044: 14}\\
  \texttt{ 3936343349167819: 14}\\
  \texttt{ 3940907532272944: 14}\\
  \texttt{ 3942344458957372: 14}\\
  \texttt{ 3943150245206144: 14}\\
  \texttt{ 3944366221447239: 14}\\
  \texttt{ 3947746206526064: 14}\\
  \texttt{ 3951587863114064: 14}\\
  \texttt{ 3952925185990112: 14}

  \texttt{ 3953913720492688: 14}\\
  \texttt{ 3956525354472244: 14}\\
  \texttt{ 3959677598576812: 14}\\
  \texttt{ 3962831318687643: 14}\\
  \texttt{ 3964392152623312: 14}\\
  \texttt{ 3964665665999968: 14}\\
  \texttt{ 3967433732320324: 14}\\
  \texttt{ 3970522426978419: 14}\\
  \texttt{ 3972650742148147: 14}\\
  \texttt{ 3973337966965372: 14}\\
  \texttt{ 3974452352205748: 14}\\
  \texttt{ 3974807049722620: 14}\\
  \texttt{ 3976560457762120: 14}\\
  \texttt{ 3976686135625615: 14}\\
  \texttt{ 3980306586341319: 14}\\
  \texttt{ 3980434944139567: 14}\\
  \texttt{ 3981811574202244: 14}\\
  \texttt{ 3981859807176872: 14}\\
  \texttt{ 3984193595863275: 14}\\
  \texttt{ 3984801455833072: 14}

  \texttt{ 3986426214831068: 14}\\
  \texttt{ 3992052919269643: 14}\\
  \texttt{ 3992686646913916: 14}\\
  \texttt{ 3996094649254539: 14}\\
  \texttt{ 3996194726215672: 14}\\
  \texttt{ 3996261317235572: 14}\\
  \texttt{ 3998783654604339: 14}\\
  \texttt{ 4001368682033864: 14}\\
  \texttt{ 4002482524164171: 14}\\
  \texttt{ 4005072589122338: 15}\\
  \texttt{ 4005678100006843: 14}\\
  \texttt{ 4006078354174840: 14}\\
  \texttt{ 4006085283251216: 14}\\
  \texttt{ 4009626499077475: 14}\\
  \texttt{ 4009686803434267: 14}\\
  \texttt{ 4014140209990820: 14}\\
  \texttt{ 4017103519974567: 14}\\
  \texttt{ 4017379227480475: 14}\\
  \texttt{ 4017403791798848: 14}\\
  \texttt{ 4019261587159120: 14}

  \texttt{ 4019350417152415: 14}\\
  \texttt{ 4019855602594839: 14}\\
  \texttt{ 4023216977182748: 14}\\
  \texttt{ 4024893000126472: 14}\\
  \texttt{ 4032627811954567: 14}\\
  \texttt{ 4034733371084816: 14}\\
  \texttt{ 4039244593121563: 14}\\
  \texttt{ 4039632702153915: 14}\\
  \texttt{ 4040904548858012: 14}\\
  \texttt{ 4045601736258919: 14}\\
  \texttt{ 4047213783530167: 14}\\
  \texttt{ 4051507752115323: 14}\\
  \texttt{ 4051803592529647: 14}\\
  \texttt{ 4053053283802719: 14}\\
  \texttt{ 4053783700276216: 14}\\
  \texttt{ 4053809045956975: 14}\\
  \texttt{ 4056863101381448: 14}\\
  \texttt{ 4058223614548070: 15}\\
  \texttt{ 4059527580809864: 14}\\
  \texttt{ 4060738357329615: 14}

  \texttt{ 4062124786924240: 14}\\
  \texttt{ 4062858795180512: 14}\\
  \texttt{ 4064801792366815: 14}\\
  \texttt{ 4065239598037215: 14}\\
  \texttt{ 4065670135260367: 14}\\
  \texttt{ 4068391627124144: 14}\\
  \texttt{ 4070303365705063: 14}\\
  \texttt{ 4071078565708184: 14}\\
  \texttt{ 4074008423031872: 14}\\
  \texttt{ 4077376143450640: 14}\\
  \texttt{ 4077721351725848: 14}\\
  \texttt{ 4078295973276215: 14}\\
  \texttt{ 4079079240102800: 14}\\
  \texttt{ 4079534019773047: 14}\\
  \texttt{ 4080125919846220: 14}\\
  \texttt{ 4085395978939324: 14}\\
  \texttt{ 4088003437114720: 14}\\
  \texttt{ 4088196867693543: 14}\\
  \texttt{ 4089536314521928: 14}\\
  \texttt{ 4091153666928415: 14}

  \texttt{ 4093725904234323: 14}\\
  \texttt{ 4094308695479775: 14}\\
  \texttt{ 4095389764559347: 14}\\
  \texttt{ 4095647774315224: 14}\\
  \texttt{ 4096191089264048: 14}\\
  \texttt{ 4097056684813071: 14}\\
  \texttt{ 4102653962069864: 14}\\
  \texttt{ 4103071949930415: 14}\\
  \texttt{ 4108799499739719: 14}\\
  \texttt{ 4109183867293360: 14}\\
  \texttt{ 4110973137543067: 14}\\
  \texttt{ 4112181220461048: 14}\\
  \texttt{ 4112550314611172: 14}\\
  \texttt{ 4113729390733172: 14}\\
  \texttt{ 4113921746472343: 14}\\
  \texttt{ 4114837303671643: 14}\\
  \texttt{ 4118867947872624: 14}\\
  \texttt{ 4120030395123344: 14}\\
  \texttt{ 4127667133101668: 14}\\
  \texttt{ 4140607095481767: 14}

  \texttt{ 4141318772292919: 15}\\
  \texttt{ 4142196506877415: 14}\\
  \texttt{ 4143323873466547: 14}\\
  \texttt{ 4144133680280719: 14}\\
  \texttt{ 4145379818462216: 14}\\
  \texttt{ 4147872438384620: 14}\\
  \texttt{ 4148766359303371: 14}\\
  \texttt{ 4149874657059272: 14}\\
  \texttt{ 4151642724477339: 14}\\
  \texttt{ 4152381684467416: 14}\\
  \texttt{ 4152911488004139: 14}\\
  \texttt{ 4153319253548872: 14}\\
  \texttt{ 4156689614307847: 14}\\
  \texttt{ 4159074972623739: 14}\\
  \texttt{ 4160150067213475: 14}\\
  \texttt{ 4160174854925967: 14}\\
  \texttt{ 4160875005159812: 14}\\
  \texttt{ 4163067696115023: 14}\\
  \texttt{ 4165653662382715: 14}\\
  \texttt{ 4169560978014675: 14}

  \texttt{ 4169682492150643: 14}\\
  \texttt{ 4170336326116719: 14}\\
  \texttt{ 4175197987809444: 14}\\
  \texttt{ 4182314568661916: 14}\\
  \texttt{ 4183473469470019: 14}\\
  \texttt{ 4185691645841872: 14}\\
  \texttt{ 4188074217925012: 14}\\
  \texttt{ 4189154456918319: 14}\\
  \texttt{ 4189300712263016: 14}\\
  \texttt{ 4190771876257520: 14}\\
  \texttt{ 4191972415237167: 14}\\
  \texttt{ 4192113607752572: 14}\\
  \texttt{ 4192739768041315: 14}\\
  \texttt{ 4192897593807124: 14}\\
  \texttt{ 4193277113925872: 14}\\
  \texttt{ 4195225535899616: 14}\\
  \texttt{ 4196935091887215: 14}\\
  \texttt{ 4199286469536548: 14}\\
  \texttt{ 4200206279759647: 14}\\
  \texttt{ 4200283484620143: 14}

  \texttt{ 4203200912410924: 14}\\
  \texttt{ 4206049235471672: 14}\\
  \texttt{ 4211203686133472: 14}\\
  \texttt{ 4211418057132868: 14}\\
  \texttt{ 4215256385231564: 14}\\
  \texttt{ 4217146391008744: 14}\\
  \texttt{ 4218024408756891: 14}\\
  \texttt{ 4219043518298275: 14}\\
  \texttt{ 4220370838700216: 14}\\
  \texttt{ 4223010953651444: 14}\\
  \texttt{ 4226941476610544: 14}\\
  \texttt{ 4228196386446123: 14}\\
  \texttt{ 4228249602000975: 14}\\
  \texttt{ 4229166136721524: 14}\\
  \texttt{ 4229603860031944: 14}\\
  \texttt{ 4229660325805072: 14}\\
  \texttt{ 4230284672805316: 14}\\
  \texttt{ 4231557780690247: 14}\\
  \texttt{ 4232327887818043: 14}\\
  \texttt{ 4233328030097223: 14}

  \texttt{ 4233850623606219: 14}\\
  \texttt{ 4234781965797475: 14}\\
  \texttt{ 4236277574136615: 14}\\
  \texttt{ 4239120621307912: 14}\\
  \texttt{ 4240902035857364: 14}\\
  \texttt{ 4240980950542567: 14}\\
  \texttt{ 4241666373604171: 14}\\
  \texttt{ 4243985855991616: 15}\\
  \texttt{ 4251861731588768: 14}\\
  \texttt{ 4252465946106547: 14}\\
  \texttt{ 4253094960453267: 14}\\
  \texttt{ 4254592002735212: 14}\\
  \texttt{ 4258537525482416: 14}\\
  \texttt{ 4260299094935367: 14}\\
  \texttt{ 4262857654716339: 14}\\
  \texttt{ 4265175428324720: 14}\\
  \texttt{ 4266155279746672: 14}\\
  \texttt{ 4266253440193112: 15}\\
  \texttt{ 4266261354788944: 14}\\
  \texttt{ 4267900897523144: 14}

  \texttt{ 4270036586218468: 14}\\
  \texttt{ 4271816560454944: 14}\\
  \texttt{ 4272909493879916: 14}\\
  \texttt{ 4274424734380540: 14}\\
  \texttt{ 4278578312689575: 14}\\
  \texttt{ 4278616884019420: 14}\\
  \texttt{ 4278958404173224: 14}\\
  \texttt{ 4286531442633020: 14}\\
  \texttt{ 4291182815054443: 14}\\
  \texttt{ 4291204356077515: 14}\\
  \texttt{ 4291449693142443: 14}\\
  \texttt{ 4294356047725364: 14}\\
  \texttt{ 4295002193131764: 14}\\
  \texttt{ 4300924418206244: 14}\\
  \texttt{ 4302386171466415: 14}\\
  \texttt{ 4303811791208563: 14}\\
  \texttt{ 4306521221102672: 14}\\
  \texttt{ 4310101652528967: 14}\\
  \texttt{ 4314616208699415: 14}\\
  \texttt{ 4314782548438316: 14}

  \texttt{ 4315964916789123: 14}\\
  \texttt{ 4317351570252843: 14}\\
  \texttt{ 4318207233330067: 14}\\
  \texttt{ 4319427175491715: 14}\\
  \texttt{ 4323578311105171: 14}\\
  \texttt{ 4324017727637439: 14}\\
  \texttt{ 4332186255032224: 14}\\
  \texttt{ 4332469722403923: 14}\\
  \texttt{ 4335486754336819: 14}\\
  \texttt{ 4339105492392368: 14}\\
  \texttt{ 4342176630734515: 14}\\
  \texttt{ 4347302309758447: 14}\\
  \texttt{ 4349627050799975: 14}\\
  \texttt{ 4352296991350319: 14}\\
  \texttt{ 4352355084459463: 14}\\
  \texttt{ 4355791673702371: 14}\\
  \texttt{ 4358068208937119: 14}\\
  \texttt{ 4362191270781639: 14}\\
  \texttt{ 4363301404986344: 14}\\
  \texttt{ 4363707330770667: 14}

  \texttt{ 4366155295993767: 14}\\
  \texttt{ 4367114495704275: 14}\\
  \texttt{ 4370139613114323: 14}\\
  \texttt{ 4370729265627272: 14}\\
  \texttt{ 4370940635562520: 14}\\
  \texttt{ 4371533093380839: 14}\\
  \texttt{ 4371557279360812: 14}\\
  \texttt{ 4374613893395743: 15}\\
  \texttt{ 4376791649399371: 14}\\
  \texttt{ 4378259850082718: 15}\\
  \texttt{ 4382225392525315: 14}\\
  \texttt{ 4382613780816467: 14}\\
  \texttt{ 4384694040979219: 14}\\
  \texttt{ 4384907141185771: 14}\\
  \texttt{ 4386921632927216: 14}\\
  \texttt{ 4387693218189164: 14}\\
  \texttt{ 4388003998667224: 14}\\
  \texttt{ 4388350557229215: 14}\\
  \texttt{ 4391105530124571: 14}\\
  \texttt{ 4395620991526612: 14}

  \texttt{ 4396942387436344: 14}\\
  \texttt{ 4397258445569116: 14}\\
  \texttt{ 4400448106426243: 14}\\
  \texttt{ 4402098537360920: 14}\\
  \texttt{ 4403797127689768: 14}\\
  \texttt{ 4410240574792112: 15}\\
  \texttt{ 4410287989098164: 14}\\
  \texttt{ 4412740814516612: 14}\\
  \texttt{ 4413227518876768: 14}\\
  \texttt{ 4413890323427775: 14}\\
  \texttt{ 4415228672707268: 14}\\
  \texttt{ 4417684920497643: 14}\\
  \texttt{ 4419442561514120: 14}\\
  \texttt{ 4419669519874843: 14}\\
  \texttt{ 4420745511172468: 14}\\
  \texttt{ 4422123003381567: 14}\\
  \texttt{ 4422346678633364: 14}\\
  \texttt{ 4422655867209172: 14}\\
  \texttt{ 4424483184229923: 14}\\
  \texttt{ 4427187076388312: 14}

  \texttt{ 4427713646432619: 14}\\
  \texttt{ 4428185204145039: 14}\\
  \texttt{ 4431915098020663: 14}\\
  \texttt{ 4439095316242315: 14}\\
  \texttt{ 4439685430929040: 14}\\
  \texttt{ 4441013979770572: 14}\\
  \texttt{ 4442041807898344: 14}\\
  \texttt{ 4443538063290315: 14}\\
  \texttt{ 4444008108358616: 14}\\
  \texttt{ 4444042212950343: 14}\\
  \texttt{ 4446183328042515: 14}\\
  \texttt{ 4446542609704167: 14}\\
  \texttt{ 4452241159996324: 14}\\
  \texttt{ 4452360909543123: 14}\\
  \texttt{ 4453176381564472: 14}\\
  \texttt{ 4455510962025319: 14}\\
  \texttt{ 4456930876907812: 14}\\
  \texttt{ 4457072319677043: 14}\\
  \texttt{ 4457701611313244: 14}\\
  \texttt{ 4459623857845572: 14}

  \texttt{ 4462612560895868: 14}\\
  \texttt{ 4470188459417863: 14}\\
  \texttt{ 4475971209357267: 14}\\
  \texttt{ 4476999427048412: 14}\\
  \texttt{ 4477698265127048: 14}\\
  \texttt{ 4478441283422468: 14}\\
  \texttt{ 4479059921100243: 14}\\
  \texttt{ 4479264443710815: 14}\\
  \texttt{ 4479509212489412: 14}\\
  \texttt{ 4482164732139124: 14}\\
  \texttt{ 4482546582117220: 14}\\
  \texttt{ 4483559226202323: 14}\\
  \texttt{ 4489294846091020: 14}\\
  \texttt{ 4492112034232771: 14}\\
  \texttt{ 4492355647730740: 14}\\
  \texttt{ 4495716752691940: 14}\\
  \texttt{ 4496367632944012: 14}\\
  \texttt{ 4500659291889644: 14}\\
  \texttt{ 4503766986376964: 14}\\
  \texttt{ 4505562349077644: 14}

  \texttt{ 4506201959502472: 14}\\
  \texttt{ 4508969580839571: 14}\\
  \texttt{ 4510849153747444: 14}\\
  \texttt{ 4512659672445415: 14}\\
  \texttt{ 4514077965620975: 14}\\
  \texttt{ 4514102583251020: 14}\\
  \texttt{ 4515393915568015: 14}\\
  \texttt{ 4516344337070420: 14}\\
  \texttt{ 4516893068034175: 14}\\
  \texttt{ 4518487286420416: 14}\\
  \texttt{ 4520858276662215: 14}\\
  \texttt{ 4522999201634347: 14}\\
  \texttt{ 4523238197945115: 14}\\
  \texttt{ 4524409057386867: 14}\\
  \texttt{ 4524605131351768: 14}\\
  \texttt{ 4525732369241971: 14}\\
  \texttt{ 4526788386247820: 14}\\
  \texttt{ 4528635981902043: 14}\\
  \texttt{ 4528740291388467: 14}\\
  \texttt{ 4532358695741468: 14}

  \texttt{ 4534270962820468: 14}\\
  \texttt{ 4535627483571164: 14}\\
  \texttt{ 4536252086626772: 14}\\
  \texttt{ 4537466720779875: 14}\\
  \texttt{ 4538636394915544: 14}\\
  \texttt{ 4539277467931767: 14}\\
  \texttt{ 4541009743049739: 14}\\
  \texttt{ 4541583117454468: 14}\\
  \texttt{ 4541829279909639: 14}\\
  \texttt{ 4543669958220123: 14}\\
  \texttt{ 4545874488197815: 14}\\
  \texttt{ 4549075449077968: 14}\\
  \texttt{ 4550759910069344: 14}\\
  \texttt{ 4553159275109167: 14}\\
  \texttt{ 4553177196968143: 14}\\
  \texttt{ 4555831219504324: 14}\\
  \texttt{ 4556291379267942: 15}\\
  \texttt{ 4556865796385812: 14}\\
  \texttt{ 4556875760417572: 14}\\
  \texttt{ 4558133758997420: 14}

  \texttt{ 4559293516803272: 14}\\
  \texttt{ 4559435040722775: 14}\\
  \texttt{ 4559605578782648: 14}\\
  \texttt{ 4563420795786247: 14}\\
  \texttt{ 4563946198737044: 14}\\
  \texttt{ 4567950523467771: 14}\\
  \texttt{ 4568457296511267: 14}\\
  \texttt{ 4570749927642867: 14}\\
  \texttt{ 4571247293859715: 14}\\
  \texttt{ 4572280060278819: 14}\\
  \texttt{ 4573729091597048: 14}\\
  \texttt{ 4579495785992824: 14}\\
  \texttt{ 4580655407318671: 14}\\
  \texttt{ 4581368673256519: 14}\\
  \texttt{ 4586320612042540: 14}\\
  \texttt{ 4589032506626067: 14}\\
  \texttt{ 4590942205110171: 14}\\
  \texttt{ 4591572692013412: 14}\\
  \texttt{ 4593211093331368: 14}\\
  \texttt{ 4593721351665039: 14}

  \texttt{ 4598043124135768: 14}\\
  \texttt{ 4599496860301447: 14}\\
  \texttt{ 4600603637395171: 14}\\
  \texttt{ 4601649928283647: 14}\\
  \texttt{ 4604553683061163: 14}\\
  \texttt{ 4607160035749468: 14}\\
  \texttt{ 4608617994397844: 14}\\
  \texttt{ 4611846334099975: 14}\\
  \texttt{ 4612651407918520: 14}\\
  \texttt{ 4614004045357216: 14}\\
  \texttt{ 4614863545164320: 15}\\
  \texttt{ 4614975064757212: 14}\\
  \texttt{ 4615489978578339: 14}\\
  \texttt{ 4616035138667815: 14}\\
  \texttt{ 4616489562809047: 15}\\
  \texttt{ 4617627485579864: 14}\\
  \texttt{ 4620921169521832: 14}\\
  \texttt{ 4621173239683419: 14}\\
  \texttt{ 4622167924351071: 14}\\
  \texttt{ 4624852422920048: 14}

  \texttt{ 4626299014555464: 14}\\
  \texttt{ 4626456046392212: 14}\\
  \texttt{ 4629730297246215: 14}\\
  \texttt{ 4630574005536571: 14}\\
  \texttt{ 4632875938858615: 14}\\
  \texttt{ 4632981437142412: 14}\\
  \texttt{ 4634808453846619: 14}\\
  \texttt{ 4634835041199724: 14}\\
  \texttt{ 4636205894673763: 14}\\
  \texttt{ 4638323140004823: 14}\\
  \texttt{ 4639103189253548: 14}\\
  \texttt{ 4639577579905875: 14}\\
  \texttt{ 4644842491340440: 14}\\
  \texttt{ 4646740740801712: 14}\\
  \texttt{ 4654152915607420: 14}\\
  \texttt{ 4654364880809671: 14}\\
  \texttt{ 4657449569729372: 14}\\
  \texttt{ 4660927271464867: 14}\\
  \texttt{ 4661322105815872: 14}\\
  \texttt{ 4661569755213243: 14}

  \texttt{ 4663130928213464: 14}\\
  \texttt{ 4663656830372540: 14}\\
  \texttt{ 4663817123411603: 14}\\
  \texttt{ 4665480752711563: 14}\\
  \texttt{ 4665505839474463: 14}\\
  \texttt{ 4665904684949824: 14}\\
  \texttt{ 4666889469573943: 14}\\
  \texttt{ 4667818542457543: 14}\\
  \texttt{ 4667921709103771: 14}\\
  \texttt{ 4672222060157971: 14}\\
  \texttt{ 4675254500014540: 14}\\
  \texttt{ 4679497443684868: 14}\\
  \texttt{ 4679847500576348: 14}\\
  \texttt{ 4679922875116616: 14}\\
  \texttt{ 4680494566971712: 14}\\
  \texttt{ 4681481250682743: 14}\\
  \texttt{ 4683265283137875: 14}\\
  \texttt{ 4683317448473339: 14}\\
  \texttt{ 4685661454725115: 14}\\
  \texttt{ 4686831270127924: 14}

  \texttt{ 4687043595947440: 14}\\
  \texttt{ 4691114700688168: 14}\\
  \texttt{ 4691287768470171: 14}\\
  \texttt{ 4699364880802867: 14}\\
  \texttt{ 4699784901168475: 14}\\
  \texttt{ 4702327522825119: 14}\\
  \texttt{ 4703713699817312: 14}\\
  \texttt{ 4703794005821371: 14}\\
  \texttt{ 4704042098850919: 14}\\
  \texttt{ 4704526521556064: 14}\\
  \texttt{ 4708095695322272: 14}\\
  \texttt{ 4709306377423239: 14}\\
  \texttt{ 4709901520597471: 14}\\
  \texttt{ 4714237269224223: 14}\\
  \texttt{ 4715277540442540: 14}\\
  \texttt{ 4715922064598648: 14}\\
  \texttt{ 4717538912254220: 14}\\
  \texttt{ 4719734643005072: 14}\\
  \texttt{ 4722583112901475: 14}\\
  \texttt{ 4726290668827219: 14}

  \texttt{ 4726344852400219: 14}\\
  \texttt{ 4727162938671339: 14}\\
  \texttt{ 4727445121141712: 14}\\
  \texttt{ 4728149670583815: 14}\\
  \texttt{ 4730822860371244: 14}\\
  \texttt{ 4733281854330848: 14}\\
  \texttt{ 4735944917481975: 14}\\
  \texttt{ 4736593585127212: 14}\\
  \texttt{ 4738224853611412: 14}\\
  \texttt{ 4741288060976440: 15}\\
  \texttt{ 4742237220512443: 14}\\
  \texttt{ 4743605805173072: 14}\\
  \texttt{ 4745587160008663: 14}\\
  \texttt{ 4751261970957016: 14}\\
  \texttt{ 4751634926279212: 14}\\
  \texttt{ 4755333277156543: 14}\\
  \texttt{ 4755648422524867: 14}\\
  \texttt{ 4756655823263468: 14}\\
  \texttt{ 4758802842488343: 14}\\
  \texttt{ 4758893159486343: 14}

  \texttt{ 4759941175770175: 14}\\
  \texttt{ 4759972360988740: 14}\\
  \texttt{ 4763765569577371: 14}\\
  \texttt{ 4765423493578616: 14}\\
  \texttt{ 4766043874166440: 14}\\
  \texttt{ 4766381581832648: 14}\\
  \texttt{ 4766415952476647: 14}\\
  \texttt{ 4768912682435319: 14}\\
  \texttt{ 4772083231508047: 14}\\
  \texttt{ 4772779590112215: 14}\\
  \texttt{ 4775358747030244: 14}\\
  \texttt{ 4777602082849772: 14}\\
  \texttt{ 4779973343542615: 15}\\
  \texttt{ 4780469810957263: 14}\\
  \texttt{ 4785620614768239: 14}\\
  \texttt{ 4787036662123420: 14}\\
  \texttt{ 4788090005614972: 14}\\
  \texttt{ 4790550286133215: 14}\\
  \texttt{ 4791506339670603: 14}\\
  \texttt{ 4794067169983816: 14}

  \texttt{ 4794363588138667: 14}\\
  \texttt{ 4794632484793243: 14}\\
  \texttt{ 4795034885954943: 14}\\
  \texttt{ 4798968422055927: 14}\\
  \texttt{ 4802108511287967: 14}\\
  \texttt{ 4802320774473668: 14}\\
  \texttt{ 4803369840961744: 15}\\
  \texttt{ 4804759244092119: 14}\\
  \texttt{ 4807199279428575: 14}\\
  \texttt{ 4807974482262844: 14}\\
  \texttt{ 4809808741779367: 14}\\
  \texttt{ 4810278630931972: 14}\\
  \texttt{ 4811280888459616: 14}\\
  \texttt{ 4814338358832164: 14}\\
  \texttt{ 4818077250129843: 14}\\
  \texttt{ 4822391741224075: 14}\\
  \texttt{ 4823633886028171: 14}\\
  \texttt{ 4825672528401124: 14}\\
  \texttt{ 4826450542087447: 14}\\
  \texttt{ 4826560311048319: 14}

  \texttt{ 4829323239437515: 14}\\
  \texttt{ 4830220823824015: 14}\\
  \texttt{ 4832093867865847: 14}\\
  \texttt{ 4832405975163715: 14}\\
  \texttt{ 4833968595150643: 14}\\
  \texttt{ 4834555624131667: 14}\\
  \texttt{ 4834865166698175: 14}\\
  \texttt{ 4835092393376475: 14}\\
  \texttt{ 4835107134036772: 14}\\
  \texttt{ 4835762606753944: 14}\\
  \texttt{ 4839175141622647: 14}\\
  \texttt{ 4844226281847271: 14}\\
  \texttt{ 4844646052408375: 14}\\
  \texttt{ 4845509678805415: 14}\\
  \texttt{ 4845867666993771: 14}\\
  \texttt{ 4850028510479043: 14}\\
  \texttt{ 4858036665262748: 14}\\
  \texttt{ 4859388698779168: 14}\\
  \texttt{ 4861587225340267: 14}\\
  \texttt{ 4861800321427240: 14}

  \texttt{ 4862553993152864: 14}\\
  \texttt{ 4865193725028572: 14}\\
  \texttt{ 4867028082874167: 14}\\
  \texttt{ 4867431485757271: 14}\\
  \texttt{ 4869596497936371: 14}\\
  \texttt{ 4870872252081764: 14}\\
  \texttt{ 4875272330042212: 14}\\
  \texttt{ 4875803874944468: 14}\\
  \texttt{ 4878207842669944: 14}\\
  \texttt{ 4879314387285164: 14}\\
  \texttt{ 4879390684607368: 14}\\
  \texttt{ 4879752947537344: 14}\\
  \texttt{ 4880750936119467: 14}\\
  \texttt{ 4881889445552043: 14}\\
  \texttt{ 4884042698897515: 14}\\
  \texttt{ 4884612994076343: 14}\\
  \texttt{ 4886791918360172: 14}\\
  \texttt{ 4889033513593964: 14}\\
  \texttt{ 4889534536205072: 14}\\
  \texttt{ 4889649556188519: 14}

  \texttt{ 4894079801130043: 14}\\
  \texttt{ 4896307222235367: 14}\\
  \texttt{ 4898120583889916: 14}\\
  \texttt{ 4899895684259948: 14}\\
  \texttt{ 4901208578201824: 15}\\
  \texttt{ 4901700736134367: 14}\\
  \texttt{ 4904366051894264: 14}\\
  \texttt{ 4905658080190743: 14}\\
  \texttt{ 4905711179468467: 14}\\
  \texttt{ 4906830959069739: 14}\\
  \texttt{ 4907130490172348: 14}\\
  \texttt{ 4907866858213443: 14}\\
  \texttt{ 4912190886313312: 14}\\
  \texttt{ 4916142608605820: 14}\\
  \texttt{ 4916925416400216: 14}\\
  \texttt{ 4918241655536416: 14}\\
  \texttt{ 4918588944666020: 14}\\
  \texttt{ 4919168728239640: 14}\\
  \texttt{ 4921040783952368: 14}\\
  \texttt{ 4921683567326863: 14}

  \texttt{ 4921860698325620: 14}\\
  \texttt{ 4927173673528843: 14}\\
  \texttt{ 4928097597182871: 14}\\
  \texttt{ 4929038302986843: 14}\\
  \texttt{ 4930738812264616: 14}\\
  \texttt{ 4932206869587639: 14}\\
  \texttt{ 4932778564395844: 14}\\
  \texttt{ 4933252009193524: 14}\\
  \texttt{ 4933390830228764: 14}\\
  \texttt{ 4933880439286616: 14}\\
  \texttt{ 4940283017484016: 14}\\
  \texttt{ 4946649121567015: 14}\\
  \texttt{ 4948525855332620: 14}\\
  \texttt{ 4949327202453944: 14}\\
  \texttt{ 4949715192775612: 15}\\
  \texttt{ 4949716591624623: 14}\\
  \texttt{ 4951161303933212: 14}\\
  \texttt{ 4951353519074175: 14}\\
  \texttt{ 4952340909538467: 14}\\
  \texttt{ 4954690297705443: 14}

  \texttt{ 4956455513274019: 14}\\
  \texttt{ 4958311426238572: 14}\\
  \texttt{ 4958593076774068: 14}\\
  \texttt{ 4959187310228912: 14}\\
  \texttt{ 4960478451947668: 14}\\
  \texttt{ 4961338630914715: 14}\\
  \texttt{ 4963202101206764: 14}\\
  \texttt{ 4963675534727371: 14}\\
  \texttt{ 4963785075770236: 15}\\
  \texttt{ 4964089945901668: 14}\\
  \texttt{ 4964733953132912: 14}\\
  \texttt{ 4965231157714447: 14}\\
  \texttt{ 4965521475624171: 14}\\
  \texttt{ 4965779166461440: 14}\\
  \texttt{ 4965932603903667: 14}\\
  \texttt{ 4966342324410967: 14}\\
  \texttt{ 4968150420718816: 14}\\
  \texttt{ 4970205566749924: 14}\\
  \texttt{ 4970618047415564: 14}\\
  \texttt{ 4971036687820963: 14}

  \texttt{ 4971571310175667: 14}\\
  \texttt{ 4972523044557075: 14}\\
  \texttt{ 4972684791538815: 14}\\
  \texttt{ 4974046565725672: 14}\\
  \texttt{ 4977657283939712: 14}\\
  \texttt{ 4982935875443640: 14}\\
  \texttt{ 4983296831129715: 14}\\
  \texttt{ 4983357231040867: 14}\\
  \texttt{ 4984113759118712: 14}\\
  \texttt{ 4991834439406168: 14}\\
  \texttt{ 4992022967144012: 14}\\
  \texttt{ 4997318137620604: 14}\\
  \texttt{ 4997445200449748: 14}\\
  \texttt{ 4998638816275840: 14}\\
  \texttt{ 4998881305807540: 14}\\
  \texttt{ 4999276165496047: 14}
\end{flushright}
}
\normalsize
\end{multicols}

\clearpage

\section*{Appendix C}

  Square-free gaps and their length $\geq 15$, up to $125\,870\,000\,000\,000\,000$

\setlength{\columnseprule}{0.1pt}
\begin{multicols}{6}
{
\tiny
\begin{flushright}
 \texttt{    79180770078548: 15}\\
 \texttt{   145065154350544: 15}\\
 \texttt{   222400208995718: 15}\\
 \texttt{   314475532576612: 15}\\
 \texttt{   522927463507963: 15}\\
 \texttt{   578323272578775: 15}\\
 \texttt{   597567054281944: 15}\\
 \texttt{   600648730444048: 15}\\
 \texttt{   613760263741842: 15}\\
 \texttt{   653292828286964: 15}\\
 \texttt{   756143836911943: 15}\\
 \texttt{   767966021500362: 15}\\
 \texttt{   799368239412412: 15}\\
 \texttt{   832835838515415: 15}\\
 \texttt{   845891516833240: 15}\\
 \texttt{   879322694418366: 15}\\
 \texttt{   963244642044616: 15}\\
 \texttt{   972256492049738: 15}\\
 \texttt{   978405629137168: 15}\\
 \texttt{  1033783215814412: 15}

 \texttt{  1133510837444575: 15}\\
 \texttt{  1186126171926711: 15}\\
 \texttt{  1204986929396019: 15}\\
 \texttt{  1407999508293715: 15}\\
 \texttt{  1413326736000616: 15}\\
 \texttt{  1482442771986122: 15}\\
 \texttt{  1640741786210644: 15}\\
 \texttt{  1647423474094275: 15}\\
 \texttt{  1666895115919911: 15}\\
 \texttt{  1728562371923214: 15}\\
 \texttt{  1778730977125411: 15}\\
 \texttt{  1887179922247112: 15}\\
 \texttt{  1891143421462772: 15}\\
 \texttt{  2045972062427371: 15}\\
 \texttt{  2063694280707338: 15}\\
 \texttt{  2100613839529216: 15}\\
 \texttt{  2227348282719068: 15}\\
 \texttt{  2255046185932722: 15}\\
 \texttt{  2347972209573171: 15}\\
 \texttt{  2369907350225071: 15}

 \texttt{  2492317253656916: 15}\\
 \texttt{  2590615352654367: 15}\\
 \texttt{  2644991500984502: 15}\\
 \texttt{  2692776998404743: 15}\\
 \texttt{  2702229285212824: 15}\\
 \texttt{  2832442820527864: 15}\\
 \texttt{  2844943274028974: 15}\\
 \texttt{  2909961461164446: 15}\\
 \texttt{  2913111930559912: 15}\\
 \texttt{  2936807265402368: 15}\\
 \texttt{  2960866051123122: 15}\\
 \texttt{  3035065451166546: 15}\\
 \texttt{  3199658816229018: 15}\\
 \texttt{  3212805940336311: 15}\\
 \texttt{  3213097301495211: 15}\\
 \texttt{  3215226335143218: 16}\\
 \texttt{  3277263843303664: 15}\\
 \texttt{  3320563792818544: 15}\\
 \texttt{  3418009019023518: 15}\\
 \texttt{  3460363358703640: 15}

 \texttt{  3480702313680522: 15}\\
 \texttt{  3558748762728724: 15}\\
 \texttt{  3584157346962774: 15}\\
 \texttt{  3689017823664916: 15}\\
 \texttt{  3751014691621671: 15}\\
 \texttt{  3755504514040316: 15}\\
 \texttt{  3796267726981840: 15}\\
 \texttt{  3870556876108662: 15}\\
 \texttt{  3889861917379168: 15}\\
 \texttt{  3898971358941020: 15}\\
 \texttt{  4005072589122338: 15}\\
 \texttt{  4058223614548070: 15}\\
 \texttt{  4141318772292919: 15}\\
 \texttt{  4243985855991616: 15}\\
 \texttt{  4266253440193112: 15}\\
 \texttt{  4374613893395743: 15}\\
 \texttt{  4378259850082718: 15}\\
 \texttt{  4410240574792112: 15}\\
 \texttt{  4556291379267942: 15}\\
 \texttt{  4614863545164320: 15}

 \texttt{  4616489562809047: 15}\\
 \texttt{  4741288060976440: 15}\\
 \texttt{  4779973343542615: 15}\\
 \texttt{  4803369840961744: 15}\\
 \texttt{  4901208578201824: 15}\\
 \texttt{  4949715192775612: 15}\\
 \texttt{  4963785075770236: 15}\\
 \texttt{  5025199393613364: 15}\\
 \texttt{  5078936944878664: 15}\\
 \texttt{  5096634238198736: 15}\\
 \texttt{  5171925979705436: 15}\\
 \texttt{  5249193979343924: 15}\\
 \texttt{  5275768505520318: 15}\\
 \texttt{  5284432503811662: 15}\\
 \texttt{  5301656238720748: 15}\\
 \texttt{  5332171161640768: 15}\\
 \texttt{  5407818505280186: 15}\\
 \texttt{  5418214590485740: 15}\\
 \texttt{  5451045885520470: 15}\\
 \texttt{  5482046019239575: 15}

 \texttt{  5521807633558612: 15}\\
 \texttt{  5546802608765972: 15}\\
 \texttt{  5603222007191174: 15}\\
 \texttt{  5724277311046912: 15}\\
 \texttt{  5725811634818668: 15}\\
 \texttt{  5768852637964468: 15}\\
 \texttt{  5823894340613318: 15}\\
 \texttt{  5873838440657274: 15}\\
 \texttt{  5880710799040623: 15}\\
 \texttt{  6004127166886167: 15}\\
 \texttt{  6033735506469975: 15}\\
 \texttt{  6128885790854863: 15}\\
 \texttt{  6139576463974467: 15}\\
 \texttt{  6260553216799142: 15}\\
 \texttt{  6270284687577314: 15}\\
 \texttt{  6338967598608847: 15}\\
 \texttt{  6356196623384223: 15}\\
 \texttt{  6413775247990374: 15}\\
 \texttt{  6546811442608814: 15}\\
 \texttt{  6608802595963540: 15}

 \texttt{  6660957244082440: 15}\\
 \texttt{  6705784903596916: 15}\\
 \texttt{  6734750493340847: 15}\\
 \texttt{  6753950108132450: 15}\\
 \texttt{  6772948118469267: 15}\\
 \texttt{  6782670508120242: 15}\\
 \texttt{  6783748676789370: 15}\\
 \texttt{  6889494851636020: 15}\\
 \texttt{  7005055514315022: 15}\\
 \texttt{  7147241034408339: 15}\\
 \texttt{  7181673258212946: 15}\\
 \texttt{  7194347867584818: 15}\\
 \texttt{  7223287409723715: 15}\\
 \texttt{  7227571492291370: 15}\\
 \texttt{  7241759899849419: 15}\\
 \texttt{  7409549709923766: 15}\\
 \texttt{  7423235661096340: 15}\\
 \texttt{  7478410973211340: 15}\\
 \texttt{  7600925746759363: 15}\\
 \texttt{  7722363759245375: 15}

 \texttt{  7747653274059355: 15}\\
 \texttt{  7800412268663068: 15}\\
 \texttt{  7831532030941022: 15}\\
 \texttt{  7851912406342818: 15}\\
 \texttt{  7963937542390874: 15}\\
 \texttt{  8007955950410862: 15}\\
 \texttt{  8038699633110614: 15}\\
 \texttt{  8105157391248738: 15}\\
 \texttt{  8119594326617574: 15}\\
 \texttt{  8143987353234642: 15}\\
 \texttt{  8228701310684948: 15}\\
 \texttt{  8231521666115066: 15}\\
 \texttt{  8274076055001414: 15}\\
 \texttt{  8385421261742468: 15}\\
 \texttt{  8408867829823519: 15}\\
 \texttt{  8460703525630623: 15}\\
 \texttt{  8552920743802922: 15}\\
 \texttt{  8598688930184514: 15}\\
 \texttt{  8642039316381340: 15}\\
 \texttt{  8646363445574822: 15}

 \texttt{  8692754988610374: 15}\\
 \texttt{  8692842376367366: 15}\\
 \texttt{  8754639324102594: 15}\\
 \texttt{  8815854482611136: 15}\\
 \texttt{  8828839658181319: 15}\\
 \texttt{  8851853171026670: 15}\\
 \texttt{  8944795273928418: 15}\\
 \texttt{  8944910252078948: 15}\\
 \texttt{  8968645235212712: 15}\\
 \texttt{  9059546704494112: 15}\\
 \texttt{  9098247937109766: 15}\\
 \texttt{  9122773061481890: 15}\\
 \texttt{  9148751903593962: 15}\\
 \texttt{  9188323833374144: 15}\\
 \texttt{  9224093706369174: 15}\\
 \texttt{  9240221180323948: 15}\\
 \texttt{  9255896220394012: 15}\\
 \texttt{  9269193715143844: 15}\\
 \texttt{  9370751847296344: 15}\\
 \texttt{  9418551826916970: 15}

 \texttt{  9425465131282970: 15}\\
 \texttt{  9449519987818274: 15}\\
 \texttt{  9459503607180938: 15}\\
 \texttt{  9522923089143650: 15}\\
 \texttt{  9548621543639815: 15}\\
 \texttt{  9553372617854214: 15}\\
 \texttt{  9561170688115912: 15}\\
 \texttt{  9569635309334744: 15}\\
 \texttt{  9693938017810266: 15}\\
 \texttt{  9697883496328840: 15}\\
 \texttt{  9703874726945367: 15}\\
 \texttt{  9720767045578662: 15}\\
 \texttt{  9746404581938516: 15}\\
 \texttt{  9771037670100964: 15}\\
 \texttt{  9784931845628167: 15}\\
 \texttt{  9786138763169740: 15}\\
 \texttt{ 10015478227399011: 15}\\
 \texttt{ 10054431373893268: 15}\\
 \texttt{ 10080856956208148: 15}\\
 \texttt{ 10157932869382516: 15}

 \texttt{ 10217133791578636: 15}\\
 \texttt{ 10224074215382415: 15}\\
 \texttt{ 10249660978060840: 15}\\
 \texttt{ 10264823636423046: 15}\\
 \texttt{ 10269365081502942: 15}\\
 \texttt{ 10436122048826967: 15}\\
 \texttt{ 10515472197571663: 15}\\
 \texttt{ 10599509748425370: 15}\\
 \texttt{ 10647033819921547: 15}\\
 \texttt{ 10712940119544724: 15}\\
 \texttt{ 10723720959559275: 15}\\
 \texttt{ 10910041665689171: 15}\\
 \texttt{ 11072698600531122: 15}\\
 \texttt{ 11106817540102267: 15}\\
 \texttt{ 11115337587206047: 15}\\
 \texttt{ 11184601388486343: 15}\\
 \texttt{ 11218039339205164: 15}\\
 \texttt{ 11393381856846248: 15}\\
 \texttt{ 11418750794871520: 15}\\
 \texttt{ 11421982381058943: 15}

 \texttt{ 11463619995663423: 15}\\
 \texttt{ 11511916894618144: 15}\\
 \texttt{ 11542074883336770: 15}\\
 \texttt{ 11551478389888842: 15}\\
 \texttt{ 11571880481903742: 15}\\
 \texttt{ 11607074026141839: 15}\\
 \texttt{ 11652210251638720: 15}\\
 \texttt{ 11655019290392104: 15}\\
 \texttt{ 11908493804078370: 15}\\
 \texttt{ 11917664983427874: 15}\\
 \texttt{ 11952459036194116: 15}\\
 \texttt{ 11980005638274150: 15}\\
 \texttt{ 12020043595275412: 15}\\
 \texttt{ 12070061517763638: 15}\\
 \texttt{ 12119356531568140: 15}\\
 \texttt{ 12134293522106344: 15}\\
 \texttt{ 12214828761579772: 15}\\
 \texttt{ 12309902683501614: 15}\\
 \texttt{ 12321016079861564: 15}\\
 \texttt{ 12395701659599511: 15}

 \texttt{ 12415247540294516: 15}\\
 \texttt{ 12474105683477548: 15}\\
 \texttt{ 12546928319301843: 15}\\
 \texttt{ 12564328703647623: 15}\\
 \texttt{ 12621411821832664: 15}\\
 \texttt{ 12647624412536036: 15}\\
 \texttt{ 12665795853207868: 15}\\
 \texttt{ 12733569847150816: 15}\\
 \texttt{ 12745921565686768: 15}\\
 \texttt{ 12897479773507372: 15}\\
 \texttt{ 12900978231034972: 15}\\
 \texttt{ 12902377369413714: 15}\\
 \texttt{ 12967616693496942: 15}\\
 \texttt{ 12992275182527043: 15}\\
 \texttt{ 13092643519157419: 15}\\
 \texttt{ 13096779538545872: 15}\\
 \texttt{ 13119286020601650: 15}\\
 \texttt{ 13155529105064646: 15}\\
 \texttt{ 13156029324308642: 15}\\
 \texttt{ 13179454288365438: 15}

 \texttt{ 13230423200166642: 15}\\
 \texttt{ 13336235581408718: 15}\\
 \texttt{ 13367207076954307: 15}\\
 \texttt{ 13438746536847639: 15}\\
 \texttt{ 13485386513086072: 15}\\
 \texttt{ 13525625146720543: 15}\\
 \texttt{ 13607061807470522: 15}\\
 \texttt{ 13671384224051322: 15}\\
 \texttt{ 13718679611368074: 15}\\
 \texttt{ 13735088865755464: 15}\\
 \texttt{ 13835243834496763: 15}\\
 \texttt{ 13868041458469219: 15}\\
 \texttt{ 13973445220973018: 15}\\
 \texttt{ 14016066772804719: 15}\\
 \texttt{ 14079678087240187: 15}\\
 \texttt{ 14093573407227619: 15}\\
 \texttt{ 14098062086586366: 15}\\
 \texttt{ 14129252536617519: 15}\\
 \texttt{ 14191841695029163: 15}\\
 \texttt{ 14195272337432439: 15}

 \texttt{ 14309193668023240: 15}\\
 \texttt{ 14414464482748064: 15}\\
 \texttt{ 14507425899460962: 15}\\
 \texttt{ 14592292717722916: 15}\\
 \texttt{ 14676018384621643: 15}\\
 \texttt{ 14703612931303538: 15}\\
 \texttt{ 14720072414163170: 15}\\
 \texttt{ 14743351289011543: 15}\\
 \texttt{ 14769759863509816: 15}\\
 \texttt{ 14820127450780522: 15}\\
 \texttt{ 14879994360551871: 15}\\
 \texttt{ 14964746628637519: 15}\\
 \texttt{ 15045258413838211: 15}\\
 \texttt{ 15125440271645319: 15}\\
 \texttt{ 15213474122037243: 15}\\
 \texttt{ 15255761957438562: 15}\\
 \texttt{ 15293299937571724: 15}\\
 \texttt{ 15306198501236038: 15}\\
 \texttt{ 15354943464598575: 15}\\
 \texttt{ 15436928812383122: 15}

 \texttt{ 15443272697890446: 15}\\
 \texttt{ 15446800488184611: 15}\\
 \texttt{ 15447925378094670: 15}\\
 \texttt{ 15467048654444511: 15}\\
 \texttt{ 15646127795635471: 15}\\
 \texttt{ 15812203411096936: 15}\\
 \texttt{ 15872652634820047: 15}\\
 \texttt{ 15877488632684775: 15}\\
 \texttt{ 15883710354941224: 15}\\
 \texttt{ 15892635212087116: 15}\\
 \texttt{ 15909914997051674: 15}\\
 \texttt{ 15925092119500148: 15}\\
 \texttt{ 15949165417649574: 15}\\
 \texttt{ 15982488080689515: 15}\\
 \texttt{ 16011099811371066: 15}\\
 \texttt{ 16018839700850420: 15}\\
 \texttt{ 16032998156389792: 15}\\
 \texttt{ 16062176516460018: 15}\\
 \texttt{ 16067590106426871: 15}\\
 \texttt{ 16159218557178967: 15}

 \texttt{ 16236484838391544: 15}\\
 \texttt{ 16326870389341171: 15}\\
 \texttt{ 16338795878261967: 15}\\
 \texttt{ 16442609815208642: 15}\\
 \texttt{ 16497034382476814: 15}\\
 \texttt{ 16534038720652962: 15}\\
 \texttt{ 16763600249356444: 15}\\
 \texttt{ 16792860919488340: 15}\\
 \texttt{ 16871480124229311: 15}\\
 \texttt{ 16931409436935718: 15}\\
 \texttt{ 17113295679116766: 15}\\
 \texttt{ 17141979195426211: 15}\\
 \texttt{ 17151647578703720: 15}\\
 \texttt{ 17271804352011272: 15}\\
 \texttt{ 17280995138681274: 15}\\
 \texttt{ 17312603811478924: 15}\\
 \texttt{ 17333611288866763: 15}\\
 \texttt{ 17390116480129467: 15}\\
 \texttt{ 17460758160664038: 15}\\
 \texttt{ 17486223162781472: 15}

 \texttt{ 17486509233212467: 15}\\
 \texttt{ 17538956205656450: 15}\\
 \texttt{ 17544570884760615: 15}\\
 \texttt{ 17615826320835218: 15}\\
 \texttt{ 17765933370110871: 15}\\
 \texttt{ 17814838248562614: 15}\\
 \texttt{ 17873483158158342: 15}\\
 \texttt{ 17915805816941970: 15}\\
 \texttt{ 17952067794942018: 15}\\
 \texttt{ 17966344154636202: 15}\\
 \texttt{ 17997386326099074: 15}\\
 \texttt{ 18034366041788774: 15}\\
 \texttt{ 18089805640839822: 15}\\
 \texttt{ 18096806116986638: 15}\\
 \texttt{ 18116404429562216: 15}\\
 \texttt{ 18145304128677858: 15}\\
 \texttt{ 18154189822092416: 15}\\
 \texttt{ 18162933018392275: 15}\\
 \texttt{ 18239484944936340: 15}\\
 \texttt{ 18253355948504223: 15}

 \texttt{ 18260068551623872: 15}\\
 \texttt{ 18277373261269170: 15}\\
 \texttt{ 18288478350163736: 15}\\
 \texttt{ 18311971095518344: 15}\\
 \texttt{ 18332817515814571: 15}\\
 \texttt{ 18336189498865244: 15}\\
 \texttt{ 18375105176997162: 15}\\
 \texttt{ 18427288407259112: 15}\\
 \texttt{ 18429307057531672: 15}\\
 \texttt{ 18522608320698315: 15}\\
 \texttt{ 18555711153616674: 15}\\
 \texttt{ 18562210020271962: 15}\\
 \texttt{ 18575152279483864: 15}\\
 \texttt{ 18669891894522415: 15}\\
 \texttt{ 18724389862220262: 15}\\
 \texttt{ 18780616325786224: 15}\\
 \texttt{ 18787018968456775: 15}\\
 \texttt{ 18942958591820140: 15}\\
 \texttt{ 19207259606015368: 15}\\
 \texttt{ 19250332158522938: 15}

 \texttt{ 19271601923133543: 15}\\
 \texttt{ 19469517582240770: 15}\\
 \texttt{ 19536369378031062: 15}\\
 \texttt{ 19632091839856612: 15}\\
 \texttt{ 19636010709862072: 15}\\
 \texttt{ 19641445488546668: 15}\\
 \texttt{ 19664630616010616: 15}\\
 \texttt{ 19744635747012419: 15}\\
 \texttt{ 19775543284106719: 15}\\
 \texttt{ 19801588969954772: 15}\\
 \texttt{ 20071640031499736: 15}\\
 \texttt{ 20174431009347871: 15}\\
 \texttt{ 20177436695558019: 15}\\
 \texttt{ 20217899713449474: 15}\\
 \texttt{ 20396004281418572: 15}\\
 \texttt{ 20607307423661211: 15}\\
 \texttt{ 20727456415950074: 15}\\
 \texttt{ 20797493969707722: 15}\\
 \texttt{ 21140678272325368: 15}\\
 \texttt{ 21154813252229671: 15}

 \texttt{ 21290453932009218: 15}\\
 \texttt{ 21299530772891848: 15}\\
 \texttt{ 21315776264314818: 15}\\
 \texttt{ 21317332225643440: 15}\\
 \texttt{ 21439404970202142: 15}\\
 \texttt{ 21558677963307163: 15}\\
 \texttt{ 21560243519765418: 15}\\
 \texttt{ 21615804955958238: 15}\\
 \texttt{ 21628873538625614: 15}\\
 \texttt{ 21730943669770816: 15}\\
 \texttt{ 21832170429023143: 15}\\
 \texttt{ 21854478153869512: 15}\\
 \texttt{ 21903245395083116: 15}\\
 \texttt{ 21920105944311772: 15}\\
 \texttt{ 21942177566659816: 15}\\
 \texttt{ 21988275643067370: 15}\\
 \texttt{ 22033472803818918: 15}\\
 \texttt{ 22102900229210742: 15}\\
 \texttt{ 22209155578235140: 15}\\
 \texttt{ 22364692421673964: 15}

 \texttt{ 22380398534828740: 15}\\
 \texttt{ 22441630557281948: 15}\\
 \texttt{ 22457604407752420: 15}\\
 \texttt{ 22607928160985142: 15}\\
 \texttt{ 22654249221424374: 15}\\
 \texttt{ 22660912169860768: 15}\\
 \texttt{ 22712112344884518: 15}\\
 \texttt{ 22788702656697938: 15}\\
 \texttt{ 22790709882713642: 15}\\
 \texttt{ 22930984969340370: 15}\\
 \texttt{ 22954310080682418: 15}\\
 \texttt{ 22990095744478214: 15}\\
 \texttt{ 22995838937064247: 15}\\
 \texttt{ 23002410619741650: 15}\\
 \texttt{ 23060713753485914: 15}\\
 \texttt{ 23100086632764464: 15}\\
 \texttt{ 23115250421553939: 15}\\
 \texttt{ 23137057786622922: 15}\\
 \texttt{ 23173760867496018: 15}\\
 \texttt{ 23253659258767322: 15}

 \texttt{ 23347344489036243: 15}\\
 \texttt{ 23355295706140540: 15}\\
 \texttt{ 23358166368715912: 15}\\
 \texttt{ 23366878412928342: 15}\\
 \texttt{ 23550306148552838: 15}\\
 \texttt{ 23742453640900972: 17}\\
 \texttt{ 23812719364497075: 15}\\
 \texttt{ 23907609837614812: 15}\\
 \texttt{ 24019835094672342: 15}\\
 \texttt{ 24045900131677911: 15}\\
 \texttt{ 24051846024013471: 15}\\
 \texttt{ 24061292077720112: 15}\\
 \texttt{ 24143987317861050: 15}\\
 \texttt{ 24147626084787543: 15}\\
 \texttt{ 24151472002697522: 15}\\
 \texttt{ 24210738319978818: 15}\\
 \texttt{ 24211042123955264: 15}\\
 \texttt{ 24328487973992416: 15}\\
 \texttt{ 24381700900416872: 15}\\
 \texttt{ 24410934322913166: 15}

 \texttt{ 24423992997599140: 15}\\
 \texttt{ 24458086947835664: 15}\\
 \texttt{ 24485703864968632: 15}\\
 \texttt{ 24612362868025422: 15}\\
 \texttt{ 24622794062062072: 15}\\
 \texttt{ 24800784835835823: 15}\\
 \texttt{ 24837120189886072: 15}\\
 \texttt{ 24863115124207419: 15}\\
 \texttt{ 25032933903214670: 15}\\
 \texttt{ 25075098149546814: 15}\\
 \texttt{ 25081128698243336: 15}\\
 \texttt{ 25130103590698323: 15}\\
 \texttt{ 25153891525058944: 15}\\
 \texttt{ 25309333599521019: 15}\\
 \texttt{ 25309741937349064: 15}\\
 \texttt{ 25337693872297062: 15}\\
 \texttt{ 25369224817395040: 15}\\
 \texttt{ 25521917210723968: 15}\\
 \texttt{ 25638221576271064: 15}\\
 \texttt{ 25742413154880522: 15}

 \texttt{ 25810271112747138: 15}\\
 \texttt{ 25837546519722112: 15}\\
 \texttt{ 25906363158313372: 15}\\
 \texttt{ 25911663727121140: 15}\\
 \texttt{ 25949472393215644: 15}\\
 \texttt{ 25959922482580746: 15}\\
 \texttt{ 26002501473378470: 15}\\
 \texttt{ 26032476279109372: 15}\\
 \texttt{ 26088186234136075: 15}\\
 \texttt{ 26157296040785943: 15}\\
 \texttt{ 26185269034847416: 15}\\
 \texttt{ 26212867563525566: 15}\\
 \texttt{ 26253151694945150: 15}\\
 \texttt{ 26336501110587068: 15}\\
 \texttt{ 26388812512737150: 15}\\
 \texttt{ 26397000977534467: 15}\\
 \texttt{ 26416413080290948: 15}\\
 \texttt{ 26530125053387524: 15}\\
 \texttt{ 26596128217708923: 15}\\
 \texttt{ 26629952422662270: 15}

 \texttt{ 26699026280141670: 15}\\
 \texttt{ 26768781765764870: 15}\\
 \texttt{ 26779729800215066: 15}\\
 \texttt{ 26797570858482423: 15}\\
 \texttt{ 26808922329332766: 15}\\
 \texttt{ 26859187640357440: 15}\\
 \texttt{ 26888335201779450: 15}\\
 \texttt{ 26891820002571639: 15}\\
 \texttt{ 26945845035445674: 15}\\
 \texttt{ 27012735129114266: 15}\\
 \texttt{ 27017962989711462: 15}\\
 \texttt{ 27018184352218036: 15}\\
 \texttt{ 27032078415230114: 15}\\
 \texttt{ 27059335915841971: 15}\\
 \texttt{ 27126775437359224: 15}\\
 \texttt{ 27151632589742863: 15}\\
 \texttt{ 27309652121685570: 15}\\
 \texttt{ 27320732279851219: 15}\\
 \texttt{ 27356472820211512: 15}\\
 \texttt{ 27426371896236967: 15}

 \texttt{ 27454617329867068: 15}\\
 \texttt{ 27510974056481668: 15}\\
 \texttt{ 27549771006451875: 15}\\
 \texttt{ 27575133626412914: 15}\\
 \texttt{ 27600657593899436: 15}\\
 \texttt{ 27643996724021620: 15}\\
 \texttt{ 27673460878960264: 15}\\
 \texttt{ 27702210444412216: 15}\\
 \texttt{ 27707514102357462: 15}\\
 \texttt{ 27818956600257966: 15}\\
 \texttt{ 27857546333578242: 15}\\
 \texttt{ 27861449808318344: 15}\\
 \texttt{ 27924125942300450: 15}\\
 \texttt{ 28122857316285774: 15}\\
 \texttt{ 28254577733866412: 15}\\
 \texttt{ 28401126711070542: 15}\\
 \texttt{ 28421445515029566: 15}\\
 \texttt{ 28440244022892170: 15}\\
 \texttt{ 28447694163960918: 15}\\
 \texttt{ 28491631156615554: 15}

 \texttt{ 28509842601921650: 15}\\
 \texttt{ 28541488369763671: 15}\\
 \texttt{ 28566509244601444: 15}\\
 \texttt{ 28607164656350263: 15}\\
 \texttt{ 28644499886254712: 15}\\
 \texttt{ 28696958943616635: 16}\\
 \texttt{ 28697368516657243: 15}\\
 \texttt{ 28710396216477339: 15}\\
 \texttt{ 28722752581542772: 15}\\
 \texttt{ 28727968196758471: 15}\\
 \texttt{ 28775565362901068: 15}\\
 \texttt{ 28807071604957915: 15}\\
 \texttt{ 28810500521414044: 15}\\
 \texttt{ 28821646951850815: 15}\\
 \texttt{ 28925496074317444: 15}\\
 \texttt{ 28992793553607750: 15}\\
 \texttt{ 29023878409728318: 15}\\
 \texttt{ 29085609822273342: 15}\\
 \texttt{ 29098807367433812: 15}\\
 \texttt{ 29172249093745016: 15}

 \texttt{ 29179776157889740: 15}\\
 \texttt{ 29213456404340848: 15}\\
 \texttt{ 29220104210271811: 15}\\
 \texttt{ 29239405515409471: 15}\\
 \texttt{ 29240843822524448: 15}\\
 \texttt{ 29265426373447922: 15}\\
 \texttt{ 29314449222741512: 15}\\
 \texttt{ 29353228389308318: 15}\\
 \texttt{ 29402764296638775: 15}\\
 \texttt{ 29419390541454236: 15}\\
 \texttt{ 29503076040463766: 15}\\
 \texttt{ 29518142275240062: 15}\\
 \texttt{ 29548932602145138: 15}\\
 \texttt{ 29600273366809214: 15}\\
 \texttt{ 29601797902527716: 15}\\
 \texttt{ 29631931072747364: 15}\\
 \texttt{ 29642213575544262: 15}\\
 \texttt{ 29860161899379244: 15}\\
 \texttt{ 29873364911595772: 15}\\
 \texttt{ 29888562122792836: 15}

 \texttt{ 30011419590498736: 15}\\
 \texttt{ 30020039635915170: 15}\\
 \texttt{ 30070756955947766: 15}\\
 \texttt{ 30142513038680136: 15}\\
 \texttt{ 30241302600598214: 15}\\
 \texttt{ 30245541941284936: 15}\\
 \texttt{ 30250552293772016: 15}\\
 \texttt{ 30269528205805374: 15}\\
 \texttt{ 30288068888283512: 15}\\
 \texttt{ 30289919656561719: 15}\\
 \texttt{ 30329937653070471: 15}\\
 \texttt{ 30331853653893414: 15}\\
 \texttt{ 30388269903725419: 15}\\
 \texttt{ 30581913246599070: 15}\\
 \texttt{ 30653457721909362: 15}\\
 \texttt{ 30670566119307842: 15}\\
 \texttt{ 30820208253974719: 15}\\
 \texttt{ 30854241600116272: 15}\\
 \texttt{ 31037322636819150: 15}\\
 \texttt{ 31042625278344316: 15}

 \texttt{ 31047831901884114: 15}\\
 \texttt{ 31071527884060719: 15}\\
 \texttt{ 31083992760597819: 15}\\
 \texttt{ 31114811774579262: 15}\\
 \texttt{ 31123469975715836: 15}\\
 \texttt{ 31329752635278923: 15}\\
 \texttt{ 31333594366534336: 15}\\
 \texttt{ 31338044837959874: 15}\\
 \texttt{ 31401976920688950: 16}\\
 \texttt{ 31407188925679072: 15}\\
 \texttt{ 31514437999941272: 15}\\
 \texttt{ 31681764314356112: 15}\\
 \texttt{ 31697640625204850: 15}\\
 \texttt{ 31719814758871719: 15}\\
 \texttt{ 31749362657943171: 15}\\
 \texttt{ 31780094832667948: 15}\\
 \texttt{ 31790557927371474: 15}\\
 \texttt{ 31848298269038140: 15}\\
 \texttt{ 31902762590034570: 15}\\
 \texttt{ 31906568448917474: 15}

 \texttt{ 31938681992551215: 15}\\
 \texttt{ 32022834208745920: 15}\\
 \texttt{ 32051215736321739: 15}\\
 \texttt{ 32083581481004750: 15}\\
 \texttt{ 32093793808215270: 15}\\
 \texttt{ 32169850858644324: 15}\\
 \texttt{ 32256826315427970: 15}\\
 \texttt{ 32261589088376236: 15}\\
 \texttt{ 32324069045074014: 15}\\
 \texttt{ 32366866355032442: 15}\\
 \texttt{ 32371368943515711: 15}\\
 \texttt{ 32425029799201519: 15}\\
 \texttt{ 32442350234014923: 15}\\
 \texttt{ 32535618569275144: 15}\\
 \texttt{ 32537683646125820: 15}\\
 \texttt{ 32672434143501018: 15}\\
 \texttt{ 32760302921667619: 15}\\
 \texttt{ 32766544519827775: 15}\\
 \texttt{ 32821489279322138: 15}\\
 \texttt{ 32828412976160272: 15}

 \texttt{ 32854045723488367: 15}\\
 \texttt{ 32864574615973124: 15}\\
 \texttt{ 32957869900202836: 15}\\
 \texttt{ 33012332375884242: 15}\\
 \texttt{ 33128245851463014: 15}\\
 \texttt{ 33134211337961272: 15}\\
 \texttt{ 33138331698855042: 15}\\
 \texttt{ 33347813013777470: 15}\\
 \texttt{ 33384339110904412: 15}\\
 \texttt{ 33415816275260071: 15}\\
 \texttt{ 33590664281468644: 15}\\
 \texttt{ 33728989962295743: 15}\\
 \texttt{ 33792965701629162: 15}\\
 \texttt{ 33827524338460616: 15}\\
 \texttt{ 33829129778098875: 15}\\
 \texttt{ 33848817396108570: 15}\\
 \texttt{ 33933571295391114: 15}\\
 \texttt{ 33978439614311620: 15}\\
 \texttt{ 33992045852905375: 15}\\
 \texttt{ 34039561634998540: 15}

 \texttt{ 34298534579894799: 15}\\
 \texttt{ 34345483222146343: 15}\\
 \texttt{ 34537123786904212: 15}\\
 \texttt{ 34783674272812854: 15}\\
 \texttt{ 34833029806505312: 15}\\
 \texttt{ 34845036422155024: 15}\\
 \texttt{ 34846508676278346: 15}\\
 \texttt{ 34865208041499843: 15}\\
 \texttt{ 34889888998506174: 15}\\
 \texttt{ 34927299667114540: 15}\\
 \texttt{ 35001738610786448: 15}\\
 \texttt{ 35083937718284368: 15}\\
 \texttt{ 35101951090828722: 15}\\
 \texttt{ 35137853211332066: 15}\\
 \texttt{ 35140666007242050: 15}\\
 \texttt{ 35180073578101119: 15}\\
 \texttt{ 35270461029080347: 15}\\
 \texttt{ 35395419325403824: 15}\\
 \texttt{ 35440190399837715: 15}\\
 \texttt{ 35459526277026520: 15}

 \texttt{ 35486832857974316: 15}\\
 \texttt{ 35530064725976942: 15}\\
 \texttt{ 35548348331637518: 15}\\
 \texttt{ 35558686735015471: 15}\\
 \texttt{ 35596667337372267: 15}\\
 \texttt{ 35769379870018448: 15}\\
 \texttt{ 35833348501486744: 15}\\
 \texttt{ 35836757301666968: 15}\\
 \texttt{ 35850737764219770: 15}\\
 \texttt{ 35851111659622050: 15}\\
 \texttt{ 35857018303942122: 15}\\
 \texttt{ 35903366635583318: 15}\\
 \texttt{ 35960019264002572: 15}\\
 \texttt{ 36070897172081775: 15}\\
 \texttt{ 36109562408007272: 15}\\
 \texttt{ 36140538847155847: 15}\\
 \texttt{ 36236164881859747: 15}\\
 \texttt{ 36300127501436346: 15}\\
 \texttt{ 36422496339766470: 15}\\
 \texttt{ 36430733674573315: 15}

 \texttt{ 36516690595994812: 15}\\
 \texttt{ 36621341284722916: 15}\\
 \texttt{ 36650024229972642: 15}\\
 \texttt{ 36658238016995150: 15}\\
 \texttt{ 36679434688064538: 15}\\
 \texttt{ 36727088592737738: 15}\\
 \texttt{ 36743380242388747: 15}\\
 \texttt{ 36749464943648666: 15}\\
 \texttt{ 36793262656611868: 15}\\
 \texttt{ 36861653772304568: 15}\\
 \texttt{ 36896079687645572: 15}\\
 \texttt{ 36932094756515966: 15}\\
 \texttt{ 36944744083041039: 15}\\
 \texttt{ 36985881099122836: 17}\\
 \texttt{ 37016664275747419: 15}\\
 \texttt{ 37044993394045238: 15}\\
 \texttt{ 37082382811628944: 15}\\
 \texttt{ 37088392787051346: 15}\\
 \texttt{ 37100109324362011: 15}\\
 \texttt{ 37168850540803363: 15}

 \texttt{ 37195563651024642: 15}\\
 \texttt{ 37228865130724650: 15}\\
 \texttt{ 37261516054346968: 15}\\
 \texttt{ 37281051710896864: 15}\\
 \texttt{ 37297941262432323: 15}\\
 \texttt{ 37512805360455266: 15}\\
 \texttt{ 37520410817995047: 15}\\
 \texttt{ 37574766213275070: 15}\\
 \texttt{ 37579569995324070: 15}\\
 \texttt{ 37584136123280838: 15}\\
 \texttt{ 37668051748062616: 15}\\
 \texttt{ 37722602323707412: 15}\\
 \texttt{ 37742650751803698: 15}\\
 \texttt{ 37763081783080818: 15}\\
 \texttt{ 37931789143757522: 15}\\
 \texttt{ 38026989982107338: 15}\\
 \texttt{ 38101581390121566: 15}\\
 \texttt{ 38168787379560711: 15}\\
 \texttt{ 38293058303461912: 15}\\
 \texttt{ 38293077065161864: 15}

 \texttt{ 38316908336947746: 15}\\
 \texttt{ 38318171436267344: 15}\\
 \texttt{ 38663748860565724: 15}\\
 \texttt{ 38733865521271323: 15}\\
 \texttt{ 38987678717571248: 15}\\
 \texttt{ 38991642075913867: 15}\\
 \texttt{ 39033534094150622: 15}\\
 \texttt{ 39285301374402774: 15}\\
 \texttt{ 39356184408385470: 15}\\
 \texttt{ 39399808392972871: 15}\\
 \texttt{ 39457704763084768: 15}\\
 \texttt{ 39466961312814014: 15}\\
 \texttt{ 39528959396578711: 15}\\
 \texttt{ 39618867271380124: 15}\\
 \texttt{ 39858487272021112: 15}\\
 \texttt{ 39866926779486774: 15}\\
 \texttt{ 39913711088163594: 15}\\
 \texttt{ 39922043541990640: 15}\\
 \texttt{ 39942416165743744: 15}\\
 \texttt{ 39985367046248714: 15}

 \texttt{ 39996207167242544: 15}\\
 \texttt{ 40025017442723512: 15}\\
 \texttt{ 40042966430149288: 15}\\
 \texttt{ 40046501588099670: 15}\\
 \texttt{ 40113660053551766: 15}\\
 \texttt{ 40211500765180072: 15}\\
 \texttt{ 40254074476964216: 15}\\
 \texttt{ 40278834818305611: 15}\\
 \texttt{ 40584211742994462: 15}\\
 \texttt{ 40658224720242723: 15}\\
 \texttt{ 40676231408710112: 15}\\
 \texttt{ 40737707154564315: 15}\\
 \texttt{ 40794173778854166: 15}\\
 \texttt{ 41082714122308719: 15}\\
 \texttt{ 41105031108427422: 15}\\
 \texttt{ 41239488680485350: 15}\\
 \texttt{ 41242072094626840: 15}\\
 \texttt{ 41287030720393443: 15}\\
 \texttt{ 41325813280195314: 15}\\
 \texttt{ 41364247705928684: 15}

 \texttt{ 41540544034559215: 15}\\
 \texttt{ 41577383888947886: 15}\\
 \texttt{ 41718582639584142: 15}\\
 \texttt{ 41741409385977819: 15}\\
 \texttt{ 41771264218246043: 15}\\
 \texttt{ 41786923926601636: 15}\\
 \texttt{ 41801162245406442: 15}\\
 \texttt{ 41862893832335812: 15}\\
 \texttt{ 41918874237509372: 15}\\
 \texttt{ 41976465859566270: 15}\\
 \texttt{ 41987547331017272: 15}\\
 \texttt{ 41994755257110847: 15}\\
 \texttt{ 41994886796260218: 15}\\
 \texttt{ 42112665530003211: 15}\\
 \texttt{ 42125654316189111: 15}\\
 \texttt{ 42133297780621638: 15}\\
 \texttt{ 42398484182106811: 15}\\
 \texttt{ 42458993185515819: 15}\\
 \texttt{ 42575370339647348: 15}\\
 \texttt{ 42648782263932616: 15}

 \texttt{ 42765981683116362: 15}\\
 \texttt{ 42874925556367448: 15}\\
 \texttt{ 42898913473562143: 15}\\
 \texttt{ 42947799704572370: 15}\\
 \texttt{ 42954398320167520: 15}\\
 \texttt{ 42963278000732466: 15}\\
 \texttt{ 43102957685226268: 15}\\
 \texttt{ 43199995064686747: 15}\\
 \texttt{ 43243198919039815: 15}\\
 \texttt{ 43256610102020766: 15}\\
 \texttt{ 43363297944627174: 15}\\
 \texttt{ 43365605879716074: 15}\\
 \texttt{ 43515636851318020: 15}\\
 \texttt{ 43519634612796040: 15}\\
 \texttt{ 43565033665976336: 15}\\
 \texttt{ 43613234159148124: 15}\\
 \texttt{ 43683641503007764: 15}\\
 \texttt{ 43710658547111838: 15}\\
 \texttt{ 44032372793765644: 15}\\
 \texttt{ 44077855801604575: 15}

 \texttt{ 44089100776547716: 15}\\
 \texttt{ 44097841179043971: 15}\\
 \texttt{ 44152036721855620: 15}\\
 \texttt{ 44166642572534739: 15}\\
 \texttt{ 44185160217286844: 15}\\
 \texttt{ 44231647607101816: 15}\\
 \texttt{ 44258185136784546: 15}\\
 \texttt{ 44318937512286016: 15}\\
 \texttt{ 44356024815281214: 15}\\
 \texttt{ 44395871435792942: 15}\\
 \texttt{ 44441210940374572: 15}\\
 \texttt{ 44450778933960448: 15}\\
 \texttt{ 44578464997174820: 15}\\
 \texttt{ 44604882569241940: 15}\\
 \texttt{ 44660848984051146: 15}\\
 \texttt{ 44690941269552212: 15}\\
 \texttt{ 44779620545323719: 15}\\
 \texttt{ 44867420486511812: 15}\\
 \texttt{ 44907834326016571: 15}\\
 \texttt{ 44995771484250839: 15}

 \texttt{ 45003097306196824: 15}\\
 \texttt{ 45094076375835219: 15}\\
 \texttt{ 45106417470130614: 15}\\
 \texttt{ 45142707680856970: 15}\\
 \texttt{ 45197543497390036: 15}\\
 \texttt{ 45236776825376936: 15}\\
 \texttt{ 45296742255433718: 15}\\
 \texttt{ 45297438681527150: 15}\\
 \texttt{ 45334936394362912: 15}\\
 \texttt{ 45366332491482314: 15}\\
 \texttt{ 45424616871865147: 15}\\
 \texttt{ 45433454215802046: 15}\\
 \texttt{ 45563651545328512: 15}\\
 \texttt{ 45564281306214174: 15}\\
 \texttt{ 45656223343034948: 15}\\
 \texttt{ 45737374457364018: 15}\\
 \texttt{ 45759093207424170: 15}\\
 \texttt{ 45788565450314715: 15}\\
 \texttt{ 45928118208442470: 15}\\
 \texttt{ 46150403205368672: 15}

 \texttt{ 46251247686506114: 15}\\
 \texttt{ 46260555076043367: 15}\\
 \texttt{ 46261523154502314: 15}\\
 \texttt{ 46266705470024214: 15}\\
 \texttt{ 46268272623068862: 15}\\
 \texttt{ 46430738818317724: 15}\\
 \texttt{ 46443945749292364: 15}\\
 \texttt{ 46521953621107974: 15}\\
 \texttt{ 46527837194398624: 15}\\
 \texttt{ 46552954475165742: 15}\\
 \texttt{ 46666436046452572: 15}\\
 \texttt{ 46674273852990014: 15}\\
 \texttt{ 46689209965182543: 15}\\
 \texttt{ 46705368080549275: 15}\\
 \texttt{ 46714261083691974: 15}\\
 \texttt{ 46717595829767167: 16}\\
 \texttt{ 46818722469684175: 15}\\
 \texttt{ 46826401168029568: 15}\\
 \texttt{ 46848223239995164: 15}\\
 \texttt{ 46870383245259174: 15}

 \texttt{ 46897032169507540: 15}\\
 \texttt{ 46946348885058847: 15}\\
 \texttt{ 46952235251863064: 15}\\
 \texttt{ 46974552063007612: 15}\\
 \texttt{ 47111607700354624: 15}\\
 \texttt{ 47136233615475220: 15}\\
 \texttt{ 47213432969532343: 15}\\
 \texttt{ 47241155769148311: 15}\\
 \texttt{ 47312196262109511: 15}\\
 \texttt{ 47373745878270474: 15}\\
 \texttt{ 47400312139889022: 15}\\
 \texttt{ 47422612443528090: 15}\\
 \texttt{ 47473869213994847: 15}\\
 \texttt{ 47504780490883468: 15}\\
 \texttt{ 47578814235275272: 15}\\
 \texttt{ 47667112824818970: 15}\\
 \texttt{ 47671702947153411: 15}\\
 \texttt{ 47680561045966216: 15}\\
 \texttt{ 47746129339385442: 15}\\
 \texttt{ 47763280707177175: 15}

 \texttt{ 47821505388469540: 15}\\
 \texttt{ 47846234156181422: 15}\\
 \texttt{ 47928686106459464: 15}\\
 \texttt{ 47984977353514362: 15}\\
 \texttt{ 48185712609941824: 15}\\
 \texttt{ 48200518279591170: 15}\\
 \texttt{ 48209676624475143: 15}\\
 \texttt{ 48258735466329270: 15}\\
 \texttt{ 48269340085335340: 15}\\
 \texttt{ 48290732533084866: 15}\\
 \texttt{ 48483719416283738: 15}\\
 \texttt{ 48543818656792866: 15}\\
 \texttt{ 48640162692714664: 15}\\
 \texttt{ 48648490271487764: 15}\\
 \texttt{ 48651026117602147: 15}\\
 \texttt{ 48704884458190467: 15}\\
 \texttt{ 48727855380799518: 15}\\
 \texttt{ 48740330490714772: 15}\\
 \texttt{ 48772582754041310: 16}\\
 \texttt{ 49052553162537716: 15}

 \texttt{ 49063090272332630: 15}\\
 \texttt{ 49065002239587775: 15}\\
 \texttt{ 49137810966803168: 15}\\
 \texttt{ 49189122290278916: 15}\\
 \texttt{ 49192968988913442: 15}\\
 \texttt{ 49245830414985542: 15}\\
 \texttt{ 49284144880352118: 15}\\
 \texttt{ 49303396567331114: 15}\\
 \texttt{ 49428341049041863: 16}\\
 \texttt{ 49483757062318550: 15}\\
 \texttt{ 49495190008274043: 15}\\
 \texttt{ 49508097871686414: 15}\\
 \texttt{ 49579547251689842: 15}\\
 \texttt{ 49651625818595750: 15}\\
 \texttt{ 49657003800469142: 15}\\
 \texttt{ 49671066783761214: 15}\\
 \texttt{ 49852815533185648: 15}\\
 \texttt{ 49885646646543712: 15}\\
 \texttt{ 49905035615258622: 15}\\
 \texttt{ 49992577108634671: 15}

 \texttt{ 49994558050298523: 15}\\
 \texttt{ 50007507293532940: 15}\\
 \texttt{ 50010324944922950: 15}\\
 \texttt{ 50011847799468448: 17}\\
 \texttt{ 50130995740832668: 15}\\
 \texttt{ 50200287833699250: 15}\\
 \texttt{ 50246536452069114: 15}\\
 \texttt{ 50266689738341920: 15}\\
 \texttt{ 50271152007476642: 15}\\
 \texttt{ 50301811159421767: 15}\\
 \texttt{ 50374378394286240: 17}\\
 \texttt{ 50384012728202223: 15}\\
 \texttt{ 50436893115685923: 15}\\
 \texttt{ 50451606536644816: 15}\\
 \texttt{ 50469829117714520: 15}\\
 \texttt{ 50471978921413062: 15}\\
 \texttt{ 50609154636078818: 15}\\
 \texttt{ 50628741774093616: 15}\\
 \texttt{ 50645718279117115: 15}\\
 \texttt{ 50784620705204440: 15}

 \texttt{ 50860930696354516: 15}\\
 \texttt{ 50937162725245616: 15}\\
 \texttt{ 50966693034899272: 15}\\
 \texttt{ 50976006434453166: 15}\\
 \texttt{ 51025556457669616: 15}\\
 \texttt{ 51041039064675667: 15}\\
 \texttt{ 51075640231971236: 15}\\
 \texttt{ 51121707820830512: 15}\\
 \texttt{ 51124246855703144: 15}\\
 \texttt{ 51132983243769338: 15}\\
 \texttt{ 51141075376599520: 15}\\
 \texttt{ 51143602021113424: 15}\\
 \texttt{ 51217555299446718: 15}\\
 \texttt{ 51258560709511771: 15}\\
 \texttt{ 51264502597521616: 15}\\
 \texttt{ 51285870150281116: 15}\\
 \texttt{ 51291227536909674: 15}\\
 \texttt{ 51310011947937775: 15}\\
 \texttt{ 51342482904825616: 15}\\
 \texttt{ 51480625541339150: 15}

 \texttt{ 51515724838001224: 15}\\
 \texttt{ 51530294447940966: 15}\\
 \texttt{ 51534087307069623: 15}\\
 \texttt{ 51638923642698764: 15}\\
 \texttt{ 51666807949559620: 15}\\
 \texttt{ 51682711939053448: 15}\\
 \texttt{ 51728481943907738: 15}\\
 \texttt{ 51793596865070824: 15}\\
 \texttt{ 51804029744226691: 15}\\
 \texttt{ 51813687518360440: 15}\\
 \texttt{ 51913388402510019: 15}\\
 \texttt{ 52016503819446412: 15}\\
 \texttt{ 52062742469325075: 15}\\
 \texttt{ 52207877090443542: 15}\\
 \texttt{ 52212579525320418: 15}\\
 \texttt{ 52333276489323812: 15}\\
 \texttt{ 52610894794081239: 15}\\
 \texttt{ 52613862012903320: 15}\\
 \texttt{ 52682186279691522: 15}\\
 \texttt{ 52683061651453048: 15}

 \texttt{ 52703127452519442: 15}\\
 \texttt{ 52800091583360944: 15}\\
 \texttt{ 52806946957186660: 17}\\
 \texttt{ 52831033131184136: 15}\\
 \texttt{ 52948683164730318: 15}\\
 \texttt{ 53043310384120264: 15}\\
 \texttt{ 53055450644018767: 15}\\
 \texttt{ 53115244486628116: 15}\\
 \texttt{ 53116005982442874: 15}\\
 \texttt{ 53121000869839736: 15}\\
 \texttt{ 53171979430861915: 15}\\
 \texttt{ 53240371440827740: 15}\\
 \texttt{ 53249016307005519: 15}\\
 \texttt{ 53254490519872050: 15}\\
 \texttt{ 53272258482102315: 15}\\
 \texttt{ 53277738675398442: 15}\\
 \texttt{ 53359503084389320: 15}\\
 \texttt{ 53390246849888767: 15}\\
 \texttt{ 53415775590068563: 15}\\
 \texttt{ 53423879318062070: 15}

 \texttt{ 53426924506696442: 15}\\
 \texttt{ 53489747499263370: 15}\\
 \texttt{ 53599990680042844: 15}\\
 \texttt{ 53609373883407615: 15}\\
 \texttt{ 53825558769893524: 15}\\
 \texttt{ 53829362555976762: 15}\\
 \texttt{ 54023551524457240: 15}\\
 \texttt{ 54026414558777348: 15}\\
 \texttt{ 54082171897919946: 15}\\
 \texttt{ 54098160595448046: 15}\\
 \texttt{ 54149968212342664: 15}\\
 \texttt{ 54159712948497247: 15}\\
 \texttt{ 54298900060894412: 15}\\
 \texttt{ 54354504591428739: 15}\\
 \texttt{ 54426061592322172: 15}\\
 \texttt{ 54458000109034314: 15}\\
 \texttt{ 54477082515396571: 15}\\
 \texttt{ 54481157419373223: 15}\\
 \texttt{ 54580615826752422: 15}\\
 \texttt{ 54622748412555711: 15}

 \texttt{ 54699646449311968: 15}\\
 \texttt{ 54840102968767443: 15}\\
 \texttt{ 54845003384783043: 15}\\
 \texttt{ 54945549223917667: 15}\\
 \texttt{ 54952124369202944: 15}\\
 \texttt{ 55013231979275236: 15}\\
 \texttt{ 55039310568335610: 16}\\
 \texttt{ 55041049437122070: 15}\\
 \texttt{ 55050350646531570: 15}\\
 \texttt{ 55268311069387168: 15}\\
 \texttt{ 55276897875358516: 15}\\
 \texttt{ 55305969458439570: 15}\\
 \texttt{ 55307173409085474: 15}\\
 \texttt{ 55324067891176242: 15}\\
 \texttt{ 55366078821443466: 15}\\
 \texttt{ 55381771686304674: 15}\\
 \texttt{ 55397541227341336: 15}\\
 \texttt{ 55399550755147682: 15}\\
 \texttt{ 55408460658385866: 15}\\
 \texttt{ 55417900362838216: 15}

 \texttt{ 55450752106870412: 15}\\
 \texttt{ 55627266724331263: 15}\\
 \texttt{ 55641272682368743: 15}\\
 \texttt{ 55657028139036640: 15}\\
 \texttt{ 55668016068709443: 15}\\
 \texttt{ 55675692693011743: 15}\\
 \texttt{ 55684058445655314: 15}\\
 \texttt{ 55709816911946514: 15}\\
 \texttt{ 55769575544820915: 15}\\
 \texttt{ 55854873720082843: 15}\\
 \texttt{ 55886721710265122: 15}\\
 \texttt{ 55924044366826518: 15}\\
 \texttt{ 55957979208496144: 15}\\
 \texttt{ 56000955426423318: 15}\\
 \texttt{ 56008670442279570: 15}\\
 \texttt{ 56032473416608274: 15}\\
 \texttt{ 56177891552706592: 15}\\
 \texttt{ 56229019211708071: 15}\\
 \texttt{ 56261262636823770: 15}\\
 \texttt{ 56306717239347267: 15}

 \texttt{ 56354897875155023: 15}\\
 \texttt{ 56460825795781771: 15}\\
 \texttt{ 56524452893621215: 15}\\
 \texttt{ 56595999499097668: 15}\\
 \texttt{ 56636938801836616: 15}\\
 \texttt{ 56700152290660614: 15}\\
 \texttt{ 56830464839442320: 15}\\
 \texttt{ 56843448337955643: 15}\\
 \texttt{ 56850509272602474: 15}\\
 \texttt{ 56954939264751319: 15}\\
 \texttt{ 57025651911546536: 15}\\
 \texttt{ 57131578048600671: 15}\\
 \texttt{ 57234564997801242: 15}\\
 \texttt{ 57400052731390743: 15}\\
 \texttt{ 57446040607796815: 15}\\
 \texttt{ 57462426100311162: 15}\\
 \texttt{ 57562017961096746: 15}\\
 \texttt{ 57567757646305950: 15}\\
 \texttt{ 57568161423622362: 15}\\
 \texttt{ 57614677573412966: 15}

 \texttt{ 57653088133834519: 15}\\
 \texttt{ 57697464690660318: 15}\\
 \texttt{ 57714635747127975: 15}\\
 \texttt{ 57747774530799114: 15}\\
 \texttt{ 57807222269043212: 15}\\
 \texttt{ 57846289667276714: 15}\\
 \texttt{ 57873214014346312: 15}\\
 \texttt{ 57877154987121836: 15}\\
 \texttt{ 57976811566774467: 15}\\
 \texttt{ 58002516320502939: 15}\\
 \texttt{ 58015359029385196: 15}\\
 \texttt{ 58021648308306943: 15}\\
 \texttt{ 58042999008997036: 17}\\
 \texttt{ 58044833828993670: 15}\\
 \texttt{ 58050277958215375: 15}\\
 \texttt{ 58249819915428111: 15}\\
 \texttt{ 58293081955792663: 15}\\
 \texttt{ 58333761253039866: 15}\\
 \texttt{ 58362563896531024: 15}\\
 \texttt{ 58385723018640487: 15}

 \texttt{ 58426550323422848: 15}\\
 \texttt{ 58468277153121223: 15}\\
 \texttt{ 58511100456350360: 17}\\
 \texttt{ 58667010347335864: 15}\\
 \texttt{ 58734880592584218: 15}\\
 \texttt{ 58743231569322136: 15}\\
 \texttt{ 58752554645178870: 15}\\
 \texttt{ 58809549166284316: 15}\\
 \texttt{ 58810342105607012: 15}\\
 \texttt{ 58868977124812688: 15}\\
 \texttt{ 59104522776517468: 15}\\
 \texttt{ 59129105703789068: 15}\\
 \texttt{ 59132828502561318: 15}\\
 \texttt{ 59135771993725771: 15}\\
 \texttt{ 59139941114298942: 15}\\
 \texttt{ 59141227937625944: 15}\\
 \texttt{ 59168348646557444: 15}\\
 \texttt{ 59179406612798562: 15}\\
 \texttt{ 59220720401839312: 15}\\
 \texttt{ 59246863434106544: 15}

 \texttt{ 59304512808005872: 15}\\
 \texttt{ 59318844555472243: 15}\\
 \texttt{ 59388070175201272: 15}\\
 \texttt{ 59449612644271447: 15}\\
 \texttt{ 59490849381707668: 15}\\
 \texttt{ 59570832496630014: 15}\\
 \texttt{ 59609716403188622: 15}\\
 \texttt{ 59706995931058875: 15}\\
 \texttt{ 59722289239357062: 15}\\
 \texttt{ 59747514878060864: 15}\\
 \texttt{ 59799204036014018: 15}\\
 \texttt{ 59819571408895718: 15}\\
 \texttt{ 59840330433438470: 15}\\
 \texttt{ 59868813823527724: 15}\\
 \texttt{ 60031808667566042: 15}\\
 \texttt{ 60061595478122642: 15}\\
 \texttt{ 60083830279497570: 15}\\
 \texttt{ 60124909261145666: 15}\\
 \texttt{ 60311023031617748: 15}\\
 \texttt{ 60374247022186250: 15}

 \texttt{ 60409875885949015: 15}\\
 \texttt{ 60523971078075666: 15}\\
 \texttt{ 60637835946882714: 15}\\
 \texttt{ 60658614309672115: 15}\\
 \texttt{ 60689959617311250: 15}\\
 \texttt{ 60702428766110047: 15}\\
 \texttt{ 60728649750540268: 15}\\
 \texttt{ 60771798415156696: 15}\\
 \texttt{ 60819636169676671: 15}\\
 \texttt{ 60891580830225220: 15}\\
 \texttt{ 60949414267745812: 15}\\
 \texttt{ 60973823573108575: 15}\\
 \texttt{ 60999484771039119: 15}\\
 \texttt{ 61008136631600738: 15}\\
 \texttt{ 61096502858053542: 15}\\
 \texttt{ 61120406757901324: 15}\\
 \texttt{ 61259451066049472: 15}\\
 \texttt{ 61320408626839075: 15}\\
 \texttt{ 61364343061306743: 15}\\
 \texttt{ 61486005691206018: 15}

 \texttt{ 61499536265900572: 15}\\
 \texttt{ 61500639903869274: 15}\\
 \texttt{ 61515491777822744: 15}\\
 \texttt{ 61591842200514571: 15}\\
 \texttt{ 61655863231719043: 15}\\
 \texttt{ 61659519164989963: 15}\\
 \texttt{ 61754240718529720: 15}\\
 \texttt{ 61830140784313220: 15}\\
 \texttt{ 61864662560550848: 15}\\
 \texttt{ 61868861412675414: 15}\\
 \texttt{ 61877332218739911: 15}\\
 \texttt{ 61908547218579268: 15}\\
 \texttt{ 61945333905888367: 15}\\
 \texttt{ 61960234216377650: 15}\\
 \texttt{ 61994710263378570: 15}\\
 \texttt{ 62024149528795916: 15}\\
 \texttt{ 62098080372146175: 15}\\
 \texttt{ 62101919985182118: 15}\\
 \texttt{ 62102417356691214: 15}\\
 \texttt{ 62142332452425868: 15}

 \texttt{ 62151394208276970: 15}\\
 \texttt{ 62234285982032175: 15}\\
 \texttt{ 62245971382433418: 15}\\
 \texttt{ 62258067900696172: 15}\\
 \texttt{ 62324844708892119: 15}\\
 \texttt{ 62358954274169046: 15}\\
 \texttt{ 62494536371132848: 15}\\
 \texttt{ 62510410734324315: 15}\\
 \texttt{ 62672546619728514: 15}\\
 \texttt{ 62832754553033524: 15}\\
 \texttt{ 62882804962084867: 15}\\
 \texttt{ 62891511411098367: 15}\\
 \texttt{ 62897295747902012: 15}\\
 \texttt{ 62936684470196118: 15}\\
 \texttt{ 62967715754720048: 15}\\
 \texttt{ 62989055999032864: 15}\\
 \texttt{ 63055803065517472: 15}\\
 \texttt{ 63057303299988150: 16}\\
 \texttt{ 63106598255437050: 15}\\
 \texttt{ 63138810652880120: 15}

 \texttt{ 63240724828173968: 15}\\
 \texttt{ 63254008808639620: 15}\\
 \texttt{ 63292075337673843: 15}\\
 \texttt{ 63303330245461971: 15}\\
 \texttt{ 63330676923099474: 15}\\
 \texttt{ 63338930332592739: 15}\\
 \texttt{ 63506833147176020: 15}\\
 \texttt{ 63523602121157563: 15}\\
 \texttt{ 63562801721611220: 15}\\
 \texttt{ 63612115127485964: 15}\\
 \texttt{ 63649693258635068: 15}\\
 \texttt{ 63670665903735248: 15}\\
 \texttt{ 63687316084239042: 15}\\
 \texttt{ 63792675702595012: 15}\\
 \texttt{ 63826999867398748: 15}\\
 \texttt{ 63834418429205740: 15}\\
 \texttt{ 63883239542413466: 15}\\
 \texttt{ 63892312988179112: 15}\\
 \texttt{ 63934785094653940: 15}\\
 \texttt{ 63996787359185414: 15}

 \texttt{ 64009201946002467: 15}\\
 \texttt{ 64035499277922711: 15}\\
 \texttt{ 64136044898366774: 15}\\
 \texttt{ 64147337158215536: 15}\\
 \texttt{ 64185545363544512: 15}\\
 \texttt{ 64228754735675574: 15}\\
 \texttt{ 64263810872350816: 15}\\
 \texttt{ 64274590593597571: 15}\\
 \texttt{ 64287162889072035: 16}\\
 \texttt{ 64442519727188272: 15}\\
 \texttt{ 64450670144708348: 15}\\
 \texttt{ 64485964990197871: 15}\\
 \texttt{ 64489167225842666: 15}\\
 \texttt{ 64525373426376115: 15}\\
 \texttt{ 64553738714621439: 15}\\
 \texttt{ 64582757599319438: 15}\\
 \texttt{ 64631020540039468: 15}\\
 \texttt{ 64639057636145750: 15}\\
 \texttt{ 64733139900787914: 15}\\
 \texttt{ 64787651045229074: 15}

 \texttt{ 64945000069986436: 15}\\
 \texttt{ 64950148467042668: 15}\\
 \texttt{ 64955021578865467: 15}\\
 \texttt{ 65062039815754574: 15}\\
 \texttt{ 65191494685146343: 16}\\
 \texttt{ 65259195468701262: 15}\\
 \texttt{ 65259944219850666: 15}\\
 \texttt{ 65301179081146074: 15}\\
 \texttt{ 65355816628701150: 15}\\
 \texttt{ 65357142370157319: 15}\\
 \texttt{ 65376405191513823: 15}\\
 \texttt{ 65559810073393275: 15}\\
 \texttt{ 65594092715498320: 15}\\
 \texttt{ 65619125094307374: 15}\\
 \texttt{ 65689394214867566: 15}\\
 \texttt{ 65864475128807838: 15}\\
 \texttt{ 65925353150769434: 15}\\
 \texttt{ 65975161293002348: 15}\\
 \texttt{ 66085780947347814: 15}\\
 \texttt{ 66113300465177954: 15}

 \texttt{ 66436181744686143: 15}\\
 \texttt{ 66500765340058420: 15}\\
 \texttt{ 66584382865044416: 15}\\
 \texttt{ 66592296588966378: 15}\\
 \texttt{ 66636520784830314: 15}\\
 \texttt{ 66651714409827724: 15}\\
 \texttt{ 66680835681611320: 15}\\
 \texttt{ 66765132016515870: 15}\\
 \texttt{ 66766179135229214: 15}\\
 \texttt{ 66888100597644148: 15}\\
 \texttt{ 66910957050367170: 15}\\
 \texttt{ 66989479222816567: 15}\\
 \texttt{ 67066973091838950: 15}\\
 \texttt{ 67185320317421514: 15}\\
 \texttt{ 67197768753498018: 15}\\
 \texttt{ 67311189745330444: 15}\\
 \texttt{ 67355895488398611: 15}\\
 \texttt{ 67365892085809550: 15}\\
 \texttt{ 67377244521365612: 15}\\
 \texttt{ 67382268812526942: 15}

 \texttt{ 67427126944101940: 15}\\
 \texttt{ 67720322517294218: 15}\\
 \texttt{ 67741037600048611: 15}\\
 \texttt{ 67796230563442624: 15}\\
 \texttt{ 67815103757964820: 15}\\
 \texttt{ 67837063201568250: 15}\\
 \texttt{ 67846774711558616: 15}\\
 \texttt{ 67863068359647138: 15}\\
 \texttt{ 67897421419357016: 15}\\
 \texttt{ 67985681324844940: 15}\\
 \texttt{ 68087785327923571: 15}\\
 \texttt{ 68119473897900771: 15}\\
 \texttt{ 68146222385502738: 15}\\
 \texttt{ 68158070938035944: 15}\\
 \texttt{ 68310910135190815: 15}\\
 \texttt{ 68419028125353066: 15}\\
 \texttt{ 68499241765183166: 15}\\
 \texttt{ 68545652326083475: 15}\\
 \texttt{ 68552880776786812: 15}\\
 \texttt{ 68740025046912270: 15}

 \texttt{ 68791740125353322: 15}\\
 \texttt{ 68801600541655023: 15}\\
 \texttt{ 68870082063613974: 15}\\
 \texttt{ 68889701715890020: 15}\\
 \texttt{ 68905151387738970: 15}\\
 \texttt{ 68974533661410266: 15}\\
 \texttt{ 69027426748339820: 15}\\
 \texttt{ 69042083368428872: 15}\\
 \texttt{ 69048961970438275: 15}\\
 \texttt{ 69125816771272374: 15}\\
 \texttt{ 69148703948993118: 15}\\
 \texttt{ 69244903606700562: 15}\\
 \texttt{ 69314733178523116: 15}\\
 \texttt{ 69341562575974814: 15}\\
 \texttt{ 69373980523163646: 15}\\
 \texttt{ 69459956044406072: 15}\\
 \texttt{ 69664011359899911: 15}\\
 \texttt{ 69712564526922074: 15}\\
 \texttt{ 69764116618052222: 15}\\
 \texttt{ 69870942244448366: 15}

 \texttt{ 69923204733729174: 15}\\
 \texttt{ 69934206662330738: 15}\\
 \texttt{ 70098789997489767: 15}\\
 \texttt{ 70162571061484674: 15}\\
 \texttt{ 70164569022316368: 15}\\
 \texttt{ 70165921214193570: 15}\\
 \texttt{ 70197553877493870: 15}\\
 \texttt{ 70198652200856418: 15}\\
 \texttt{ 70241313269469111: 15}\\
 \texttt{ 70263084995589424: 15}\\
 \texttt{ 70311324937013475: 15}\\
 \texttt{ 70353698380496222: 15}\\
 \texttt{ 70379758143364362: 15}\\
 \texttt{ 70388181091003363: 15}\\
 \texttt{ 70415088258705639: 15}\\
 \texttt{ 70510888531744063: 15}\\
 \texttt{ 70514375477350616: 15}\\
 \texttt{ 70592478828069536: 15}\\
 \texttt{ 70661842511871472: 15}\\
 \texttt{ 70688916814126119: 15}

 \texttt{ 70753884157622138: 15}\\
 \texttt{ 70868963267771044: 15}\\
 \texttt{ 70923332925194120: 15}\\
 \texttt{ 70924031565683874: 15}\\
 \texttt{ 70948660359475372: 15}\\
 \texttt{ 70985413244700871: 15}\\
 \texttt{ 71030112805218039: 15}\\
 \texttt{ 71036049308878422: 15}\\
 \texttt{ 71039223554647219: 15}\\
 \texttt{ 71080931734411468: 15}\\
 \texttt{ 71137529980510624: 15}\\
 \texttt{ 71175159206473623: 15}\\
 \texttt{ 71256059376291738: 15}\\
 \texttt{ 71266959141290971: 15}\\
 \texttt{ 71328948441455068: 15}\\
 \texttt{ 71378340533424943: 15}\\
 \texttt{ 71420104746696346: 15}\\
 \texttt{ 71434062189325843: 15}\\
 \texttt{ 71451564346959664: 15}\\
 \texttt{ 71555039193896720: 15}

 \texttt{ 71785474741820714: 15}\\
 \texttt{ 71828253464081920: 15}\\
 \texttt{ 71853005484544518: 15}\\
 \texttt{ 71887822686059224: 15}\\
 \texttt{ 71961523044010314: 15}\\
 \texttt{ 72034182923882523: 15}\\
 \texttt{ 72071445980830219: 15}\\
 \texttt{ 72076086349430418: 15}\\
 \texttt{ 72079044618262362: 15}\\
 \texttt{ 72126967209963138: 15}\\
 \texttt{ 72196709726404663: 15}\\
 \texttt{ 72207655418618944: 15}\\
 \texttt{ 72269227218102016: 15}\\
 \texttt{ 72287444029292466: 15}\\
 \texttt{ 72299498475211972: 15}\\
 \texttt{ 72333448635307112: 15}\\
 \texttt{ 72392776826842274: 15}\\
 \texttt{ 72398686797930222: 15}\\
 \texttt{ 72403801848077012: 15}\\
 \texttt{ 72436396090059770: 15}

 \texttt{ 72481805570508844: 15}\\
 \texttt{ 72525246672526323: 15}\\
 \texttt{ 72561609885909543: 15}\\
 \texttt{ 72648315801713066: 15}\\
 \texttt{ 72671889286943248: 15}\\
 \texttt{ 72768129519054222: 15}\\
 \texttt{ 72827491656186267: 15}\\
 \texttt{ 72829587081826303: 15}\\
 \texttt{ 72830122019578419: 15}\\
 \texttt{ 72841122942416370: 15}\\
 \texttt{ 72984173046297450: 15}\\
 \texttt{ 72993074538994432: 15}\\
 \texttt{ 73215748756305774: 15}\\
 \texttt{ 73266551726579071: 15}\\
 \texttt{ 73304675547602740: 15}\\
 \texttt{ 73399979793504870: 15}\\
 \texttt{ 73439154820016944: 15}\\
 \texttt{ 73441210523643268: 15}\\
 \texttt{ 73554743370479019: 15}\\
 \texttt{ 73568437421281468: 15}

 \texttt{ 73616914054500570: 15}\\
 \texttt{ 73686037987640572: 15}\\
 \texttt{ 73688991459959743: 15}\\
 \texttt{ 73845764619475122: 15}\\
 \texttt{ 73884588909314372: 15}\\
 \texttt{ 73886580201025820: 15}\\
 \texttt{ 74126973652020366: 15}\\
 \texttt{ 74178247975701812: 15}\\
 \texttt{ 74250902742515862: 15}\\
 \texttt{ 74391382621433563: 15}\\
 \texttt{ 74400998525077364: 15}\\
 \texttt{ 74479268774961512: 15}\\
 \texttt{ 74522340968064174: 15}\\
 \texttt{ 74537770175845964: 15}\\
 \texttt{ 74650684057311940: 15}\\
 \texttt{ 74694016137939664: 15}\\
 \texttt{ 74703454027579419: 15}\\
 \texttt{ 74748953785739947: 15}\\
 \texttt{ 74876521278724168: 15}\\
 \texttt{ 74950996988289464: 15}

 \texttt{ 75023267446115250: 15}\\
 \texttt{ 75037403118037544: 15}\\
 \texttt{ 75045575076363870: 15}\\
 \texttt{ 75140430308677843: 15}\\
 \texttt{ 75227150281044222: 15}\\
 \texttt{ 75236767618864612: 15}\\
 \texttt{ 75282330357271050: 15}\\
 \texttt{ 75332560546323374: 15}\\
 \texttt{ 75339133836866439: 15}\\
 \texttt{ 75401264497397919: 15}\\
 \texttt{ 75403124264456647: 15}\\
 \texttt{ 75543425389246472: 15}\\
 \texttt{ 75572945945086615: 15}\\
 \texttt{ 75661377694582419: 15}\\
 \texttt{ 75725911009789070: 15}\\
 \texttt{ 75771748138771050: 15}\\
 \texttt{ 75776270114738822: 15}\\
 \texttt{ 75870193517751114: 15}\\
 \texttt{ 75906849847005868: 15}\\
 \texttt{ 75941165304926224: 15}

 \texttt{ 75950072906541243: 15}\\
 \texttt{ 75988630762062316: 15}\\
 \texttt{ 76039411928537824: 15}\\
 \texttt{ 76069675122458371: 15}\\
 \texttt{ 76126812061811120: 15}\\
 \texttt{ 76272780153298374: 15}\\
 \texttt{ 76379186797822768: 15}\\
 \texttt{ 76480043813734875: 15}\\
 \texttt{ 76502028204033919: 15}\\
 \texttt{ 76521099054045472: 15}\\
 \texttt{ 76526656236606074: 15}\\
 \texttt{ 76644376256466914: 15}\\
 \texttt{ 76671212398564766: 15}\\
 \texttt{ 76709812211100340: 15}\\
 \texttt{ 76771903621698520: 15}\\
 \texttt{ 76837167777764744: 15}\\
 \texttt{ 76867688625389872: 15}\\
 \texttt{ 76873166258274944: 15}\\
 \texttt{ 76929757136180619: 15}\\
 \texttt{ 76973712441907516: 15}

 \texttt{ 76975922974793874: 15}\\
 \texttt{ 76991179262218243: 15}\\
 \texttt{ 77063548006138075: 15}\\
 \texttt{ 77103629149744239: 15}\\
 \texttt{ 77123679678282536: 15}\\
 \texttt{ 77196442668079338: 15}\\
 \texttt{ 77203941977533370: 15}\\
 \texttt{ 77237886455090872: 15}\\
 \texttt{ 77272070099927047: 15}\\
 \texttt{ 77272199602826344: 15}\\
 \texttt{ 77306782608980716: 15}\\
 \texttt{ 77310685026667850: 15}\\
 \texttt{ 77379613017736472: 15}\\
 \texttt{ 77387705591365062: 15}\\
 \texttt{ 77545450119182120: 15}\\
 \texttt{ 77677210623243714: 15}\\
 \texttt{ 77716209008260311: 15}\\
 \texttt{ 77893297911700972: 15}\\
 \texttt{ 77920259890990615: 15}\\
 \texttt{ 77950846140689348: 15}

 \texttt{ 77985331739551516: 15}\\
 \texttt{ 78052055320428243: 15}\\
 \texttt{ 78125139578552416: 15}\\
 \texttt{ 78143473374444774: 15}\\
 \texttt{ 78188774924868162: 15}\\
 \texttt{ 78208510476703864: 15}\\
 \texttt{ 78275673603858536: 15}\\
 \texttt{ 78319752877486770: 15}\\
 \texttt{ 78384056643761944: 15}\\
 \texttt{ 78395804898521768: 15}\\
 \texttt{ 78408633527282414: 15}\\
 \texttt{ 78444302214210819: 15}\\
 \texttt{ 78586889425858372: 15}\\
 \texttt{ 78707600306956516: 15}\\
 \texttt{ 78715273675998715: 15}\\
 \texttt{ 78719095140521371: 15}\\
 \texttt{ 78727886594290062: 15}\\
 \texttt{ 78766394161069924: 15}\\
 \texttt{ 78829039313235871: 15}\\
 \texttt{ 78833469260029514: 15}

 \texttt{ 78910088196649219: 15}\\
 \texttt{ 78914207144032240: 15}\\
 \texttt{ 78970910234599311: 15}\\
 \texttt{ 78991400056073238: 15}\\
 \texttt{ 79001305672292936: 15}\\
 \texttt{ 79033947777062870: 15}\\
 \texttt{ 79080448714122470: 15}\\
 \texttt{ 79130183658016436: 15}\\
 \texttt{ 79132612264348838: 16}\\
 \texttt{ 79274114795974766: 15}\\
 \texttt{ 79310578413628744: 15}\\
 \texttt{ 79320560666575544: 15}\\
 \texttt{ 79338857895828546: 15}\\
 \texttt{ 79506876792173720: 15}\\
 \texttt{ 79578982382494823: 15}\\
 \texttt{ 79640613249330462: 15}\\
 \texttt{ 79651759763608971: 15}\\
 \texttt{ 79660579391621646: 15}\\
 \texttt{ 79712171276790064: 15}\\
 \texttt{ 79738292361950248: 15}

 \texttt{ 79761062925520468: 15}\\
 \texttt{ 79773265820747164: 15}\\
 \texttt{ 79778519954451611: 15}\\
 \texttt{ 79801124924366140: 15}\\
 \texttt{ 79920457776297016: 15}\\
 \texttt{ 79942856358426016: 15}\\
 \texttt{ 79964560142531224: 15}\\
 \texttt{ 79966220353272512: 15}\\
 \texttt{ 80050670566083375: 15}\\
 \texttt{ 80114395170458966: 15}\\
 \texttt{ 80122470314446671: 15}\\
 \texttt{ 80163259936372624: 15}\\
 \texttt{ 80253834769688740: 15}\\
 \texttt{ 80282881422984568: 15}\\
 \texttt{ 80284388518295114: 15}\\
 \texttt{ 80472995884872112: 15}\\
 \texttt{ 80473079045305022: 15}\\
 \texttt{ 80535858982751816: 15}\\
 \texttt{ 80574529505935540: 15}\\
 \texttt{ 80630606013479620: 15}

 \texttt{ 80741295655898515: 15}\\
 \texttt{ 80790555717228870: 15}\\
 \texttt{ 80824109182899162: 15}\\
 \texttt{ 80841513118764219: 15}\\
 \texttt{ 80845646929296714: 15}\\
 \texttt{ 80853128412819848: 15}\\
 \texttt{ 80871244910396368: 15}\\
 \texttt{ 80963107356659044: 15}\\
 \texttt{ 81021172063452350: 15}\\
 \texttt{ 81036065869361468: 15}\\
 \texttt{ 81046517308093016: 15}\\
 \texttt{ 81053516791218472: 15}\\
 \texttt{ 81093319562898764: 15}\\
 \texttt{ 81117039221058442: 15}\\
 \texttt{ 81133196646485384: 15}\\
 \texttt{ 81140884391339262: 15}\\
 \texttt{ 81195611577700648: 15}\\
 \texttt{ 81233208692983447: 15}\\
 \texttt{ 81243699645580146: 15}\\
 \texttt{ 81330501241881868: 15}

 \texttt{ 81340362957346072: 15}\\
 \texttt{ 81370941225317442: 15}\\
 \texttt{ 81487181373156844: 15}\\
 \texttt{ 81544437589465314: 15}\\
 \texttt{ 81718745656778574: 15}\\
 \texttt{ 81788832447639316: 15}\\
 \texttt{ 81797880023369250: 15}\\
 \texttt{ 81883049102280870: 15}\\
 \texttt{ 81884521819709936: 15}\\
 \texttt{ 81911112950891418: 15}\\
 \texttt{ 81922946441525512: 15}\\
 \texttt{ 82022277783189715: 15}\\
 \texttt{ 82085837945562842: 15}\\
 \texttt{ 82118390966166570: 15}\\
 \texttt{ 82242319608021546: 15}\\
 \texttt{ 82331042577308072: 15}\\
 \texttt{ 82383385887204268: 15}\\
 \texttt{ 82388837963719014: 15}\\
 \texttt{ 82417246358495414: 15}\\
 \texttt{ 82476285120268711: 15}

 \texttt{ 82486531694735044: 15}\\
 \texttt{ 82511046299056075: 15}\\
 \texttt{ 82601906747179074: 15}\\
 \texttt{ 82643231198169915: 15}\\
 \texttt{ 82692542628779523: 15}\\
 \texttt{ 82703522837349546: 15}\\
 \texttt{ 82729651804356068: 15}\\
 \texttt{ 82739213954746675: 15}\\
 \texttt{ 82771124727613543: 15}\\
 \texttt{ 82908742012651995: 15}\\
 \texttt{ 83054949545036672: 15}\\
 \texttt{ 83065669664138622: 15}\\
 \texttt{ 83323651821930364: 15}\\
 \texttt{ 83341020520815271: 15}\\
 \texttt{ 83483449691648416: 15}\\
 \texttt{ 83589663845753836: 15}\\
 \texttt{ 83673424311418266: 15}\\
 \texttt{ 83677879578856119: 15}\\
 \texttt{ 83825033576869915: 15}\\
 \texttt{ 83882378642902064: 15}

 \texttt{ 83927516647902018: 15}\\
 \texttt{ 83998468571486370: 15}\\
 \texttt{ 83998660178103340: 15}\\
 \texttt{ 84001821977511844: 15}\\
 \texttt{ 84009888510249212: 15}\\
 \texttt{ 84029697851057618: 15}\\
 \texttt{ 84046144118866564: 15}\\
 \texttt{ 84186769093306324: 15}\\
 \texttt{ 84395199346680040: 15}\\
 \texttt{ 84455915928039942: 15}\\
 \texttt{ 84462622835349942: 15}\\
 \texttt{ 84471001055602070: 15}\\
 \texttt{ 84523236998214320: 15}\\
 \texttt{ 84530347705111624: 15}\\
 \texttt{ 84598056262108772: 15}\\
 \texttt{ 84617478939043840: 15}\\
 \texttt{ 84670461353471018: 15}\\
 \texttt{ 84753648803283124: 15}\\
 \texttt{ 84798370926371319: 15}\\
 \texttt{ 84808302962435368: 15}

 \texttt{ 84828247602097868: 15}\\
 \texttt{ 84885177125882211: 15}\\
 \texttt{ 84954880056917822: 15}\\
 \texttt{ 84957760850602612: 15}\\
 \texttt{ 84982145042756347: 15}\\
 \texttt{ 85045807560525343: 15}\\
 \texttt{ 85062276153101824: 15}\\
 \texttt{ 85124833000026775: 15}\\
 \texttt{ 85131387622007648: 15}\\
 \texttt{ 85190074974274875: 15}\\
 \texttt{ 85226833568282515: 15}\\
 \texttt{ 85232347511290323: 15}\\
 \texttt{ 85267227337730468: 15}\\
 \texttt{ 85327250616866740: 15}\\
 \texttt{ 85343564112048872: 15}\\
 \texttt{ 85345711745492668: 15}\\
 \texttt{ 85403504274681967: 15}\\
 \texttt{ 85409607801859243: 15}\\
 \texttt{ 85438879330047342: 15}\\
 \texttt{ 85474647562249412: 15}

 \texttt{ 85515750570320942: 15}\\
 \texttt{ 85518108285636244: 15}\\
 \texttt{ 85641971505801448: 15}\\
 \texttt{ 85698882237104467: 15}\\
 \texttt{ 85705158022737847: 15}\\
 \texttt{ 85722802344150844: 15}\\
 \texttt{ 85771281686326216: 15}\\
 \texttt{ 85824905851817618: 15}\\
 \texttt{ 85887976244644012: 15}\\
 \texttt{ 85906911676642146: 15}\\
 \texttt{ 86007185174652568: 15}\\
 \texttt{ 86033788045695520: 15}\\
 \texttt{ 86093801931789943: 15}\\
 \texttt{ 86170356514623172: 15}\\
 \texttt{ 86171134280660943: 15}\\
 \texttt{ 86404137865328212: 15}\\
 \texttt{ 86586503337903246: 15}\\
 \texttt{ 86616104764185547: 15}\\
 \texttt{ 86622356220566120: 15}\\
 \texttt{ 86669989983414519: 15}

 \texttt{ 86775131513309512: 15}\\
 \texttt{ 86799725566675172: 15}\\
 \texttt{ 86802701766574912: 15}\\
 \texttt{ 86866683840724072: 15}\\
 \texttt{ 86936908121854566: 15}\\
 \texttt{ 87041376330802971: 15}\\
 \texttt{ 87119770778333368: 15}\\
 \texttt{ 87171162051973363: 15}\\
 \texttt{ 87256620324639618: 15}\\
 \texttt{ 87301939348619072: 15}\\
 \texttt{ 87357701245454223: 15}\\
 \texttt{ 87362646745012216: 15}\\
 \texttt{ 87551228022722574: 15}\\
 \texttt{ 87627878127138546: 15}\\
 \texttt{ 87632245950836368: 15}\\
 \texttt{ 87753117126190864: 15}\\
 \texttt{ 87857342668513372: 15}\\
 \texttt{ 87902356007527564: 15}\\
 \texttt{ 87926673590003564: 15}\\
 \texttt{ 88002863247648870: 15}

 \texttt{ 88099112224087540: 15}\\
 \texttt{ 88114480285053174: 15}\\
 \texttt{ 88159321624255315: 15}\\
 \texttt{ 88266774630022218: 15}\\
 \texttt{ 88351312201938114: 15}\\
 \texttt{ 88537353970596315: 15}\\
 \texttt{ 88570955496750246: 15}\\
 \texttt{ 88588580647506639: 15}\\
 \texttt{ 88750824954148674: 15}\\
 \texttt{ 88925087295914415: 15}\\
 \texttt{ 89059275811914416: 15}\\
 \texttt{ 89191625620242968: 15}\\
 \texttt{ 89200247127828124: 15}\\
 \texttt{ 89208629628654536: 15}\\
 \texttt{ 89224936222992199: 15}\\
 \texttt{ 89251700263599246: 15}\\
 \texttt{ 89321747864581166: 15}\\
 \texttt{ 89336627160003411: 15}\\
 \texttt{ 89446161335205067: 15}\\
 \texttt{ 89450050385020338: 15}

 \texttt{ 89454024832841468: 15}\\
 \texttt{ 89468848842450415: 15}\\
 \texttt{ 89482042720444218: 15}\\
 \texttt{ 89500838422960922: 15}\\
 \texttt{ 89532034728582474: 16}\\
 \texttt{ 89582397542648215: 15}\\
 \texttt{ 89632868533339144: 15}\\
 \texttt{ 89657338841784042: 15}\\
 \texttt{ 89700753521973470: 15}\\
 \texttt{ 89731395397935339: 15}\\
 \texttt{ 89735615301302464: 15}\\
 \texttt{ 90015479162231512: 15}\\
 \texttt{ 90034930452817516: 15}\\
 \texttt{ 90046757936723875: 15}\\
 \texttt{ 90100488670824246: 15}\\
 \texttt{ 90216255335044318: 15}\\
 \texttt{ 90227495446536750: 15}\\
 \texttt{ 90236293430539914: 15}\\
 \texttt{ 90342932135490918: 15}\\
 \texttt{ 90411210795397468: 15}

 \texttt{ 90521865508956424: 15}\\
 \texttt{ 90588581057127615: 15}\\
 \texttt{ 90620818422055467: 15}\\
 \texttt{ 90679631785851438: 15}\\
 \texttt{ 90681754124205364: 15}\\
 \texttt{ 90691489403092744: 15}\\
 \texttt{ 90697902842969250: 15}\\
 \texttt{ 90728497340262475: 15}\\
 \texttt{ 90824795340476236: 15}\\
 \texttt{ 90864176392665546: 15}\\
 \texttt{ 90925800722023155: 15}\\
 \texttt{ 91005904238896818: 15}\\
 \texttt{ 91052644049594367: 15}\\
 \texttt{ 91144463524687770: 15}\\
 \texttt{ 91190000635283318: 15}\\
 \texttt{ 91198976795242672: 15}\\
 \texttt{ 91338422032653272: 15}\\
 \texttt{ 91382147464469368: 15}\\
 \texttt{ 91448547291063547: 15}\\
 \texttt{ 91493281107736216: 15}

 \texttt{ 91505845789273672: 15}\\
 \texttt{ 91530445190597622: 15}\\
 \texttt{ 91563803477542448: 15}\\
 \texttt{ 91595550147902215: 15}\\
 \texttt{ 91878996068777067: 15}\\
 \texttt{ 91928340955985666: 15}\\
 \texttt{ 91967019661485724: 15}\\
 \texttt{ 92081643818517472: 15}\\
 \texttt{ 92326659040605570: 15}\\
 \texttt{ 92360834859868172: 15}\\
 \texttt{ 92395069704939111: 15}\\
 \texttt{ 92459147534019846: 15}\\
 \texttt{ 92631745123641638: 15}\\
 \texttt{ 92651508487205367: 15}\\
 \texttt{ 92774971157441416: 15}\\
 \texttt{ 92826575288112712: 15}\\
 \texttt{ 92985145895209444: 15}\\
 \texttt{ 93025732303864411: 15}\\
 \texttt{ 93083665797652219: 15}\\
 \texttt{ 93158902867286068: 15}

 \texttt{ 93376657725366064: 15}\\
 \texttt{ 93381908667821872: 15}\\
 \texttt{ 93389642677582348: 15}\\
 \texttt{ 93525941546634939: 15}\\
 \texttt{ 93586976981933011: 15}\\
 \texttt{ 93590819730173816: 15}\\
 \texttt{ 93610965897318067: 15}\\
 \texttt{ 93619639472259919: 15}\\
 \texttt{ 93647958568725472: 15}\\
 \texttt{ 93673204014664636: 15}\\
 \texttt{ 94030732829200442: 15}\\
 \texttt{ 94079410310211819: 15}\\
 \texttt{ 94083126733529871: 15}\\
 \texttt{ 94105347979036914: 15}\\
 \texttt{ 94135635521165347: 15}\\
 \texttt{ 94168587103488150: 15}\\
 \texttt{ 94206503762610998: 15}\\
 \texttt{ 94213398464277822: 15}\\
 \texttt{ 94238082471980572: 17}\\
 \texttt{ 94259770479522914: 15}

 \texttt{ 94274984973738940: 15}\\
 \texttt{ 94318747452905766: 15}\\
 \texttt{ 94331181742535418: 15}\\
 \texttt{ 94370720416130740: 15}\\
 \texttt{ 94380011766624067: 15}\\
 \texttt{ 94409589507063848: 15}\\
 \texttt{ 94430638739940450: 15}\\
 \texttt{ 94463037526023522: 15}\\
 \texttt{ 94479833701510146: 15}\\
 \texttt{ 94585615860507570: 15}\\
 \texttt{ 94597054520846343: 15}\\
 \texttt{ 94626280668938272: 15}\\
 \texttt{ 94716055445939475: 15}\\
 \texttt{ 94801616834745462: 15}\\
 \texttt{ 94815509837934039: 15}\\
 \texttt{ 94876057588922824: 15}\\
 \texttt{ 94936986266735822: 15}\\
 \texttt{ 94962914376400238: 15}\\
 \texttt{ 95025462122055614: 16}\\
 \texttt{ 95025486981367670: 15}

 \texttt{ 95029323570520072: 15}\\
 \texttt{ 95030399519454640: 15}\\
 \texttt{ 95058988856527616: 15}\\
 \texttt{ 95144640617193136: 15}\\
 \texttt{ 95218459451062038: 15}\\
 \texttt{ 95323152286051718: 15}\\
 \texttt{ 95343341746184046: 15}\\
 \texttt{ 95460708718066362: 15}\\
 \texttt{ 95468371913960474: 15}\\
 \texttt{ 95480546269994872: 15}\\
 \texttt{ 95493351763824315: 15}\\
 \texttt{ 95557327936823066: 15}\\
 \texttt{ 95766688451419568: 15}\\
 \texttt{ 95907483510140347: 15}\\
 \texttt{ 96024948662541067: 15}\\
 \texttt{ 96062461954112168: 15}\\
 \texttt{ 96070009576302342: 15}\\
 \texttt{ 96120283146056862: 15}\\
 \texttt{ 96196094107008642: 15}\\
 \texttt{ 96264256179810274: 15}

 \texttt{ 96317475626518315: 15}\\
 \texttt{ 96490555041259519: 15}\\
 \texttt{ 96504816442146171: 15}\\
 \texttt{ 96516399031061702: 15}\\
 \texttt{ 96634478102710864: 17}\\
 \texttt{ 96650671366627820: 15}\\
 \texttt{ 96650709411705200: 15}\\
 \texttt{ 96883242985782375: 15}\\
 \texttt{ 96896845699443420: 15}\\
 \texttt{ 96918220988503747: 15}\\
 \texttt{ 97019500356003942: 15}\\
 \texttt{ 97036422988160139: 15}\\
 \texttt{ 97084573205024740: 15}\\
 \texttt{ 97119407251982816: 15}\\
 \texttt{ 97212830726385115: 15}\\
 \texttt{ 97220851989651546: 15}\\
 \texttt{ 97316379891082672: 15}\\
 \texttt{ 97547280782515218: 15}\\
 \texttt{ 97586840586451614: 15}\\
 \texttt{ 97638015181766812: 15}

 \texttt{ 97639597347615548: 15}\\
 \texttt{ 97639820693567466: 15}\\
 \texttt{ 97656400930822720: 15}\\
 \texttt{ 97657803306979350: 15}\\
 \texttt{ 97730695259836719: 15}\\
 \texttt{ 97751725005839947: 15}\\
 \texttt{ 97772490880774540: 15}\\
 \texttt{ 97808916771905514: 15}\\
 \texttt{ 97822017539948419: 15}\\
 \texttt{ 97904167003041812: 15}\\
 \texttt{ 97976293652300814: 15}\\
 \texttt{ 98000034006702750: 15}\\
 \texttt{ 98094122041590014: 15}\\
 \texttt{ 98146802967943324: 15}\\
 \texttt{ 98148564390356836: 15}\\
 \texttt{ 98159126095570542: 15}\\
 \texttt{ 98162243654316786: 15}\\
 \texttt{ 98291473863440767: 15}\\
 \texttt{ 98365618599026919: 15}\\
 \texttt{ 98562460820538748: 15}

 \texttt{ 98632375874723347: 15}\\
 \texttt{ 98679447857390414: 15}\\
 \texttt{ 98684210344750612: 15}\\
 \texttt{ 98698011656858042: 15}\\
 \texttt{ 98830691468252943: 15}\\
 \texttt{ 98843765736965044: 15}\\
 \texttt{ 98862273984431346: 15}\\
 \texttt{ 98900672502719572: 15}\\
 \texttt{ 98969936462527119: 15}\\
 \texttt{ 99012201716047324: 15}\\
 \texttt{ 99013200425149218: 15}\\
 \texttt{ 99013665587649350: 15}\\
 \texttt{ 99080969055232216: 15}\\
 \texttt{ 99144136175539730: 15}\\
 \texttt{ 99153398083358967: 15}\\
 \texttt{ 99202504960798722: 15}\\
 \texttt{ 99316272608306870: 15}\\
 \texttt{ 99320555694178420: 15}\\
 \texttt{ 99332682900593548: 15}\\
 \texttt{ 99351225472676516: 15}

 \texttt{ 99361347526628971: 15}\\
 \texttt{ 99397893430359763: 15}\\
 \texttt{ 99578178292763174: 15}\\
 \texttt{ 99671989498331020: 15}\\
 \texttt{ 99681041294060344: 15}\\
 \texttt{ 99704325863238775: 15}\\
 \texttt{ 99758508008787016: 15}\\
 \texttt{ 99814131649982346: 15}\\
 \texttt{ 99885122419671940: 15}\\
 \texttt{ 99922640094450374: 15}\\
 \texttt{ 99943148937906044: 15}\\
 \texttt{ 99943992614301424: 15}\\
 \texttt{ 99954168293155540: 15}\\
 \texttt{100000813228553138: 16}\\
 \texttt{100023821443331464: 15}\\
 \texttt{100032470743935642: 15}\\
 \texttt{100051676864120668: 15}\\
 \texttt{100095893058102271: 15}\\
 \texttt{100236069553442270: 16}\\
 \texttt{100320096805665316: 15}

 \texttt{100342855933776736: 15}\\
 \texttt{100418427911006836: 15}\\
 \texttt{100539248490548646: 15}\\
 \texttt{100593433900509316: 15}\\
 \texttt{100599997794884262: 15}\\
 \texttt{100626691679688642: 15}\\
 \texttt{100648071097073874: 15}\\
 \texttt{100666895091758042: 15}\\
 \texttt{100778634267595964: 15}\\
 \texttt{100804783613385866: 15}\\
 \texttt{100814442996918044: 17}\\
 \texttt{100896421551965872: 15}\\
 \texttt{100918561435941571: 15}\\
 \texttt{100970055695135318: 15}\\
 \texttt{100997401648855362: 15}\\
 \texttt{101118155298947836: 15}\\
 \texttt{101208369718220344: 15}\\
 \texttt{101305539603939211: 15}\\
 \texttt{101364236922322418: 15}\\
 \texttt{101364367268289542: 15}

 \texttt{101372233031498366: 15}\\
 \texttt{101412779238631988: 15}\\
 \texttt{101510432264887672: 15}\\
 \texttt{101522620166847020: 15}\\
 \texttt{101525843706229411: 15}\\
 \texttt{101530253526347812: 15}\\
 \texttt{101576688351200619: 15}\\
 \texttt{101597484893026875: 15}\\
 \texttt{101624521527025143: 15}\\
 \texttt{101694920307422912: 15}\\
 \texttt{101893390436798864: 15}\\
 \texttt{101908694598425944: 15}\\
 \texttt{101958243384895719: 15}\\
 \texttt{102118776330076614: 15}\\
 \texttt{102165340081068374: 15}\\
 \texttt{102196902886622466: 15}\\
 \texttt{102198318298980774: 15}\\
 \texttt{102311197172138814: 15}\\
 \texttt{102317276279880844: 15}\\
 \texttt{102370514720260839: 15}

 \texttt{102406601269995711: 15}\\
 \texttt{102516194192516838: 15}\\
 \texttt{102517742742112167: 15}\\
 \texttt{102541634917868718: 15}\\
 \texttt{102594980251671666: 15}\\
 \texttt{102637026725713746: 15}\\
 \texttt{102678423403143460: 17}\\
 \texttt{102689921287660040: 15}\\
 \texttt{102693916581048914: 15}\\
 \texttt{102711959652336920: 15}\\
 \texttt{102712601782872512: 15}\\
 \texttt{102734633521323870: 15}\\
 \texttt{102964066777959568: 15}\\
 \texttt{103014361506177014: 15}\\
 \texttt{103020652769548974: 15}\\
 \texttt{103029128555110142: 15}\\
 \texttt{103059286192875038: 15}\\
 \texttt{103066178869271212: 15}\\
 \texttt{103330219105798312: 15}\\
 \texttt{103347413079108416: 15}

 \texttt{103481264344875567: 15}\\
 \texttt{103523950533890136: 17}\\
 \texttt{103590492060997267: 15}\\
 \texttt{103643871869093763: 16}\\
 \texttt{103701018058390948: 15}\\
 \texttt{103706073036340244: 15}\\
 \texttt{103719534278264942: 15}\\
 \texttt{103743967572719214: 15}\\
 \texttt{103770466881748719: 15}\\
 \texttt{103773045878464443: 15}\\
 \texttt{103812206658527814: 15}\\
 \texttt{103815436375145322: 15}\\
 \texttt{103831198909332016: 15}\\
 \texttt{103855585176213174: 15}\\
 \texttt{103865085864452067: 15}\\
 \texttt{103904487331114766: 15}\\
 \texttt{103984460223451564: 15}\\
 \texttt{104095122277303924: 15}\\
 \texttt{104109678579195770: 15}\\
 \texttt{104164883518978239: 15}

 \texttt{104174056232078118: 15}\\
 \texttt{104216316297665550: 15}\\
 \texttt{104381172991011316: 15}\\
 \texttt{104461955046045915: 15}\\
 \texttt{104473835464832612: 15}\\
 \texttt{104533047812271211: 15}\\
 \texttt{104547699761042574: 15}\\
 \texttt{104619754694694938: 15}\\
 \texttt{104646773018475412: 17}\\
 \texttt{104662091541563212: 15}\\
 \texttt{104715652882433274: 15}\\
 \texttt{104723647156103767: 15}\\
 \texttt{104768002166688618: 15}\\
 \texttt{104804333608257407: 15}\\
 \texttt{104814409411357672: 15}\\
 \texttt{104923316303272016: 15}\\
 \texttt{104977947296089744: 15}\\
 \texttt{105063400715677468: 15}\\
 \texttt{105087363206066368: 15}\\
 \texttt{105089409877031566: 15}

 \texttt{105102806578608712: 15}\\
 \texttt{105263654056097742: 15}\\
 \texttt{105398422577380743: 15}\\
 \texttt{105433995443054443: 15}\\
 \texttt{105469490813505570: 15}\\
 \texttt{105492284565429038: 15}\\
 \texttt{105575681161535912: 15}\\
 \texttt{105597491502528762: 15}\\
 \texttt{105786899718949022: 15}\\
 \texttt{105807012817708314: 15}\\
 \texttt{105892291708611772: 15}\\
 \texttt{105906004728990423: 15}\\
 \texttt{106165845703936314: 15}\\
 \texttt{106182851744008766: 15}\\
 \texttt{106251632201309120: 15}\\
 \texttt{106297825780786338: 15}\\
 \texttt{106398063529192916: 15}\\
 \texttt{106426841986012219: 15}\\
 \texttt{106632386892607912: 15}\\
 \texttt{106662257813815036: 15}

 \texttt{106729096330610212: 15}\\
 \texttt{106751366472701862: 15}\\
 \texttt{106806967957715812: 15}\\
 \texttt{106884395245368842: 15}\\
 \texttt{106908521297673716: 15}\\
 \texttt{106915107328578923: 15}\\
 \texttt{106982538080818767: 15}\\
 \texttt{107151896898766420: 15}\\
 \texttt{107187989690392070: 15}\\
 \texttt{107255052492057375: 15}\\
 \texttt{107272999815507116: 15}\\
 \texttt{107326922439628950: 15}\\
 \texttt{107382636330970012: 15}\\
 \texttt{107389112306983816: 15}\\
 \texttt{107416061701667618: 15}\\
 \texttt{107429784974726742: 15}\\
 \texttt{107502739803888416: 15}\\
 \texttt{107529934923347872: 15}\\
 \texttt{107599450754862736: 15}\\
 \texttt{107663625785504716: 15}

 \texttt{107663780754684618: 15}\\
 \texttt{107692589174465311: 15}\\
 \texttt{107695472253385508: 15}\\
 \texttt{107699895658905714: 15}\\
 \texttt{107843674859204270: 15}\\
 \texttt{107860564271941711: 15}\\
 \texttt{107916430370926144: 15}\\
 \texttt{107979224123671964: 15}\\
 \texttt{108135755066976674: 15}\\
 \texttt{108184120888710471: 15}\\
 \texttt{108185626651234443: 15}\\
 \texttt{108244072481128242: 15}\\
 \texttt{108259024970125672: 15}\\
 \texttt{108293215530262370: 15}\\
 \texttt{108410412955207614: 15}\\
 \texttt{108413819211999424: 15}\\
 \texttt{108431511787584942: 15}\\
 \texttt{108433057784238044: 15}\\
 \texttt{108500634407023722: 15}\\
 \texttt{108616379885415619: 15}

 \texttt{108709236995105574: 15}\\
 \texttt{108755331543650814: 15}\\
 \texttt{108767525479923474: 15}\\
 \texttt{108797128016875242: 15}\\
 \texttt{108962669055995072: 17}\\
 \texttt{109011814968604864: 15}\\
 \texttt{109074128128560763: 15}\\
 \texttt{109090665172141146: 15}\\
 \texttt{109121342162999120: 15}\\
 \texttt{109263621524597438: 15}\\
 \texttt{109273161941219311: 15}\\
 \texttt{109350883802229272: 15}\\
 \texttt{109393226514343012: 15}\\
 \texttt{109527053045858622: 15}\\
 \texttt{109564837294531540: 15}\\
 \texttt{109579263752321618: 15}\\
 \texttt{109646733685864770: 15}\\
 \texttt{109675835876035242: 15}\\
 \texttt{109681721267009371: 15}\\
 \texttt{109683271870397068: 15}

 \texttt{109726002618814516: 15}\\
 \texttt{109789861574349015: 15}\\
 \texttt{109832499946807072: 15}\\
 \texttt{109891582679176324: 15}\\
 \texttt{109952747319358922: 15}\\
 \texttt{109961356988350672: 15}\\
 \texttt{109997634817712743: 15}\\
 \texttt{110054179325925520: 15}\\
 \texttt{110081536957099768: 15}\\
 \texttt{110109575460420872: 15}\\
 \texttt{110131608263036118: 15}\\
 \texttt{110134285406483415: 15}\\
 \texttt{110134741592497648: 15}\\
 \texttt{110223560656580812: 15}\\
 \texttt{110256686048864319: 15}\\
 \texttt{110321950816709467: 15}\\
 \texttt{110322489195475315: 15}\\
 \texttt{110344202180543212: 15}\\
 \texttt{110416378622772247: 15}\\
 \texttt{110420353309365619: 15}

 \texttt{110463574886716219: 15}\\
 \texttt{110658715644468644: 15}\\
 \texttt{110697774353284111: 15}\\
 \texttt{110742061070114150: 15}\\
 \texttt{110746382434347064: 15}\\
 \texttt{110771754519522074: 15}\\
 \texttt{110828473170801570: 15}\\
 \texttt{110895091798932762: 15}\\
 \texttt{110974352592489964: 15}\\
 \texttt{111130539484757418: 15}\\
 \texttt{111155513455291766: 15}\\
 \texttt{111170464635251536: 15}\\
 \texttt{111172904642643172: 15}\\
 \texttt{111221964450689911: 15}\\
 \texttt{111368358938436664: 17}\\
 \texttt{111397206589426636: 15}\\
 \texttt{111435759141847542: 15}\\
 \texttt{111520733899962643: 15}\\
 \texttt{111560189004713643: 15}\\
 \texttt{111562590095265244: 15}

 \texttt{111580285672585419: 15}\\
 \texttt{111627560286557620: 15}\\
 \texttt{111667165494441811: 15}\\
 \texttt{111687297558729568: 15}\\
 \texttt{111707598881456872: 15}\\
 \texttt{111721930234233712: 15}\\
 \texttt{111767096720580424: 15}\\
 \texttt{111880394330731243: 15}\\
 \texttt{111959274337867743: 15}\\
 \texttt{112101745223810812: 15}\\
 \texttt{112133621643836418: 15}\\
 \texttt{112179908308922622: 15}\\
 \texttt{112290946891772044: 15}\\
 \texttt{112315783962010144: 15}\\
 \texttt{112327658901763538: 15}\\
 \texttt{112368628182953642: 15}\\
 \texttt{112421515996165868: 15}\\
 \texttt{112443803977121466: 15}\\
 \texttt{112452217548034436: 15}\\
 \texttt{112477931643567964: 15}

 \texttt{112491382839434511: 15}\\
 \texttt{112511388618417867: 15}\\
 \texttt{112658527415751674: 15}\\
 \texttt{112739711547992344: 15}\\
 \texttt{112831784308227039: 15}\\
 \texttt{112893593346543112: 15}\\
 \texttt{112946140120001019: 15}\\
 \texttt{113084212412452218: 15}\\
 \texttt{113097748406125372: 15}\\
 \texttt{113135073315893222: 15}\\
 \texttt{113147325720159171: 15}\\
 \texttt{113150615765156215: 15}\\
 \texttt{113157197915197442: 15}\\
 \texttt{113442512090382162: 15}\\
 \texttt{113501951174446539: 15}\\
 \texttt{113527628783097568: 15}\\
 \texttt{113579102355194174: 15}\\
 \texttt{113647622565945246: 15}\\
 \texttt{113690551667030468: 15}\\
 \texttt{113750383728616216: 15}

 \texttt{113887641344504415: 15}\\
 \texttt{113979023587742264: 15}\\
 \texttt{114034453796426022: 15}\\
 \texttt{114089732783564414: 15}\\
 \texttt{114254353322380239: 15}\\
 \texttt{114264519143219819: 15}\\
 \texttt{114277001562980144: 15}\\
 \texttt{114345416318261947: 15}\\
 \texttt{114384750820454542: 15}\\
 \texttt{114411034189707172: 15}\\
 \texttt{114498220465719616: 15}\\
 \texttt{114579433258681622: 15}\\
 \texttt{114715487932201275: 16}\\
 \texttt{114716465093830662: 16}\\
 \texttt{114717221893544048: 15}\\
 \texttt{114799001657037822: 15}\\
 \texttt{114931823225305118: 15}\\
 \texttt{114996885230657214: 15}\\
 \texttt{115104839374104342: 15}\\
 \texttt{115189366875716967: 15}

 \texttt{115323752744531114: 15}\\
 \texttt{115361021973186774: 15}\\
 \texttt{115378075772602611: 15}\\
 \texttt{115413232723259874: 15}\\
 \texttt{115420765916719744: 15}\\
 \texttt{115547424232662772: 15}\\
 \texttt{115558526371042772: 15}\\
 \texttt{115642288995030544: 15}\\
 \texttt{115646665885941542: 15}\\
 \texttt{115691139245699644: 15}\\
 \texttt{115729032057259542: 15}\\
 \texttt{115768388313886544: 15}\\
 \texttt{115827246872420666: 15}\\
 \texttt{115827337615830642: 15}\\
 \texttt{115959210437285438: 15}\\
 \texttt{116013889660337512: 15}\\
 \texttt{116099091630702724: 15}\\
 \texttt{116118709308891618: 15}\\
 \texttt{116210076308483875: 15}\\
 \texttt{116218775561206472: 15}

 \texttt{116270652444063711: 15}\\
 \texttt{116361699791039612: 15}\\
 \texttt{116372792921975116: 15}\\
 \texttt{116424967167164347: 15}\\
 \texttt{116429321207560875: 15}\\
 \texttt{116457180819643674: 15}\\
 \texttt{116465053318911619: 15}\\
 \texttt{116507731832747718: 15}\\
 \texttt{116566810912676671: 15}\\
 \texttt{116573486422756974: 15}\\
 \texttt{116642521555198170: 15}\\
 \texttt{116781717196556466: 15}\\
 \texttt{116790121422004688: 15}\\
 \texttt{116826892557810918: 16}\\
 \texttt{116940829757485162: 15}\\
 \texttt{117030152342848146: 15}\\
 \texttt{117229819774946920: 15}\\
 \texttt{117235830841332423: 15}\\
 \texttt{117327922216111166: 15}\\
 \texttt{117346868822512746: 15}

 \texttt{117352062005564211: 15}\\
 \texttt{117498696362508247: 15}\\
 \texttt{117523033359734319: 15}\\
 \texttt{117545082421375843: 15}\\
 \texttt{117565476268182162: 15}\\
 \texttt{117617476333258671: 15}\\
 \texttt{117627589174420912: 15}\\
 \texttt{117691580499155222: 15}\\
 \texttt{117799350479862642: 15}\\
 \texttt{117880767482636318: 15}\\
 \texttt{118012945467069519: 15}\\
 \texttt{118103948636028642: 15}\\
 \texttt{118178466640989304: 15}\\
 \texttt{118188517249239164: 15}\\
 \texttt{118213288269326072: 15}\\
 \texttt{118216166012078022: 15}\\
 \texttt{118293398541525075: 15}\\
 \texttt{118307380655581711: 15}\\
 \texttt{118320449396740820: 15}\\
 \texttt{118554532231453146: 15}

 \texttt{118573620478022416: 15}\\
 \texttt{118581540217515616: 15}\\
 \texttt{118627294254766315: 15}\\
 \texttt{118634776275760316: 15}\\
 \texttt{118739159598891714: 15}\\
 \texttt{118827588298061862: 15}\\
 \texttt{118946562385557067: 15}\\
 \texttt{119070535116486474: 15}\\
 \texttt{119098091696583219: 15}\\
 \texttt{119314231539726618: 15}\\
 \texttt{119465202071349162: 15}\\
 \texttt{119583126212143219: 15}\\
 \texttt{119671317068564212: 15}\\
 \texttt{119764898104257647: 15}\\
 \texttt{119770777333229272: 15}\\
 \texttt{119799358226364112: 15}\\
 \texttt{119843710963050643: 15}\\
 \texttt{119926942826523566: 15}\\
 \texttt{119937417861271867: 15}\\
 \texttt{119938286282491264: 15}

 \texttt{119967118433698119: 15}\\
 \texttt{120048016409678572: 15}\\
 \texttt{120057471711958371: 15}\\
 \texttt{120076615162014819: 15}\\
 \texttt{120112280341895275: 15}\\
 \texttt{120132256855682824: 15}\\
 \texttt{120326905951925418: 15}\\
 \texttt{120360973365758672: 15}\\
 \texttt{120392882746784574: 15}\\
 \texttt{120443476557706348: 15}\\
 \texttt{120545276545734666: 15}\\
 \texttt{120567471859616874: 15}\\
 \texttt{120572635974305740: 15}\\
 \texttt{120596140638268724: 15}\\
 \texttt{120627825291832144: 15}\\
 \texttt{120682615976030175: 15}\\
 \texttt{120811817621851024: 15}\\
 \texttt{120812053385633948: 15}\\
 \texttt{120860001191502414: 15}\\
 \texttt{120890166667268948: 15}

 \texttt{120946396870325142: 15}\\
 \texttt{120954084376945624: 15}\\
 \texttt{121041531181928850: 15}\\
 \texttt{121050358363613852: 15}\\
 \texttt{121065712829814570: 15}\\
 \texttt{121122613073037246: 15}\\
 \texttt{121126638708400623: 15}\\
 \texttt{121133528916588474: 15}\\
 \texttt{121145741278132623: 15}\\
 \texttt{121277365509005646: 15}\\
 \texttt{121289088553033470: 15}\\
 \texttt{121389253820883375: 15}\\
 \texttt{121414638659710324: 15}\\
 \texttt{121425134650377846: 15}\\
 \texttt{121492819729048514: 15}\\
 \texttt{121495501408032340: 15}\\
 \texttt{121502422386605066: 15}\\
 \texttt{121559626215118419: 15}\\
 \texttt{121600670623772044: 15}\\
 \texttt{121778362674689224: 15}

 \texttt{121835497996787523: 15}\\
 \texttt{121900822989082564: 15}\\
 \texttt{121919277441657938: 15}\\
 \texttt{121968538046406724: 15}\\
 \texttt{122080058987609166: 15}\\
 \texttt{122155120560178962: 15}\\
 \texttt{122157220815357818: 15}\\
 \texttt{122173009293172662: 15}\\
 \texttt{122192846292710164: 15}\\
 \texttt{122315917509313072: 15}\\
 \texttt{122368010661138066: 15}\\
 \texttt{122425791332841367: 15}\\
 \texttt{122436344351726862: 15}\\
 \texttt{122618736043573614: 15}\\
 \texttt{122627039929951770: 15}\\
 \texttt{122637371121313240: 15}\\
 \texttt{122695556811480412: 15}\\
 \texttt{122716671653445414: 15}\\
 \texttt{122860132958852618: 15}\\
 \texttt{122888197102869412: 15}

 \texttt{122907764317979750: 15}\\
 \texttt{123176876431552468: 15}\\
 \texttt{123186911979684942: 15}\\
 \texttt{123195196269357812: 15}\\
 \texttt{123219324956875122: 15}\\
 \texttt{123225884738206819: 15}\\
 \texttt{123262742676898671: 15}\\
 \texttt{123266182559377436: 15}\\
 \texttt{123302130479825872: 15}\\
 \texttt{123433799870545136: 15}\\
 \texttt{123444725054075838: 15}\\
 \texttt{123494713080933220: 15}\\
 \texttt{123498892600110871: 15}\\
 \texttt{123556872795044116: 15}\\
 \texttt{123595656080315336: 15}\\
 \texttt{123634902363616624: 15}\\
 \texttt{123643304599570143: 15}\\
 \texttt{123674229588415550: 15}\\
 \texttt{123810755202623167: 15}\\
 \texttt{123827528821983602: 15}

 \texttt{123844052283936163: 15}\\
 \texttt{123954615524721572: 15}\\
 \texttt{123971186548964644: 15}\\
 \texttt{124024774568393140: 15}\\
 \texttt{124176108022485914: 15}\\
 \texttt{124215301213158424: 15}\\
 \texttt{124233783311704218: 15}\\
 \texttt{124242222886830548: 15}\\
 \texttt{124248776689872411: 15}\\
 \texttt{124248796065710943: 15}\\
 \texttt{124282917499198623: 15}\\
 \texttt{124317399137390646: 15}\\
 \texttt{124360716197171275: 15}\\
 \texttt{124397924464578211: 15}\\
 \texttt{124424079783094815: 15}\\
 \texttt{124486741719932955: 15}\\
 \texttt{124647862448970044: 15}\\
 \texttt{124661525125769863: 15}\\
 \texttt{124674918801604567: 15}\\
 \texttt{124705785294943719: 15}

 \texttt{124729893242091338: 15}\\
 \texttt{124811135639651563: 15}\\
 \texttt{124847723677822314: 15}\\
 \texttt{124876506571459718: 15}\\
 \texttt{124888477425603822: 15}\\
 \texttt{124909260301793416: 15}\\
 \texttt{124956687692632975: 15}\\
 \texttt{125036985334204023: 15}\\
 \texttt{125070186114903075: 15}\\
 \texttt{125079967166853174: 15}\\
 \texttt{125277360246767740: 15}\\
 \texttt{125308438781944011: 15}\\
 \texttt{125309709390430672: 15}\\
 \texttt{125388986111190762: 15}\\
 \texttt{125498674295222911: 15}\\
 \texttt{125597635383387968: 15}\\
 \texttt{125624172891181074: 15}\\
 \texttt{125637121664242946: 15}\\
 \texttt{125645206950685375: 15}\\
 \texttt{125674125693775758: 15}

 \texttt{125735119325648368: 15}\\
 \texttt{125777823654324211: 16}\\
 \texttt{125781000834058568: 18}\\
 \texttt{125800712844871120: 15}\\
 \texttt{125830091489745614: 15}\\
 \texttt{125861098929212524: 15}\\
 \texttt{125868726008591522: 15}
\end{flushright}
}
\normalsize
\end{multicols}

\clearpage

\section*{Appendix D}

  Square-free gaps and their length $\geq 16$, up to $125\,870\,000\,000\,000\,000$

\begin{multicols}{4}
{
\footnotesize
\begin{flushright}
\noindent
 \texttt{  3215226335143218: 16}\\
 \texttt{ 23742453640900972: 17}\\
 \texttt{ 28696958943616635: 16}\\
 \texttt{ 31401976920688950: 16}\\
 \texttt{ 36985881099122836: 17}\\
 \texttt{ 46717595829767167: 16}\\
 \texttt{ 48772582754041310: 16}\\
 \texttt{ 49428341049041863: 16}\\
 \texttt{ 50011847799468448: 17}\\
 \texttt{ 50374378394286240: 17}\\
 \texttt{ 52806946957186660: 17}\\
 \texttt{ 55039310568335610: 16}\\
 \texttt{ 58042999008997036: 17}\\
 \texttt{ 58511100456350360: 17}\\
 \texttt{ 63057303299988150: 16}\\
 \texttt{ 64287162889072035: 16}\\
 \texttt{ 65191494685146343: 16}\\
 \texttt{ 79132612264348838: 16}\\
 \texttt{ 89532034728582474: 16}\\
 \texttt{ 94238082471980572: 17}
 
 \texttt{ 95025462122055614: 16}\\
 \texttt{ 96634478102710864: 17}\\
 \texttt{100000813228553138: 16}\\
 \texttt{100236069553442270: 16}\\
 \texttt{100814442996918044: 17}\\
 \texttt{102678423403143460: 17}\\
 \texttt{103523950533890136: 17}\\
 \texttt{103643871869093763: 16}\\
 \texttt{104646773018475412: 17}\\
 \texttt{108962669055995072: 17}\\
 \texttt{111368358938436664: 17}\\
 \texttt{114715487932201275: 16}\\
 \texttt{114716465093830662: 16}\\
 \texttt{116826892557810918: 16}\\
 \texttt{125777823654324211: 16}\\
 \texttt{125781000834058568: 18}
\end{flushright}
}
\normalsize
\end{multicols}

\clearpage

\section*{Appendix E}

Completed ranges computed on 44 different processors by 14 different users.

\begin{displaymath}
  \begin{array}{llrr}
  \text{User Name}& \text{(Computer Name)}& Range& Date Completed\\
\hline
  \text{L. Marmet}& \text{(neurone5 Linux .8GHz)}& 125.7\times 10^{15} \text{ to } 125.8\times 10^{15}& \text{September $12^{th}$, 2005}\\
  \text{H. Gassmann}& \text{(Gus's beast)}& 125.3\times 10^{15} \text{ to } 125.7\times 10^{15}& \text{August $19^{th}$, 2005}\\
  \text{L. Marmet}& \text{(neurone3 .8GHz)}& 125.1\times 10^{15} \text{ to } 125.3\times 10^{15}& \text{November $28^{th}$, 2005}\\
  \text{N. Marmet}& \text{(Poisson.51)}& 124.7\times 10^{15} \text{ to } 125.1\times 10^{15}& \text{August $19^{th}$, 2005}\\
  \text{N. Marmet}& \text{(Poisson.35)}& 124.4\times 10^{15} \text{ to } 124.7\times 10^{15}& \text{August $19^{th}$, 2005}\\
  \text{N. Marmet}& \text{(Poisson.141)}& 124.1\times 10^{15} \text{ to } 124.4\times 10^{15}& \text{August $19^{th}$, 2005}\\
  \text{N. Marmet}& \text{(Poisson.142.Brain1)}& 123.8\times 10^{15} \text{ to } 124.1\times 10^{15}& \text{August $19^{th}$, 2005}\\
  \text{N. Robertson}& \text{(Laika.Droica)}& 123.4\times 10^{15} \text{ to } 123.8\times 10^{15}& \text{August $20^{th}$, 2005}\\
  \text{L. Marmet}& \text{(neurone2 .8GHz)}& 122.7\times 10^{15} \text{ to } 123.4\times 10^{15}& \text{August $10^{th}$, 2005}\\
  \text{L. Marmet}& \text{(neurone1 .8GHz)}& 122.1\times 10^{15} \text{ to } 122.7\times 10^{15}& \text{August $7^{th}$, 2005}\\
  \text{N. Marmet}& \text{(Poisson.94)}& 121.6\times 10^{15} \text{ to } 122.1\times 10^{15}& \text{August $19^{th}$, 2005}\\
  \text{N. Marmet}& \text{(Poisson.63)}& 121.3\times 10^{15} \text{ to } 121.6\times 10^{15}& \text{August $19^{th}$, 2005}\\
  \text{L. Marmet}& \text{(neurone5 Linux .8GHz)}& 120.9\times 10^{15} \text{ to } 121.3\times 10^{15}& \text{July $26^{th}$, 2005}\\
  \text{N. Robertson}& \text{(Laika.Droica)}& 120.5\times 10^{15} \text{ to } 120.9\times 10^{15}& \text{July $16^{th}$, 2005}\\
  \text{L. Marmet}& \text{(neurone4 .8GHz)}& 119.3\times 10^{15} \text{ to } 120.5\times 10^{15}& \text{August $13^{th}$, 2005}\\
  \text{N. Marmet}& \text{(Poisson.142.Brain2)}& 118.0\times 10^{15} \text{ to } 119.3\times 10^{15}& \text{September $8^{th}$, 2005}\\
  \text{N. Marmet}& \text{(Poisson.107)}& 116.3\times 10^{15} \text{ to } 118.0\times 10^{15}& \text{September $2^{nd}$, 2005}\\
  \text{N. Marmet}& \text{(Poisson.31)}& 115.7\times 10^{15} \text{ to } 116.3\times 10^{15}& \text{August $29^{th}$, 2005}\\
  \text{N. Marmet}& \text{(Poisson.30)}& 114.9\times 10^{15} \text{ to } 115.7\times 10^{15}& \text{August $28^{th}$, 2005}\\
  \text{J.-P. Bernier}& \text{(Pentium 600MHz)}& 114.0\times 10^{15} \text{ to } 114.9\times 10^{15}& \text{August $24^{th}$, 2005}\\
  \text{L. Marmet}& \text{(neurone0 W95 .233GHz)}& 113.5\times 10^{15} \text{ to } 114.0\times 10^{15}& \text{August $16^{th}$, 2005}\\
  \text{L. Marmet}& \text{(neurone5 Linux .8GHz)}& 113.1\times 10^{15} \text{ to } 113.5\times 10^{15}& \text{June $30^{th}$, 2005}\\
  \text{N. Robertson}& \text{(Laika.Droica)}& 112.7\times 10^{15} \text{ to } 113.1\times 10^{15}& \text{June $24^{th}$, 2005}\\
  \text{N. Marmet}& \text{(Poisson.29)}& 111.5\times 10^{15} \text{ to } 112.7\times 10^{15}& \text{August $28^{th}$, 2005}\\
  \text{H. Gassmann}& \text{(Gus's beast)}& 111.0\times 10^{15} \text{ to } 111.5\times 10^{15}& \text{July $27^{th}$, 2005}\\
  \text{N. Marmet}& \text{(Poisson.44)}& 110.4\times 10^{15} \text{ to } 111.0\times 10^{15}& \text{August $24^{th}$, 2005}\\
  \text{N. Marmet}& \text{(Poisson.37)}& 109.6\times 10^{15} \text{ to } 110.4\times 10^{15}& \text{September $1^{st}$, 2005}\\
  \text{N. Marmet}& \text{(Poisson.34)}& 108.8\times 10^{15} \text{ to } 109.6\times 10^{15}& \text{September $3^{rd}$, 2005}\\
  \text{J.-P. Bernier}& \text{(Athlon 2000)}& 106.3\times 10^{15} \text{ to } 108.8\times 10^{15}& \text{August $5^{th}$, 2005}\\
  \text{B. Le Tual}& \text{(Celeron 2.4GHz)}& 105.3\times 10^{15} \text{ to } 106.3\times 10^{15}& \text{July $16^{th}$, 2005}\\
  \end{array}
\end{displaymath}

\begin{displaymath}
  \begin{array}{llrr}
  \text{User Name}& \text{(Computer Name)}& Range& Date Completed\\
\hline
  \text{B. Le Tual}& \text{(Celeron .9GHz)}& 104.9\times 10^{15} \text{ to } 105.3\times 10^{15}& \text{July $16^{th}$, 2005}\\
  \text{L. Marmet}& \text{(neurone2 .8GHz)}& 103.2\times 10^{15} \text{ to } 104.9\times 10^{15}& \text{July $8^{th}$, 2005}\\
  \text{B. Le Tual}& \text{(Celeron 2.4GHz)}& 102.7\times 10^{15} \text{ to } 103.2\times 10^{15}& \text{May $13^{th}$, 2005}\\
  \text{L. Marmet}& \text{(neurone4 .8GHz)}& 101.1\times 10^{15} \text{ to } 102.7\times 10^{15}& \text{June $25^{th}$, 2005}\\
  \text{L. Marmet}& \text{(neurone1 .8GHz)}& 99.5\times 10^{15} \text{ to } 101.1\times 10^{15}& \text{July $2^{nd}$, 2005}\\
  \text{L. Marmet}& \text{(Poisson.63)}& 98.8\times 10^{15} \text{ to } 99.5\times 10^{15}& \text{July $12^{th}$, 2005}\\
  \text{H. Gassmann}& \text{(Gus's beast)}& 98.4\times 10^{15} \text{ to } 98.8\times 10^{15}& \text{June $1^{st}$, 2005}\\
  \text{N. Robertson}& \text{(Laika.Droica)}& 98.0\times 10^{15} \text{ to } 98.4\times 10^{15}& \text{May $9^{th}$, 2005}\\
  \text{J.-P. Bernier}& \text{(Athlon 2000)}& 96.2\times 10^{15} \text{ to } 98.0\times 10^{15}& \text{May $24^{th}$, 2005}\\
  \text{L. Marmet}& \text{(neurone0 W95 .233GHz)}& 95.7\times 10^{15} \text{ to } 96.2\times 10^{15}& \text{June $10^{th}$, 2005}\\
  \text{N. Marmet}& \text{(Poisson.142.Brain2)}& 94.2\times 10^{15} \text{ to } 95.7\times 10^{15}& \text{June $13^{th}$, 2005}\\
  \text{N. Marmet}& \text{(Poisson.94)}& 92.7\times 10^{15} \text{ to } 94.2\times 10^{15}& \text{July $12^{th}$, 2005}\\
  \text{J.-P. Bernier}& \text{(Pentium 600MHz)}& 91.7\times 10^{15} \text{ to } 92.7\times 10^{15}& \text{June $12^{th}$, 2005}\\
  \text{N. Marmet}& \text{(Poisson.142.Brain1)}& 89.7\times 10^{15} \text{ to } 91.7\times 10^{15}& \text{July $21^{st}$, 2005}\\
  \text{N. Marmet}& \text{(Poisson.34)}& 88.9\times 10^{15} \text{ to } 89.7\times 10^{15}& \text{May $30^{th}$, 2005}\\
  \text{N. Marmet}& \text{(Poisson.107)}& 86.8\times 10^{15} \text{ to } 88.9\times 10^{15}& \text{June $13^{th}$, 2005}\\
  \text{N. Marmet}& \text{(Poisson.29)}& 85.5\times 10^{15} \text{ to } 86.8\times 10^{15}& \text{June $2^{nd}$, 2005}\\
  \text{N. Marmet}& \text{(Poisson.31)}& 84.7\times 10^{15} \text{ to } 85.5\times 10^{15}& \text{June $13^{th}$, 2005}\\
  \text{N. Marmet}& \text{(Poisson.51)}& 82.4\times 10^{15} \text{ to } 84.7\times 10^{15}& \text{July $25^{th}$, 2005}\\
  \text{N. Marmet}& \text{(Poisson.30)}& 81.4\times 10^{15} \text{ to } 82.4\times 10^{15}& \text{May $13^{th}$, 2005}\\
  \text{N. Marmet}& \text{(Poisson.45)}& 80.8\times 10^{15} \text{ to } 81.4\times 10^{15}& \text{May $30^{th}$, 2005}\\
  \text{N. Marmet}& \text{(Poisson.37)}& 80.0\times 10^{15} \text{ to } 80.8\times 10^{15}& \text{May $30^{th}$, 2005}\\
  \text{L. Marmet}& \text{(neurone5 Linux .8GHz)}& 79.2\times 10^{15} \text{ to } 80.0\times 10^{15}& \text{June $9^{th}$, 2005}\\
  \text{L. Marmet}& \text{(neurone5 Linux .8GHz)}& 78.8\times 10^{15} \text{ to } 79.2\times 10^{15}& \text{April $21^{st}$, 2005}\\
  \text{J.-P. Bernier}& \text{(Athlon 2000)}& 78.1\times 10^{15} \text{ to } 78.8\times 10^{15}& \text{April $9^{th}$, 2005}\\
  \text{N. Marmet}& \text{(Poisson.94)}& 77.6\times 10^{15} \text{ to } 78.1\times 10^{15}& \text{April $5^{th}$, 2005}\\
  \text{N. Marmet}& \text{(Poisson.142.Brain2)}& 77.1\times 10^{15} \text{ to } 77.6\times 10^{15}& \text{April $5^{th}$, 2005}\\
  \text{N. Marmet}& \text{(Poisson.142.Brain1)}& 76.6\times 10^{15} \text{ to } 77.1\times 10^{15}& \text{March $31^{st}$, 2005}\\
  \text{N. Marmet}& \text{(Poisson.107)}& 76.1\times 10^{15} \text{ to } 76.6\times 10^{15}& \text{March $30^{th}$, 2005}\\
  \text{N. Marmet}& \text{(Poisson.51)}& 75.6\times 10^{15} \text{ to } 76.1\times 10^{15}& \text{March $29^{th}$, 2005}\\
  \text{N. Marmet}& \text{(Poisson.63)}& 75.1\times 10^{15} \text{ to } 75.6\times 10^{15}& \text{April $26^{th}$, 2005}\\
  \text{N. Marmet}& \text{(Poisson.30)}& 74.9\times 10^{15} \text{ to } 75.1\times 10^{15}& \text{March $28^{th}$, 2005}\\
  \text{N. Marmet}& \text{(Poisson.31)}& 74.7\times 10^{15} \text{ to } 74.9\times 10^{15}& \text{March $30^{th}$, 2005}\\
  \text{N. Marmet}& \text{(Poisson.34)}& 74.5\times 10^{15} \text{ to } 74.7\times 10^{15}& \text{March $30^{th}$, 2005}\\
  \text{N. Marmet}& \text{(Poisson.37)}& 74.4\times 10^{15} \text{ to } 74.5\times 10^{15}& \text{March $23^{rd}$, 2005}\\
  \text{N. Marmet}& \text{(Poisson.45)}& 74.3\times 10^{15} \text{ to } 74.4\times 10^{15}& \text{March $25^{th}$, 2005}\\
  \text{N. Robertson}& \text{(Laika.Droica)}& 73.8\times 10^{15} \text{ to } 74.3\times 10^{15}& \text{April $5^{th}$, 2005}\\
  \text{H. Gassmann}& \text{(Gus's beast)}& 73.5\times 10^{15} \text{ to } 73.8\times 10^{15}& \text{April $27^{th}$, 2005}\\
  \text{L. Marmet}& \text{(neurone4 .8GHz)}& 71.9\times 10^{15} \text{ to } 73.5\times 10^{15}& \text{April $30^{th}$, 2005}\\
  \text{L. Marmet}& \text{(neurone2 .8GHz)}& 70.2\times 10^{15} \text{ to } 71.9\times 10^{15}& \text{May $4^{th}$, 2005}\\
  \end{array}
\end{displaymath}

\begin{displaymath}
  \begin{array}{llrr}
  \text{User Name}& \text{(Computer Name)}& Range& Date Completed\\
\hline
  \text{L. Marmet}& \text{(neurone1 .8GHz)}& 68.7\times 10^{15} \text{ to } 70.2\times 10^{15}& \text{April $28^{th}$, 2005}\\
  \text{J.-P. Bernier}& \text{(Pentium 600MHz)}& 68.4\times 10^{15} \text{ to } 68.7\times 10^{15}& \text{April $3^{rd}$, 2005}\\
  \text{B. Le Tual}& \text{(Celeron .9GHz)}& 68.0\times 10^{15} \text{ to } 68.4\times 10^{15}& \text{May $15^{th}$, 2005}\\
  \text{B. Le Tual}& \text{(Celeron 2.4GHz)}& 67.6\times 10^{15} \text{ to } 68.0\times 10^{15}& \text{May $1^{st}$, 2005}\\
  \text{J.-P. Bernier}& \text{(Athlon 2000)}& 67.2\times 10^{15} \text{ to } 67.6\times 10^{15}& \text{March $17^{th}$, 2005}\\
  \text{N. Marmet}& \text{(Poisson.29)}& 66.9\times 10^{15} \text{ to } 67.2\times 10^{15}& \text{March $30^{th}$, 2005}\\
  \text{J.-P. Bernier}& \text{(Athlon 2000)}& 66.5\times 10^{15} \text{ to } 66.9\times 10^{15}& \text{January $1^{st}$, 2004}\\
  \text{J.-P. Bernier}& \text{(Pentium 600MHz)}& 66.2\times 10^{15} \text{ to } 66.5\times 10^{15}& \text{January $12^{th}$, 2004}\\
  \text{B. Le Tual}& \text{(Celeron 2.4GHz)}& 65.9\times 10^{15} \text{ to } 66.2\times 10^{15}& \text{January $4^{th}$, 2004}\\
  \text{B. Le Tual}& \text{(Celeron .9GHz)}& 65.5\times 10^{15} \text{ to } 65.9\times 10^{15}& \text{January $9^{th}$, 2004}\\
  \text{J.-P. Bernier}& \text{(Athlon 2000)}& 65.1\times 10^{15} \text{ to } 65.5\times 10^{15}& \text{December $19^{th}$, 2003}\\
  \text{J.-P. Bernier}& \text{(Pentium 600MHz)}& 64.8\times 10^{15} \text{ to } 65.1\times 10^{15}& \text{December $14^{th}$, 2003}\\
  \text{J.-P. Bernier}& \text{(Athlon 2000)}& 64.4\times 10^{15} \text{ to } 64.8\times 10^{15}& \text{December $4^{th}$, 2003}\\
  \text{B. Le Tual}& \text{(Celeron 2.4GHz)}& 64.1\times 10^{15} \text{ to } 64.4\times 10^{15}& \text{December $9^{th}$, 2003}\\
  \text{L. Marmet}& \text{(neurone0 W95 .233GHz)}& 63.8\times 10^{15} \text{ to } 64.1\times 10^{15}& \text{April $11^{th}$, 2005}\\
  \text{J.-P. Bernier}& \text{(Athlon 2000)}& 63.4\times 10^{15} \text{ to } 63.8\times 10^{15}& \text{November $13^{th}$, 2003}\\
  \text{B. Le Tual}& \text{(Celeron .9GHz)}& 63.0\times 10^{15} \text{ to } 63.4\times 10^{15}& \text{December $5^{th}$, 2003}\\
  \text{J.-P. Bernier}& \text{(Pentium 600MHz)}& 62.7\times 10^{15} \text{ to } 63.0\times 10^{15}& \text{November $24^{th}$, 2003}\\
  \text{D. Bernier}& \text{(Gecko)}& 62.4\times 10^{15} \text{ to } 62.7\times 10^{15}& \text{December $27^{th}$, 2003}\\
  \text{J.-P. Bernier}& \text{(Athlon 2000)}& 62.0\times 10^{15} \text{ to } 62.4\times 10^{15}& \text{November $1^{st}$, 2003}\\
  \text{N. Robertson}& \text{(Laika.Droica)}& 61.6\times 10^{15} \text{ to } 62.0\times 10^{15}& \text{March $9^{th}$, 2005}\\
  \text{L. Marmet}& \text{(neurone5 Linux .8GHz)}& 61.2\times 10^{15} \text{ to } 61.6\times 10^{15}& \text{March $27^{th}$, 2005}\\
  \text{J.-P. Bernier}& \text{(Athlon 2000)}& 60.9\times 10^{15} \text{ to } 61.2\times 10^{15}& \text{October $15^{th}$, 2003}\\
  \text{L. Marmet}& \text{(Riyadh)}& 60.6\times 10^{15} \text{ to } 60.9\times 10^{15}& \text{November $9^{th}$, 2003}\\
  \text{J.-P. Bernier}& \text{(Pentium 600MHz)}& 60.3\times 10^{15} \text{ to } 60.6\times 10^{15}& \text{October $31^{st}$, 2003}\\
  \text{B. Le Tual}& \text{(Celeron 2.4GHz)}& 60.0\times 10^{15} \text{ to } 60.3\times 10^{15}& \text{November $12^{th}$, 2003}\\
  \text{B. Le Tual}& \text{(Celeron .9GHz)}& 59.6\times 10^{15} \text{ to } 60.0\times 10^{15}& \text{November $1^{st}$, 2003}\\
  \text{N. Robertson}& \text{(Laika.Droica)}& 59.2\times 10^{15} \text{ to } 59.6\times 10^{15}& \text{October $1^{st}$, 2003}\\
  \text{B. Le Tual}& \text{(Celeron .9GHz)}& 58.9\times 10^{15} \text{ to } 59.2\times 10^{15}& \text{September $28^{th}$, 2003}\\
  \text{B. Le Tual}& \text{(Celeron 2.4GHz)}& 58.6\times 10^{15} \text{ to } 58.9\times 10^{15}& \text{October $5^{th}$, 2003}\\
  \text{J.-P. Bernier}& \text{(Pentium 600MHz)}& 58.3\times 10^{15} \text{ to } 58.6\times 10^{15}& \text{October $7^{th}$, 2003}\\
  \text{L. Marmet}& \text{(Riyadh)}& 58.0\times 10^{15} \text{ to } 58.3\times 10^{15}& \text{October $7^{th}$, 2003}\\
  \text{D. Bernier}& \text{(Gecko)}& 57.7\times 10^{15} \text{ to } 58.0\times 10^{15}& \text{October $29^{th}$, 2003}\\
  \text{N. Robertson}& \text{(Maya)}& 57.4\times 10^{15} \text{ to } 57.7\times 10^{15}& \text{September $19^{th}$, 2003}\\
  \text{N. Robertson}& \text{(Laika.Droica)}& 57.1\times 10^{15} \text{ to } 57.4\times 10^{15}& \text{September $10^{th}$, 2003}\\
  \text{J.-P. Bernier}& \text{(Pentium 600MHz)}& 56.8\times 10^{15} \text{ to } 57.1\times 10^{15}& \text{September $10^{th}$, 2003}\\
  \text{N. Robertson}& \text{(Maya)}& 56.5\times 10^{15} \text{ to } 56.8\times 10^{15}& \text{August $4^{th}$, 2003}\\
  \text{N. Robertson}& \text{(Laika.Droica)}& 56.2\times 10^{15} \text{ to } 56.5\times 10^{15}& \text{July $29^{th}$, 2003}\\
  \text{L. Marmet}& \text{(Riyadh)}& 55.9\times 10^{15} \text{ to } 56.2\times 10^{15}& \text{August $13^{th}$, 2003}\\
  \text{J.-P. Bernier}& \text{(Pentium 600MHz)}& 55.6\times 10^{15} \text{ to } 55.9\times 10^{15}& \text{August $4^{th}$, 2003}\\
  \end{array}
\end{displaymath}

\begin{displaymath}
  \begin{array}{llrr}
  \text{User Name}& \text{(Computer Name)}& Range& Date Completed\\
\hline
  \text{J.-P. Bernier}& \text{(Pentium 600MHz)}& 55.3\times 10^{15} \text{ to } 55.6\times 10^{15}& \text{July $10^{th}$, 2003}\\
  \text{J.-P. Bernier}& \text{(Pentium 600MHz)}& 55.0\times 10^{15} \text{ to } 55.3\times 10^{15}& \text{June $15^{th}$, 2003}\\
  \text{L. Marmet}& \text{(Riyadh)}& 54.7\times 10^{15} \text{ to } 55.0\times 10^{15}& \text{June $28^{th}$, 2003}\\
  \text{N. Robertson}& \text{(Maya)}& 54.4\times 10^{15} \text{ to } 54.7\times 10^{15}& \text{June $10^{th}$, 2003}\\
  \text{N. Robertson}& \text{(Laika.Droica)}& 54.1\times 10^{15} \text{ to } 54.4\times 10^{15}& \text{July $9^{th}$, 2003}\\
  \text{J.-P. Bernier}& \text{(Pentium 600MHz)}& 53.8\times 10^{15} \text{ to } 54.1\times 10^{15}& \text{May $14^{th}$, 2003}\\
  \text{D. Bernier}& \text{(Gecko)}& 53.3\times 10^{15} \text{ to } 53.8\times 10^{15}& \text{August $23^{rd}$, 2003}\\
  \text{N. Robertson}& \text{(Maya)}& 53.0\times 10^{15} \text{ to } 53.3\times 10^{15}& \text{May $12^{th}$, 2003}\\
  \text{N. Robertson}& \text{(Laika.Droica)}& 52.7\times 10^{15} \text{ to } 53.0\times 10^{15}& \text{May $7^{th}$, 2003}\\
  \text{L. Marmet}& \text{(Riyadh)}& 52.4\times 10^{15} \text{ to } 52.7\times 10^{15}& \text{May $16^{th}$, 2003}\\
  \text{J.-P. Bernier}& \text{(Pentium 600MHz)}& 52.1\times 10^{15} \text{ to } 52.4\times 10^{15}& \text{April $23^{rd}$, 2003}\\
  \text{N. Robertson}& \text{(Maya)}& 51.8\times 10^{15} \text{ to } 52.1\times 10^{15}& \text{March $31^{st}$, 2003}\\
  \text{L. Marmet}& \text{(Riyadh)}& 51.5\times 10^{15} \text{ to } 51.8\times 10^{15}& \text{April $6^{th}$, 2003}\\
  \text{N. Robertson}& \text{(Laika.Droica)}& 51.2\times 10^{15} \text{ to } 51.5\times 10^{15}& \text{March $22^{nd}$, 2003}\\
  \text{D. Bernier}& \text{(Gecko)}& 51.1\times 10^{15} \text{ to } 51.2\times 10^{15}& \text{April $20^{th}$, 2003}\\
  \text{D. Bernier}& \text{(Gecko)}& 50.7\times 10^{15} \text{ to } 51.1\times 10^{15}& \text{April $2^{nd}$, 2003}\\
  \text{J.-P. Bernier}& \text{(Pentium 600MHz)}& 50.4\times 10^{15} \text{ to } 50.7\times 10^{15}& \text{March $24^{th}$, 2003}\\
  \text{D. Bernier}& \text{(Gecko)}& 50.1\times 10^{15} \text{ to } 50.4\times 10^{15}& \text{February $26^{th}$, 2003}\\
  \text{L. Marmet}& \text{(Riyadh)}& 49.8\times 10^{15} \text{ to } 50.1\times 10^{15}& \text{March $9^{th}$, 2003}\\
  \text{J.-P. Bernier}& \text{(Pentium 600MHz)}& 49.5\times 10^{15} \text{ to } 49.8\times 10^{15}& \text{February $25^{th}$, 2003}\\
  \text{L. Marmet}& \text{(Riyadh)}& 49.2\times 10^{15} \text{ to } 49.5\times 10^{15}& \text{January $31^{st}$, 2003}\\
  \text{J.-P. Bernier}& \text{(Pentium 600MHz)}& 49.0\times 10^{15} \text{ to } 49.2\times 10^{15}& \text{January $25^{th}$, 2003}\\
  \text{J.-P. Bernier}& \text{(Pentium 600MHz)}& 48.8\times 10^{15} \text{ to } 49.0\times 10^{15}& \text{December $26^{th}$, 2002}\\
  \text{L. Marmet}& \text{(Riyadh)}& 48.5\times 10^{15} \text{ to } 48.8\times 10^{15}& \text{December $22^{nd}$, 2002}\\
  \text{N. Robertson}& \text{(Maya)}& 48.2\times 10^{15} \text{ to } 48.5\times 10^{15}& \text{November $28^{th}$, 2002}\\
  \text{N. Robertson}& \text{(Laika.Droica)}& 47.9\times 10^{15} \text{ to } 48.2\times 10^{15}& \text{January $11^{th}$, 2003}\\
  \text{N. Robertson}& \text{(Arthurus)}& 47.8\times 10^{15} \text{ to } 47.9\times 10^{15}& \text{February $13^{th}$, 2003}\\
  \text{J.-P. Bernier}& \text{(Pentium 600MHz)}& 47.5\times 10^{15} \text{ to } 47.8\times 10^{15}& \text{December $6^{th}$, 2002}\\
  \text{D. Bernier}& \text{(Pentium 500MHz)}& 47.2\times 10^{15} \text{ to } 47.5\times 10^{15}& \text{November $6^{th}$, 2002}\\
  \text{J.-P. Bernier}& \text{(Pentium 600MHz)}& 46.9\times 10^{15} \text{ to } 47.2\times 10^{15}& \text{November $7^{th}$, 2002}\\
  \text{D. Bernier}& \text{(Pentium 500MHz)}& 46.6\times 10^{15} \text{ to } 46.9\times 10^{15}& \text{October $14^{th}$, 2002}\\
  \text{J.-P. Bernier}& \text{(Pentium 600MHz)}& 46.3\times 10^{15} \text{ to } 46.6\times 10^{15}& \text{October $13^{th}$, 2002}\\
  \text{L. Marmet}& \text{(Riyadh)}& 46.0\times 10^{15} \text{ to } 46.3\times 10^{15}& \text{October $15^{th}$, 2002}\\
  \text{N. Robertson}& \text{(Maya)}& 45.7\times 10^{15} \text{ to } 46.0\times 10^{15}& \text{November $3^{rd}$, 2002}\\
  \text{N. Robertson}& \text{(Laika.Droica)}& 45.4\times 10^{15} \text{ to } 45.7\times 10^{15}& \text{October $30^{th}$, 2002}\\
  \text{J.-P. Bernier}& \text{(Pentium 600MHz)}& 45.1\times 10^{15} \text{ to } 45.4\times 10^{15}& \text{September $12^{th}$, 2002}\\
  \text{D. Bernier}& \text{(Pentium 500MHz)}& 44.8\times 10^{15} \text{ to } 45.1\times 10^{15}& \text{September $15^{th}$, 2002}\\
  \text{J.-P. Bernier}& \text{(Pentium 600MHz)}& 44.5\times 10^{15} \text{ to } 44.8\times 10^{15}& \text{August $13^{th}$, 2002}\\
  \text{D. Bernier}& \text{(Pentium 500MHz)}& 44.2\times 10^{15} \text{ to } 44.5\times 10^{15}& \text{August $12^{th}$, 2002}\\
  \text{L. Marmet}& \text{(Riyadh)}& 43.9\times 10^{15} \text{ to } 44.2\times 10^{15}& \text{August $7^{th}$, 2002}\\
  \end{array}
\end{displaymath}

\begin{displaymath}
  \begin{array}{llrr}
  \text{User Name}& \text{(Computer Name)}& Range& Date Completed\\
\hline
  \text{N. Robertson}& \text{(Maya)}& 43.6\times 10^{15} \text{ to } 43.9\times 10^{15}& \text{July $31^{st}$, 2002}\\
  \text{N. Robertson}& \text{(Laika.Droica)}& 43.3\times 10^{15} \text{ to } 43.6\times 10^{15}& \text{July $27^{th}$, 2002}\\
  \text{J.-P. Bernier}& \text{(Pentium 600MHz)}& 43.0\times 10^{15} \text{ to } 43.3\times 10^{15}& \text{July $23^{rd}$, 2002}\\
  \text{D. Bernier}& \text{(Pentium 500MHz)}& 42.7\times 10^{15} \text{ to } 43.0\times 10^{15}& \text{July $11^{th}$, 2002}\\
  \text{L. Marmet}& \text{(Riyadh)}& 42.4\times 10^{15} \text{ to } 42.7\times 10^{15}& \text{July $11^{th}$, 2002}\\
  \text{J.-P. Bernier}& \text{(Pentium 600MHz)}& 42.1\times 10^{15} \text{ to } 42.4\times 10^{15}& \text{June $24^{th}$, 2002}\\
  \text{N. Robertson}& \text{(Laika.Droica)}& 41.8\times 10^{15} \text{ to } 42.1\times 10^{15}& \text{June $17^{th}$, 2002}\\
  \text{L. Marmet}& \text{(Riyadh)}& 41.5\times 10^{15} \text{ to } 41.8\times 10^{15}& \text{June $14^{th}$, 2002}\\
  \text{N. Robertson}& \text{(Maya)}& 41.2\times 10^{15} \text{ to } 41.5\times 10^{15}& \text{June $9^{th}$, 2002}\\
  \text{N. Robertson}& \text{(Laika.Droica)}& 40.9\times 10^{15} \text{ to } 41.2\times 10^{15}& \text{May $28^{th}$, 2002}\\
  \text{J.-P. Bernier}& \text{(Pentium 600MHz)}& 40.6\times 10^{15} \text{ to } 40.9\times 10^{15}& \text{June $5^{th}$, 2002}\\
  \text{D. Bernier}& \text{(Pentium 500MHz)}& 40.3\times 10^{15} \text{ to } 40.6\times 10^{15}& \text{May $26^{th}$, 2002}\\
  \text{N. Robertson}& \text{(Maya)}& 40.0\times 10^{15} \text{ to } 40.3\times 10^{15}& \text{March $6^{th}$, 2003}\\
  \text{J.-P. Bernier}& \text{(Pentium 600MHz)}& 39.7\times 10^{15} \text{ to } 40.0\times 10^{15}& \text{May $2^{nd}$, 2002}\\
  \text{D. Bernier}& \text{(Pentium 500MHz)}& 39.4\times 10^{15} \text{ to } 39.7\times 10^{15}& \text{May $1^{st}$, 2002}\\
  \text{N. Robertson}& \text{(Maya)}& 39.1\times 10^{15} \text{ to } 39.4\times 10^{15}& \text{April $23^{rd}$, 2002}\\
  \text{N. Robertson}& \text{(Laika.Droica)}& 38.8\times 10^{15} \text{ to } 39.1\times 10^{15}& \text{April $18^{th}$, 2002}\\
  \text{L. Marmet}& \text{(Riyadh)}& 38.5\times 10^{15} \text{ to } 38.8\times 10^{15}& \text{April $26^{th}$, 2002}\\
  \text{N. Robertson}& \text{(Maya)}& 38.2\times 10^{15} \text{ to } 38.5\times 10^{15}& \text{March $25^{th}$, 2002}\\
  \text{N. Robertson}& \text{(Laika.Droica)}& 37.9\times 10^{15} \text{ to } 38.2\times 10^{15}& \text{March $20^{th}$, 2002}\\
  \text{L. Marmet}& \text{(Riyadh)}& 37.6\times 10^{15} \text{ to } 37.9\times 10^{15}& \text{March $27^{th}$, 2002}\\
  \text{J.-P. Bernier}& \text{(Pentium 600MHz)}& 37.3\times 10^{15} \text{ to } 37.6\times 10^{15}& \text{March $6^{th}$, 2002}\\
  \text{D. Bernier}& \text{(Pentium 500MHz)}& 37.0\times 10^{15} \text{ to } 37.3\times 10^{15}& \text{March $12^{th}$, 2002}\\
  \text{N. Robertson}& \text{(Maya)}& 36.7\times 10^{15} \text{ to } 37.0\times 10^{15}& \text{March $4^{th}$, 2002}\\
  \text{N. Robertson}& \text{(Laika.Droica)}& 36.4\times 10^{15} \text{ to } 36.7\times 10^{15}& \text{February $28^{th}$, 2002}\\ 
  \text{D. Bernier}& \text{(Pentium 500MHz)}& 36.2\times 10^{15} \text{ to } 36.4\times 10^{15}& \text{June $17^{th}$, 2002}\\
  \text{N. Robertson}& \text{(Arthurus)}& 36.0\times 10^{15} \text{ to } 36.2\times 10^{15}& \text{October $2^{nd}$, 2002}\\
  \text{L. Marmet}& \text{(Riyadh)}& 35.7\times 10^{15} \text{ to } 36.0\times 10^{15}& \text{February $18^{th}$, 2002}\\
  \text{C.R. Ward}& \text{(Cosmos)}& 35.4\times 10^{15} \text{ to } 35.7\times 10^{15}& \text{March $10^{th}$, 2002}\\
  \text{J.-P. Bernier}& \text{(Pentium 600MHz)}& 35.1\times 10^{15} \text{ to } 35.4\times 10^{15}& \text{February $5^{th}$, 2002}\\
  \text{D. Bernier}& \text{(Pentium 500MHz)}& 34.8\times 10^{15} \text{ to } 35.1\times 10^{15}& \text{February $2^{nd}$, 2002}\\
  \text{L. Marmet}& \text{(Riyadh)}& 34.5\times 10^{15} \text{ to } 34.8\times 10^{15}& \text{January $6^{th}$, 2002}\\
  \text{N. Robertson}& \text{(Maya)}& 34.2\times 10^{15} \text{ to } 34.5\times 10^{15}& \text{December $30^{th}$, 2001}\\
  \text{N. Robertson}& \text{(Laika.Droica)}& 33.9\times 10^{15} \text{ to } 34.2\times 10^{15}& \text{December $25^{th}$, 2001}\\
  \text{J.-P. Bernier}& \text{(Pentium 600MHz)}& 33.6\times 10^{15} \text{ to } 33.9\times 10^{15}& \text{January $10^{th}$, 2002}\\
  \text{C.R. Ward}& \text{(Cosmos)}& 33.3\times 10^{15} \text{ to } 33.6\times 10^{15}& \text{January $16^{th}$, 2002}\\
  \text{L. Marmet}& \text{(Riyadh)}& 33.0\times 10^{15} \text{ to } 33.3\times 10^{15}& \text{December $8^{th}$, 2001}\\
  \text{N. Robertson}& \text{(Maya)}& 32.7\times 10^{15} \text{ to } 33.0\times 10^{15}& \text{December $4^{th}$, 2001}\\
  \text{N. Robertson}& \text{(Laika.Droica)}& 32.4\times 10^{15} \text{ to } 32.7\times 10^{15}& \text{November $27^{th}$, 2001}\\
  \text{D. Bernier}& \text{(Pentium 500MHz)}& 32.1\times 10^{15} \text{ to } 32.4\times 10^{15}& \text{December $26^{th}$, 2001}\\
  \end{array}
\end{displaymath}

\begin{displaymath}
  \begin{array}{llrr}
  \text{User Name}& \text{(Computer Name)}& Range& Date Completed\\
\hline
  \text{J.-P. Bernier}& \text{(Pentium 600MHz)}& 31.8\times 10^{15} \text{ to } 32.1\times 10^{15}& \text{December $6^{th}$, 2001}\\
  \text{J.-P. Bernier}& \text{(Pentium 600MHz)}& 31.5\times 10^{15} \text{ to } 31.8\times 10^{15}& \text{November $1^{st}$, 2001}\\
  \text{C.R. Ward}& \text{(Cosmos)}& 31.2\times 10^{15} \text{ to } 31.5\times 10^{15}& \text{November $25^{th}$, 2001}\\
  \text{L. Marmet}& \text{(Riyadh)}& 30.9\times 10^{15} \text{ to } 31.2\times 10^{15}& \text{November $2^{nd}$, 2001}\\
  \text{N. Robertson}& \text{(Maya)}& 30.6\times 10^{15} \text{ to } 30.9\times 10^{15}& \text{October $26^{th}$, 2001}\\
  \text{N. Robertson}& \text{(Laika.Droica)}& 30.3\times 10^{15} \text{ to } 30.6\times 10^{15}& \text{October $20^{th}$, 2001}\\
  \text{J.-P. Bernier}& \text{(Pentium 600MHz)}& 30.0\times 10^{15} \text{ to } 30.3\times 10^{15}& \text{September $29^{th}$, 2001}\\
  \text{N. Robertson}& \text{(Maya)}& 29.7\times 10^{15} \text{ to } 30.0\times 10^{15}& \text{September $18^{th}$, 2001}\\
  \text{N. Robertson}& \text{(Laika.Droica)}& 29.4\times 10^{15} \text{ to } 29.7\times 10^{15}& \text{September $16^{th}$, 2001}\\
  \text{L. Marmet}& \text{(Riyadh)}& 29.1\times 10^{15} \text{ to } 29.4\times 10^{15}& \text{September $19^{th}$, 2001}\\
  \text{N. Robertson}& \text{(Rosette.Droica)}& 28.9\times 10^{15} \text{ to } 29.1\times 10^{15}& \text{December $12^{th}$, 2001}\\
  \text{N. Robertson}& \text{(Laika.Droica)}& 28.6\times 10^{15} \text{ to } 28.9\times 10^{15}& \text{August $28^{th}$, 2001}\\
  \text{J.-P. Bernier}& \text{(Pentium 600MHz)}& 28.3\times 10^{15} \text{ to } 28.6\times 10^{15}& \text{September $2^{nd}$, 2001}\\
  \text{N. Robertson}& \text{(Maya)}& 28.0\times 10^{15} \text{ to } 28.3\times 10^{15}& \text{August $28^{th}$, 2001}\\
  \text{C.R. Ward}& \text{(Cosmos)}& 27.7\times 10^{15} \text{ to } 28.0\times 10^{15}& \text{September $22^{nd}$, 2001}\\
  \text{L. Marmet}& \text{(Riyadh)}& 27.4\times 10^{15} \text{ to } 27.7\times 10^{15}& \text{August $22^{nd}$, 2001}\\
  \text{N. Robertson}& \text{(Maya)}& 27.1\times 10^{15} \text{ to } 27.4\times 10^{15}& \text{August $7^{th}$, 2001}\\
  \text{J.-P. Bernier}& \text{(Pentium 600MHz)}& 26.8\times 10^{15} \text{ to } 27.1\times 10^{15}& \text{August $8^{th}$, 2001}\\
  \text{D. Bernier}& \text{(Pentium 500MHz)}& 26.5\times 10^{15} \text{ to } 26.8\times 10^{15}& \text{April $6^{th}$, 2002}\\
  \text{C.R. Ward}& \text{(Cosmos)}& 26.2\times 10^{15} \text{ to } 26.5\times 10^{15}& \text{April $29^{th}$, 2002}\\
  \text{D. Bernier}& \text{(Pentium 500MHz)}& 25.9\times 10^{15} \text{ to } 26.2\times 10^{15}& \text{November $2^{nd}$, 2001}\\
  \text{N. Robertson}& \text{(Maya)}& 25.6\times 10^{15} \text{ to } 25.9\times 10^{15}& \text{July $16^{th}$, 2001}\\
  \text{J.-P. Bernier}& \text{(Pentium 600MHz)}& 25.3\times 10^{15} \text{ to } 25.6\times 10^{15}& \text{July $15^{th}$, 2001}\\
  \text{L. Marmet}& \text{(Riyadh)}& 25.0\times 10^{15} \text{ to } 25.3\times 10^{15}& \text{July $23^{rd}$, 2001}\\
  \text{C.R. Ward}& \text{(Cosmos)}& 24.7\times 10^{15} \text{ to } 25.0\times 10^{15}& \text{August $2^{nd}$, 2001}\\
  \text{J.-P. Bernier}& \text{(Pentium 600MHz)}& 24.4\times 10^{15} \text{ to } 24.7\times 10^{15}& \text{April $8^{th}$, 2002}\\
  \text{Z. McGregor-Dorsey}& \text{(Abzug)}& 24.1\times 10^{15} \text{ to } 24.4\times 10^{15}& \text{July $7^{th}$, 2001}\\
  \text{N. Robertson}& \text{(Arthurus)}& 23.9\times 10^{15} \text{ to } 24.1\times 10^{15}& \text{November $30^{th}$, 2001}\\
  \text{L. Marmet}& \text{(Riyadh)}& 23.6\times 10^{15} \text{ to } 23.9\times 10^{15}& \text{June $14^{th}$, 2001}\\
  \text{N. Robertson}& \text{(Maya)}& 23.3\times 10^{15} \text{ to } 23.6\times 10^{15}& \text{June $18^{th}$, 2001}\\
  \text{Z. McGregor-Dorsey}& \text{(Hayduke)}& 23.0\times 10^{15} \text{ to } 23.3\times 10^{15}& \text{July $7^{th}$, 2001}\\
  \text{J.-P. Bernier}& \text{(Pentium 600MHz)}& 22.7\times 10^{15} \text{ to } 23.0\times 10^{15}& \text{June $16^{th}$, 2001}\\
  \text{J.-P. Bernier}& \text{(Pentium 600MHz)}& 22.4\times 10^{15} \text{ to } 22.7\times 10^{15}& \text{May $26^{th}$, 2001}\\
  \text{D. Bernier}& \text{(Pentium 500MHz)}& 22.1\times 10^{15} \text{ to } 22.4\times 10^{15}& \text{June $19^{th}$, 2001}\\
  \text{Z. McGregor-Dorsey}& \text{(Castalia)}& 21.8\times 10^{15} \text{ to } 22.1\times 10^{15}& \text{June $5^{th}$, 2001}\\
  \text{C.R. Ward}& \text{(Cosmos)}& 21.5\times 10^{15} \text{ to } 21.8\times 10^{15}& \text{June $7^{th}$, 2001}\\
  \text{J.-P. Bernier}& \text{(Pentium 600MHz)}& 21.2\times 10^{15} \text{ to } 21.5\times 10^{15}& \text{May $8^{th}$, 2001}\\
  \text{L. Marmet}& \text{(Riyadh)}& 20.9\times 10^{15} \text{ to } 21.2\times 10^{15}& \text{May $16^{th}$, 2001}\\
  \text{D. Bernier}& \text{(Pentium 500MHz)}& 20.6\times 10^{15} \text{ to } 20.9\times 10^{15}& \text{May $3^{rd}$, 2001}\\
  \text{L. Marmet}& \text{(Fontaine)}& 20.3\times 10^{15} \text{ to } 20.6\times 10^{15}& \text{May $8^{th}$, 2001}\\
  \end{array}
\end{displaymath}

\begin{displaymath}
  \begin{array}{llrr}
  \text{User Name}& \text{(Computer Name)}& Range& Date Completed\\
\hline
  \text{J.-P. Bernier}& \text{(Pentium 600MHz)}& 20.0\times 10^{15} \text{ to } 20.3\times 10^{15}& \text{April $20^{th}$, 2001}\\
  \text{C.R. Ward}& \text{(Cosmos)}& 19.8\times 10^{15} \text{ to } 20.0\times 10^{15}& \text{April $20^{th}$, 2001}\\
  \text{L. Marmet}& \text{(Riyadh)}& 19.6\times 10^{15} \text{ to } 19.8\times 10^{15}& \text{April $18^{th}$, 2001}\\
  \text{Z. McGregor-Dorsey}& \text{(Abzug)}& 19.4\times 10^{15} \text{ to } 19.6\times 10^{15}& \text{May $27^{th}$, 2001}\\
  \text{Z. McGregor-Dorsey}& \text{(Hayduke)}& 19.2\times 10^{15} \text{ to } 19.4\times 10^{15}& \text{May $10^{th}$, 2001}\\
  \text{D. Bernier}& \text{(Pentium 500MHz)}& 19.0\times 10^{15} \text{ to } 19.2\times 10^{15}& \text{April $11^{th}$, 2001}\\
  \text{L. Marmet}& \text{(Fontaine)}& 18.8\times 10^{15} \text{ to } 19.0\times 10^{15}& \text{April $7^{th}$, 2001}\\
  \text{D. Bernier}& \text{(Pentium 500MHz)}& 18.6\times 10^{15} \text{ to } 18.8\times 10^{15}& \text{April $4^{th}$, 2001}\\
  \text{J.-P. Bernier}& \text{(Pentium 600MHz)}& 18.4\times 10^{15} \text{ to } 18.6\times 10^{15}& \text{March $27^{th}$, 2001}\\
  \text{L. Marmet}& \text{(Riyadh)}& 18.2\times 10^{15} \text{ to } 18.4\times 10^{15}& \text{April $1^{st}$, 2001}\\
  \text{D. Bernier}& \text{(Pentium 500MHz)}& 18.0\times 10^{15} \text{ to } 18.2\times 10^{15}& \text{March $18^{th}$, 2001}\\
  \text{L. Marmet}& \text{(Fontaine)}& 17.8\times 10^{15} \text{ to } 18.0\times 10^{15}& \text{March $21^{st}$, 2001}\\
  \text{L. Marmet}& \text{(Riyadh)}& 17.6\times 10^{15} \text{ to } 17.8\times 10^{15}& \text{March $15^{th}$, 2001}\\
  \text{J.-P. Bernier}& \text{(Pentium 600MHz)}& 17.4\times 10^{15} \text{ to } 17.6\times 10^{15}& \text{March $14^{th}$, 2001}\\
  \text{D. Bernier}& \text{(Pentium 500MHz)}& 17.2\times 10^{15} \text{ to } 17.4\times 10^{15}& \text{March $5^{th}$, 2001}\\
  \text{J.-P. Bernier}& \text{(Pentium 600MHz)}& 17.0\times 10^{15} \text{ to } 17.2\times 10^{15}& \text{March $2^{nd}$, 2001}\\
  \text{L. Marmet}& \text{(Fontaine)}& 16.8\times 10^{15} \text{ to } 17.0\times 10^{15}& \text{March $3^{rd}$, 2001}\\
  \text{L. Marmet}& \text{(Riyadh)}& 16.6\times 10^{15} \text{ to } 16.8\times 10^{15}& \text{February $26^{th}$, 2001}\\
  \text{D. Bernier}& \text{(Pentium 500MHz)}& 16.4\times 10^{15} \text{ to } 16.6\times 10^{15}& \text{February $19^{th}$, 2001}\\
  \text{J.-P. Bernier}& \text{(Pentium 600MHz)}& 16.2\times 10^{15} \text{ to } 16.4\times 10^{15}& \text{February $17^{th}$, 2001}\\
  \text{N. Robertson}& \text{(Arthurus)}& 16.0\times 10^{15} \text{ to } 16.2\times 10^{15}& \text{May $27^{th}$, 2001}\\
  \text{L. Marmet}& \text{(Fontaine)}& 15.8\times 10^{15} \text{ to } 16.0\times 10^{15}& \text{February $11^{th}$, 2001}\\
  \text{Z. McGregor-Dorsey}& \text{(Castalia)}& 15.6\times 10^{15} \text{ to } 15.8\times 10^{15}& \text{April $25^{th}$, 2001}\\
  \text{Z. McGregor-Dorsey}& \text{(Abzug)}& 15.4\times 10^{15} \text{ to } 15.6\times 10^{15}& \text{April $17^{th}$, 2001}\\
  \text{L. Marmet}& \text{(Riyadh)}& 15.2\times 10^{15} \text{ to } 15.4\times 10^{15}& \text{February $8^{th}$, 2001}\\
  \text{Z. McGregor-Dorsey}& \text{(Hayduke)}& 15.0\times 10^{15} \text{ to } 15.2\times 10^{15}& \text{April $5^{th}$, 2001}\\
  \text{J.-P. Bernier}& \text{(Pentium 600MHz)}& 14.8\times 10^{15} \text{ to } 15.0\times 10^{15}& \text{February $4^{th}$, 2001}\\
  \text{L. Marmet}& \text{(Riyadh)}& 14.6\times 10^{15} \text{ to } 14.8\times 10^{15}& \text{January $22^{nd}$, 2001}\\
  \text{Z. McGregor-Dorsey}& \text{(Castalia)}& 14.4\times 10^{15} \text{ to } 14.6\times 10^{15}& \text{April $3^{rd}$, 2001}\\
  \text{Z. McGregor-Dorsey}& \text{(Abzug)}& 14.2\times 10^{15} \text{ to } 14.4\times 10^{15}& \text{March $19^{th}$, 2001}\\
  \text{Z. McGregor-Dorsey}& \text{(Hayduke)}& 14.0\times 10^{15} \text{ to } 14.2\times 10^{15}& \text{March $19^{th}$, 2001}\\
  \text{L. Marmet}& \text{(Strontium)}& 13.8\times 10^{15} \text{ to } 14.0\times 10^{15}& \text{December $17^{th}$, 2000}\\
  \text{L. Marmet}& \text{(Riyadh)}& 13.6\times 10^{15} \text{ to } 13.8\times 10^{15}& \text{January $5^{th}$, 2001}\\
  \text{L. Marmet}& \text{(Strontium)}& 13.4\times 10^{15} \text{ to } 13.6\times 10^{15}& \text{December $11^{th}$, 2000}\\
  \text{L. Marmet}& \text{(Fontaine)}& 13.2\times 10^{15} \text{ to } 13.4\times 10^{15}& \text{December $13^{th}$, 2000}\\
  \text{L. Marmet}& \text{(Strontium)}& 13.0\times 10^{15} \text{ to } 13.2\times 10^{15}& \text{December $4^{th}$, 2000}\\
  \text{L. Marmet}& \text{(Riyadh)}& 12.8\times 10^{15} \text{ to } 13.0\times 10^{15}& \text{December $9^{th}$, 2000}\\
  \text{L. Marmet}& \text{(Strontium)}& 12.6\times 10^{15} \text{ to } 12.8\times 10^{15}& \text{November $26^{th}$, 2000}\\
  \text{L. Marmet}& \text{(Fontaine)}& 12.4\times 10^{15} \text{ to } 12.6\times 10^{15}& \text{November $24^{th}$, 2000}\\
  \text{Z. McGregor-Dorsey}& \text{(Abzug)}& 12.2\times 10^{15} \text{ to } 12.4\times 10^{15}& \text{January $21^{st}$, 2001}\\
  \end{array}
\end{displaymath}

\begin{displaymath}
  \begin{array}{llrr}
  \text{User Name}& \text{(Computer Name)}& Range& Date Completed\\
\hline
  \text{L. Marmet}& \text{(Riyadh)}& 12.0\times 10^{15} \text{ to } 12.2\times 10^{15}& \text{November $22^{nd}$, 2000}\\
  \text{L. Marmet}& \text{(Strontium)}& 11.8\times 10^{15} \text{ to } 12.0\times 10^{15}& \text{November $17^{th}$, 2000}\\
  \text{Z. McGregor-Dorsey}& \text{(Castalia)}& 11.6\times 10^{15} \text{ to } 11.8\times 10^{15}& \text{January $22^{nd}$, 2001}\\
  \text{L. Marmet}& \text{(Riyadh)}& 11.4\times 10^{15} \text{ to } 11.6\times 10^{15}& \text{November $7^{th}$, 2000}\\
  \text{Z. McGregor-Dorsey}& \text{(Hayduke)}& 11.2\times 10^{15} \text{ to } 11.4\times 10^{15}& \text{January $21^{st}$, 2001}\\
  \text{L. Marmet}& \text{(Riyadh)}& 11.0\times 10^{15} \text{ to } 11.2\times 10^{15}& \text{November $1^{st}$, 2000}\\
  \text{Z. McGregor-Dorsey}& \text{(Hayduke)}& 10.8\times 10^{15} \text{ to } 11.0\times 10^{15}& \text{December $10^{th}$, 2000}\\
  \text{Z. McGregor-Dorsey}& \text{(Castalia)}& 10.6\times 10^{15} \text{ to } 10.8\times 10^{15}& \text{December $10^{th}$, 2000}\\
  \text{L. Marmet}& \text{(Riyadh)}& 10.4\times 10^{15} \text{ to } 10.6\times 10^{15}& \text{October $13^{th}$, 2000}\\
  \text{N. Robertson}& \text{(Arthurus)}& 10.2\times 10^{15} \text{ to } 10.4\times 10^{15}& \text{January $27^{th}$, 2001}\\
  \text{N. Marmet}& \text{(Computer)}& 10.0\times 10^{15} \text{ to } 10.2\times 10^{15}& \text{October $21^{st}$, 2000}\\
  \text{L. Marmet}& \text{(Riyadh)}& 9.8\times 10^{15} \text{ to } 10.0\times 10^{15}& \text{September $27^{th}$, 2000}\\
  \text{Z. McGregor-Dorsey}& \text{(Castalia)}& 9.6\times 10^{15} \text{ to } 9.8\times 10^{15}& \text{November $16^{th}$, 2000}\\
  \text{Z. McGregor-Dorsey}& \text{(Hayduke)}& 9.4\times 10^{15} \text{ to } 9.6\times 10^{15}& \text{January $26^{th}$, 2001}\\
  \text{L. Marmet}& \text{(Fontaine)}& 9.2\times 10^{15} \text{ to } 9.4\times 10^{15}& \text{January $26^{th}$, 2001}\\
  \text{Z. McGregor-Dorsey}& \text{(Castalia)}& 9.0\times 10^{15} \text{ to } 9.2\times 10^{15}& \text{October $24^{th}$, 2000}\\
  \text{Z. McGregor-Dorsey}& \text{(Castalia)}& 8.8\times 10^{15} \text{ to } 9.0\times 10^{15}& \text{September $27^{th}$, 2000}\\
  \text{Z. McGregor-Dorsey}& \text{(Hayduke)}& 8.6\times 10^{15} \text{ to } 8.8\times 10^{15}& \text{September $11^{th}$, 2000}\\
  \text{L. Marmet}& \text{(Riyadh)}& 8.4\times 10^{15} \text{ to } 8.6\times 10^{15}& \text{September $11^{th}$, 2000}\\
  \text{Z. McGregor-Dorsey}& \text{(Castalia)}& 8.2\times 10^{15} \text{ to } 8.4\times 10^{15}& \text{September $3^{rd}$, 2000}\\
  \text{Z. McGregor-Dorsey}& \text{(Abzug)}& 8.0\times 10^{15} \text{ to } 8.2\times 10^{15}& \text{October $16^{th}$, 2000}\\
  \text{Z. McGregor-Dorsey}& \text{(Hayduke)}& 7.8\times 10^{15} \text{ to } 8.0\times 10^{15}& \text{October $16^{th}$, 2000}\\
  \text{Z. McGregor-Dorsey}& \text{(Abzug)}& 7.6\times 10^{15} \text{ to } 7.8\times 10^{15}& \text{September $27^{th}$, 2000}\\
  \text{Z. McGregor-Dorsey}& \text{(Hayduke)}& 7.4\times 10^{15} \text{ to } 7.6\times 10^{15}& \text{August $26^{th}$, 2000}\\
  \text{N. Marmet}& \text{(Computer)}& 7.2\times 10^{15} \text{ to } 7.4\times 10^{15}& \text{September $17^{th}$, 2000}\\
  \text{L. Marmet}& \text{(Riyadh)}& 7.0\times 10^{15} \text{ to } 7.2\times 10^{15}& \text{August $25^{th}$, 2000}\\
  \text{G. Engebreth}& \text{(Computer)}& 6.8\times 10^{15} \text{ to } 7.0\times 10^{15}& \text{August $28^{th}$, 2000}\\
  \text{Z. McGregor-Dorsey}& \text{(Hayduke)}& 6.6\times 10^{15} \text{ to } 6.8\times 10^{15}& \text{August $17^{th}$, 2000}\\
  \text{Z. McGregor-Dorsey}& \text{(Abzug)}& 6.4\times 10^{15} \text{ to } 6.6\times 10^{15}& \text{September $2^{nd}$, 2000}\\
  \text{Z. McGregor-Dorsey}& \text{(Castalia)}& 6.2\times 10^{15} \text{ to } 6.4\times 10^{15}& \text{August $26^{th}$, 2000}\\
  \text{Z. McGregor-Dorsey}& \text{(Abzug)}& 6.0\times 10^{15} \text{ to } 6.2\times 10^{15}& \text{August $7^{th}$, 2000}\\
  \text{G. Engebreth}& \text{(Computer)}& 5.8\times 10^{15} \text{ to } 6.0\times 10^{15}& \text{August $10^{th}$, 2000}\\
  \text{L. Marmet}& \text{(Riyadh)}& 5.6\times 10^{15} \text{ to } 5.8\times 10^{15}& \text{August $6^{th}$, 2000}\\
  \text{Z. McGregor-Dorsey}& \text{(Abzug)}& 5.4\times 10^{15} \text{ to } 5.6\times 10^{15}& \text{August $26^{th}$, 2000}\\
  \text{Z. McGregor-Dorsey}& \text{(Hayduke)}& 5.2\times 10^{15} \text{ to } 5.4\times 10^{15}& \text{July $28^{th}$, 2000}\\
  \text{Z. McGregor-Dorsey}& \text{(Castalia)}& 5.0\times 10^{15} \text{ to } 5.2\times 10^{15}& \text{August $16^{th}$, 2000}\\
  \text{Z. McGregor-Dorsey}& \text{(Castalia)}& 4.8\times 10^{15} \text{ to } 5.0\times 10^{15}& \text{July $28^{th}$, 2000}\\
  \text{N. Marmet}& \text{(Computer)}& 4.6\times 10^{15} \text{ to } 4.8\times 10^{15}& \text{August $15^{th}$, 2000}\\
  \text{E. Wong}& \text{(Computer)}& 4.4\times 10^{15} \text{ to } 4.6\times 10^{15}& \text{September $27^{th}$, 2000}\\
  \text{Z. McGregor-Dorsey}& \text{(Hayduke)}& 4.2\times 10^{15} \text{ to } 4.4\times 10^{15}& \text{July $14^{th}$, 2000}\\
  \end{array}
\end{displaymath}

\begin{displaymath}
  \begin{array}{llrr}
  \text{User Name}& \text{(Computer Name)}& Range& Date Completed\\
\hline
  \text{L. Marmet}& \text{(Riyadh)}& 4.0\times 10^{15} \text{ to } 4.2\times 10^{15}& \text{July $20^{th}$, 2000}\\
  \text{Z. McGregor-Dorsey}& \text{(Abzug)}& 3.8\times 10^{15} \text{ to } 4.0\times 10^{15}& \text{July $18^{th}$, 2000}\\
  \text{Z. McGregor-Dorsey}& \text{(Castalia)}& 3.6\times 10^{15} \text{ to } 3.8\times 10^{15}& \text{July $10^{th}$, 2000}\\
  \text{Z. McGregor-Dorsey}& \text{(Abzug)}& 3.4\times 10^{15} \text{ to } 3.6\times 10^{15}& \text{July $10^{th}$, 2000}\\
  \text{Z. McGregor-Dorsey}& \text{(Hayduke)}& 3.2\times 10^{15} \text{ to } 3.4\times 10^{15}& \text{July $1^{st}$, 2000}\\
  \text{L. Marmet}& \text{(Riyadh)}& 3.0\times 10^{15} \text{ to } 3.2\times 10^{15}& \text{June $29^{th}$, 2000}\\
  \text{Z. McGregor-Dorsey}& \text{(Castalia)}& 2.8\times 10^{15} \text{ to } 3.0\times 10^{15}& \text{June $18^{th}$, 2000}\\
  \text{G. Engebreth}& \text{(Computer)}& 2.6\times 10^{15} \text{ to } 2.8\times 10^{15}& \text{July $22^{nd}$, 2000}\\
  \text{Z. McGregor-Dorsey}& \text{(Castalia)}& 2.4\times 10^{15} \text{ to } 2.6\times 10^{15}& \text{June $8^{th}$, 2000}\\
  \text{N. Marmet}& \text{(Computer)}& 2.2\times 10^{15} \text{ to } 2.4\times 10^{15}& \text{July $5^{th}$, 2000}\\
  \text{L. Marmet}& \text{(Riyadh)}& 2.0\times 10^{15} \text{ to } 2.2\times 10^{15}& \text{June $9^{th}$, 2000}\\
  \text{Z. McGregor-Dorsey}& \text{(Abzug)}& 1.8\times 10^{15} \text{ to } 2.0\times 10^{15}& \text{May $27^{th}$, 2000}\\
  \text{A. Simpson}& \text{(Computer)}& 1.6\times 10^{15} \text{ to } 1.8\times 10^{15}& \text{June $30^{th}$, 2000}\\
  \text{Z. McGregor-Dorsey}& \text{(Hayduke)}& 1.4\times 10^{15} \text{ to } 1.6\times 10^{15}& \text{May $13^{th}$, 2000}\\
  \text{L. Marmet}& \text{(Riyadh)}& 1.2\times 10^{15} \text{ to } 1.4\times 10^{15}& \text{May $16^{th}$, 2000}\\
  \text{G. Engebreth}& \text{(Computer)}& 1.1\times 10^{15} \text{ to } 1.2\times 10^{15}& \text{June $1^{st}$, 2000}\\
  \text{E. Wong}& \text{(Computer)}& 1.0\times 10^{15} \text{ to } 1.1\times 10^{15}& \text{May $31^{st}$, 2000}\\
  \text{L. Marmet}& \text{(Lion)}& 9 \times 10^{14} \text{ to } 10 \times 10^{14}& \text{June $13^{th}$, 2000}\\
  \text{L. Marmet}& \text{(Riyadh)}& 8 \times 10^{14} \text{ to } 9 \times 10^{14}& \text{April $5^{th}$, 2000}\\
  \text{Z. McGregor-Dorsey}& \text{(Castalia)}& 7 \times 10^{14} \text{ to } 8 \times 10^{14}& \text{May $1^{st}$, 2000}\\
  \text{L. Marmet}& \text{(Fontaine)}& 6 \times 10^{14} \text{ to } 7 \times 10^{14}& \text{April $22^{nd}$, 2000}\\
  \text{E. Wong}& \text{(Computer)}& 5 \times 10^{14} \text{ to } 6 \times 10^{14}& \text{April $8^{th}$, 2000}\\
  \text{L. Marmet}& \text{(Fontaine)}& 4 \times 10^{14} \text{ to } 5 \times 10^{14}& \text{March $5^{th}$, 2000}\\
  \text{D. Bernier}& \text{(Pentium 500MHz)}& 3 \times 10^{14} \text{ to } 4 \times 10^{14}& \text{March $24^{th}$, 2000}\\
  \text{L. Marmet}& \text{(Riyadh)}& 2 \times 10^{14} \text{ to } 3 \times 10^{14}& \text{February $4^{th}$, 2000}\\
  \text{D. Bernier}& \text{(Pentium 500MHz)}& 1.5\times 10^{14} \text{ to } 2.0\times 10^{14}& \text{January $24^{th}$, 2000}\\
  \text{L. Marmet}& \text{(Lion)}& 4 \times 10^{0} \text{ to } 1500 \times 10^{11}& \text{December $20^{th}$, 1999}
  \end{array}
\end{displaymath}

\end{document}